\documentclass[10pt,a4paper]{article}
\usepackage[utf8]{inputenc}
\usepackage[english]{babel}
\usepackage{amsmath}
\usepackage{amsthm}
\usepackage{amsfonts}
\usepackage{amssymb}
\usepackage{color}
\usepackage{graphicx}
\usepackage{hyperref}
\usepackage{ulem}
\usepackage{mathrsfs}
\usepackage{lineno}
\usepackage{verbatim}
\usepackage{subcaption}
\usepackage{cite}

\usepackage[margin=1.2in]{geometry}

\usepackage{todonotes}

\title{Numerical schemes for a multi-species BGK model with velocity-dependent collision frequency}
\author{J. Haack, C. Hauck, C. Klingenberg, M. Pirner, S. Warnecke}

\newcommand{\dv}{\,\mathrm{d}\mathbf{v}}
\newcommand{\dr}{\,\mathrm{d}\mathbf{x}}
\newcommand{\dx}{\mathrm{d}x}
\newcommand{\mbv}{\mathbf{v}}
\newcommand{\mbx}{\mathbf{x}}
\newcommand{\mba}{\mathbf{a}}
\newcommand{\mbu}{\mathbf{u}}

\newcommand{\mbq}{\mathbf{q}}
\newcommand{\mblambda}{\boldsymbol{\lambda}}
\newcommand{\mbalpha}{\boldsymbol{\alpha}}
\newcommand{\mbmu}{\boldsymbol{\mu}}
\newcommand{\mix}{{\rm mix}}

\newtheorem{theorem}{Theorem}[section]
\newtheorem{definition}{Definition}[section]

\newtheorem{proposition}{Proposition}[section]
\newtheorem{corollary}{Corollary}[section]
\newtheorem{assumption}{Assumption}[section]

\begin{document}
\maketitle


\section*{Abstract}
We consider a kinetic description of multi-species gas mixture modeled with Bhatnagar-Gross-Krook (BGK) collision operators, in which the collision frequency varies not only in time and space but also with the microscopic velocity. 
In this model, the Maxwellians typically used in standard BGK operators are replaced by a generalization of such target functions, which are defined by a variational procedure \cite{HHKPW}. 
In this paper we present a numerical method for simulating this model, which uses an Implicit-Explicit (IMEX) scheme to minimize a certain potential function, mimicking the Lagrange functional that appears in the theoretical derivation. 
We show that theoretical properties such as conservation of mass, total momentum and total energy as well as positivity of the distribution functions are preserved by the numerical method, and illustrate its usefulness and effectiveness with numerical examples.

\tableofcontents

\section{Introduction}

In a kinetic description, the state of a dilute gas or plasma is given by a distribution function that prescribes the density of particles at each point in position-momentum phase space.  In a time-dependent setting, the evolution of this distribution function is due to a balance of particle advection and binary collisions.  Perhaps the most well-known model for collisions is the Boltzmann collision operator, an integral operator that preserves collision invariants and dissipates the mathematical entropy of the system.  Unfortunately, the expense of evaluating this operator can be prohibitive. Indeed, its evaluation requires the calculation of a five-dimensional integral at every point in phase-space. Thus even with fast spectral methods \cite{mouhot2006fast,pareschi2000numerical,gamba2009spectral,gamba2017fast}, the collision operator is typically the dominant part of a kinetic calculation.
Furthermore, grid resolution requirements for multi-species Boltzmann collision operators provide additional constraints on the computational expense \cite{munafo2014}. 

The Bhatnagar-Gross-Krook (BGK) operator is a widely used surrogate for the Boltzmann operator that models collisions by a simple relaxation mechanism.  This simplification brings significant computational advantages while also maintaining the conservation and entropy dissipation properties of the Boltzmann operator.  However, the BGK operator does not recover all of the physics of the Boltzmann operator. Most notably, it cannot recover correct viscosity and heat conduction coefficients at the same time, although this problem can be remedied with a slightly extended model \cite{holway1965kinetic,andries2000gaussian}.  Another significant limitation is that the strength of relaxation in the standard BGK model is characterized by a collision frequency that is independent of particle velocity when, in reality, the collision frequency is expected to be velocity dependent \cite{struchtrup2005macroscopic, LeeMore}.

Velocity-dependent frequencies were first incorporated into single-species BGK models in \cite{Struchtrup97}, with a subsequent numerical implementation in \cite{MieussensStruchtrup2004}.  Unlike the constant-velocity case, the target of the relaxation model does not have the same mass, momentum, and energy density as the kinetic distribution, even though it has the form of a Maxwellian distribution, (i.e., the exponential of a polynomial).  An extension of the single-species model to the multi-species setting was recently developed in \cite{HHKPW}.  There, the existence and uniqueness of well-defined target functions for the relaxation operator were established rigorously via the solution of a convex entropy minimization problem.  Again these targets are in general not Maxwellians in the classical sense; instead they match the kinetic distributions via moments that are weighted by the collision frequency.

In this paper we present a numerical implementation of the velocity-dependent, multi-species BGK model developed in \cite{HHKPW}.  The implementation is a discrete velocity method that relies on standard spatial and temporal discretizations from the literature.  The key new ingredient is a solver which enables an implicit treatment of the BGK operator. As in the analytic case, the crucial step involves the formulation of a convex entropy minimization problem.  In particular, the solver uses a numerical minimization procedure in order to determine the coefficients of the target functions. 
This construction guarantees conservation and entropy properties at the discrete level, up to numerical tolerances, even when using a discrete velocity mesh.  In this sense, it is related to the implementation in \cite{mieussens2000discrete}, which considered a single-species BGK model with velocity-independent frequency.

The optimization problem that must be solved with velocity-dependent frequencies adds considerable expense to the simulation.  We solve it numerically by applying Newton's method to the dual problem.   While the overall computational cost depends strongly on the details of the implementation  \cite{schaerer2017efficient,AlldredgeHauckOLearyTits2014,AlldredgeHauckTits2012,GarrettHauckHill2015,Abramov2007}, the quadrature approximation of integrals in the gradient and Hessian of the dual objective is the most expensive part of the calculation.

In spite of the additional expense from the optimization problem, the number of operations needed to evaluate the BGK operator with a velocity-dependent collision frequency still scales like $\mathcal{O}(N^3)$, where $N$ is the number of points in each dimension of the velocity grid. 
In comparison, the fastest algorithms for evaluating the Boltzmann collision operator are spectral methods, whose complexity for general collision kernels scales like $\mathcal{O}(MN^4\log N)$ \cite{mouhot2006fast} but for specialized kernels can be reduced to $\mathcal{O}(MN^3\log N)$ \cite{gamba2017fast}. 
Here $M$ is the number of quadrature points for approximating the integrals over the unit sphere $\mathbb{S}^2$ in $\mathbb{R}^3$. 
The size of $M$ is problem dependent, but typically $N\leq M \ll N^2$ \cite{mouhot2006fast}.  
Thus, while more expensive than the standard BGK models, the BGK model with velocity-dependent collision frequencies is still of lower computational complexity than the Boltzmann collision operator. 
Additionally, each species in the BGK model can be discretized on a separate velocity grid, while grid discretization in multi-species Boltzmann models introduces expensive grid resolution requirements for problems with significant differences in species masses.

The remainder of this paper is organized as follows. 
In Section \ref{sec:model}, we recall the multi-species BGK model from \cite{HHKPW} with velocity-dependent collision frequency. 
In Section \ref{sec:time}, we present the first- and second-order implicit-explicit time discretizations that are used in the paper.  We also introduce the optimization-based approach for the implicit evaluation of the BGK operator.  In Section \ref{sec:space}, we describe the space discretization.
In Section \ref{sec:discreteproperties}, we verify some structure preserving properties of the semi-discrete scheme.  In Section \ref{sec:velocity}, we introduce the velocity discretization and summarize the numerical implementation of the optimization algorithm introduced in Section \ref{sec:time}.
In Section \ref{sec:results}, we provide an array of numerical results that illustrate the properties of our scheme and explore the effects of velocity-dependent collision frequencies that are motivated by Coulomb interactions common to plasmas. 
Discussion and conclusions are provided in Section \ref{sec:conclusions}.

\section{A consistent multi-species BGK model with velocity-dependent collision frequency}
\label{sec:model}

For simplicity we focus on a mixture of two species; because collisions are assumed to be binary, the generalization to more species is straightforward.
Given the spatial coordinate $\mbx \in  \mathbb{R}^3$, velocity coordinate $\mbv \in \mathbb{R}^3$, and time $t \geq 0$,
we consider two scalar-valued functions $f_1 = f_1(\mbx,\mbv,t) \geq 0$ and $f_2 = f_2(\mbx,\mbv,t) \geq 0$ which give the phase space density (i.e. the density with respect to the measure $\dr \dv$) of species with masses $m_1$ and $m_2$, respectively.

In this setting, the BGK model in \cite{HHKPW} takes the form
\begin{align} \label{eq:BGK}
\begin{split}
\partial_t f_{1} + \mbv\cdot \nabla_\mbx f_1  &= \nu_{11}  (A_{11} -f_1) + \nu_{12}  (A_{12} - f_{1}),\\
\partial_t f_2+ \mbv\cdot \nabla_\mbx f_2  &= \nu_{22}  (A_{22} -f_2) + \nu_{21} (A_{21} - f_{2}),
\end{split}
\end{align}
where  $ \nu_{ij} = \nu_{ij}(\mbx,\mbv,t) \geq 0$  for $i,j=1,2$ are velocity-dependent collision frequencies. 
The intra-species target functions $A_{11}$ and $A_{22}$ take the form
\begin{align} \label{eq:shape_1}
A_{ii} = e^{\mblambda_i \cdot \mba_i(\mbv)},
\end{align}
where $\mba_i( \mbv) = m_i(1,\mbv,\vert \mbv   \vert^2)^\top$, $\mblambda_i = (\lambda_i^0, \mblambda_i^1, \lambda_i^2)^\top  \in\mathbb{R}\times \mathbb{R}^3 \times \mathbb{R}^-$, while the inter-species target functions $A_{12}$ and $A_{21}$ take the form
\begin{align} \label{eq:shape_2}
A_{ij} = e^{\mblambda_{ij} \cdot \mba_i(\mbv)},
\end{align}
where $\mblambda_{ij} = (\lambda_{ij}^0, \mblambda^1, \lambda^2)^\top \in \mathbb{R}\times \mathbb{R}^3 \times \mathbb{R}^-$.
For convenience, we suppress the dependence of $A_{ij}$ and the parameters $\mblambda_1,\mblambda_2, \mblambda_{12}$ and $\mblambda_{21}$ on $f_1$ and $f_2$.  However, these parameters are directly tied to $f_1$ and $f_2$.  Indeed they are the Lagrange multipliers associated with a minimization procedure with constraints involving $f_1$ and $f_2$ that enforce requisite conservation laws.

\begin{definition}
For $i \in \{1,2\}$, the species mass, momentum, and energy densities are given by%
\footnote{Here we suppress the Boltzmann constant $k_B$ in the definition of the temperature for ease of presentation; however, in some later formulas we include it for emphasis. }
\begin{equation}
\label{eq:conserved}
    \rho_i = \int m_i f_i \dv, 
    \quad
    \rho_i \mbu_i = \int m_i \mbv f_i \dv, \quad \text{and}\,
    \quad 
    \frac{1}{2}\rho_i |\mbu_i|^2 + \frac{3}{2} n_i T_i 
    = \frac12 \int  m_i |\mbv|^2 f_i \dv,
\end{equation} 
respectively, where the species number density $n_i$, mean velocity $\mbu_i$ and species mean temperature $T_i$ are given by
\begin{equation}
\label{eq:u_T}
    n_i = \frac{\rho_i}{m_i}
    \qquad \qquad 
    \mbu_i = \frac{\int \mbv f_i \dv }{\int f_i \dv}
    \qquad \text{and} \qquad 
    T_i = \frac{1}{3}\frac{\int m_i |\mbv-\mbu_i|^2 f_i \dv}{\int f_i \dv },
\end{equation}
respectively.   The total momentum and total energy are the sums of the individual species momenta and energy.
\end{definition}

\begin{definition}
For a given number density $n >0$, mean velocity $\mbu \in \mathbb{R}^3$, and temperature $T>0$, {a Maxwellian distribution} for a species with mass $m_i$ is given by
\begin{align} \label{eq:Maxwellian}
    M_i[n,\mbu,T](\mbv) = n \left( \frac{m_i}{2 \pi T}\right)^{3/2} \exp\left(-\frac{m_i|\mbv-\mbu|^2}{2T}\right).
\end{align}
If a distribution function $f_i$ has number density $n_i$, mean velocity $\mbu_i $, and temperature $T_i$,  then we call $M_i[n_i,\mbu_i, T_i]$ the Maxwellian of $f_i$.
\end{definition}
The target functions in \eqref{eq:shape_1} are chosen so that the individual mass, momentum, and energy densities are invariant under intra-species collisions; that is, for $i \in \{1,2 \}$,
\begin{align}
\label{eq:moment_intra}
    \int \nu_{ii} \mba_i (A_{ii} - f_i) \dv = 0.
\end{align}
When $\nu_{ii}$ is independent of velocity, these constraints recover the standard BGK model. In particular,  the coefficients $\mblambda_i$ used in the functional form \eqref{eq:shape_1} are unique and can be found analytically.  
On the other hand, when $\nu_{ii}$ is a function of $\mbv$, the integrals in \eqref{eq:moment_intra} cannot, in general, be evaluated analytically.  
However the target function can still be realized as the solution of the weighted entropy minimization problem
\begin{align}
\min_{g \in \chi_i} \int \nu_{ii}(\mbv) h(g(\mbv)) \dv, \quad i \in\{1,2\},
\label{min}
\end{align}
where 
\begin{equation}
\label{eq:h}
    h(z) = z \log(z) - z , \quad z >0
\end{equation}
and the constraint set 
\begin{align}
\chi_i= \left\lbrace g ~ \Big |~ g > 0,\, \nu_{ii} (1+ |\mbv|^2) g \in L^1(\mathbb{R}^3),\, \int \nu_{ii} \mba_i
(g - f_i) \dv = 0 \right\rbrace
\end{align}
enforces the constraints in \eqref{eq:moment_intra}. 
When $\nu_{ii}$ is velocity-independent, \eqref{min} recovers the standard Maxwellian associated to $f_i$.
Moreover, if $\Lambda = \{\mbalpha = (\alpha^0, \mbalpha^1, \alpha^2)^\top \in\mathbb{R}\times \mathbb{R}^3 \times \mathbb{R}^- \}$, then the multipliers $\mblambda_i = (\lambda_i^0, \mblambda_i^1, \lambda_i^2)^\top  \in \mathbb{R}\times \mathbb{R}^3 \times \mathbb{R}^-$ solve the dual of \eqref{min}:
\begin{align}\label{eq:dual_1}
  \mblambda_i = \operatorname*{argmin}_{\mbalpha \in \Lambda } 
  \left\{-\int \nu_{ii} e^{\mbalpha \cdot \mba_i(\mbv)} \dv + \mbalpha \cdot  \int \nu_{ii} \mba_i( \mbv) f_i  \dv \right\}.
\end{align}
For more details see \cite{HHKPW}. 

The inter-species target functions \eqref{eq:shape_2} are constrained by the conservation of species mass density and the total momentum and energy, i.e., 
\begin{alignat}{3}
\label{eq:moment_inter}
    \int \nu_{1} m_1 (A_{12} - f_1) \dv = 0,&\quad
    \quad 
    \int \nu_{2} m_2 (A_{21} - f_2) \dv = 0,& \nonumber\\
    \int \nu_{12} m_1 \mbv  (A_{12} - f_1) \dv &+ 
    \int \nu_{21} m_2 \mbv  (A_{21} - f_2) \dv = 0,&\\
    \int \nu_{12} m_1 |\mbv|^2  (A_{12} - f_1) \dv &+ 
    \int \nu_{21} m_2 |\mbv|^2  (A_{21} - f_2) \dv = 0.& \nonumber
\end{alignat}
Unlike the intra-species case, these constraints do not uniquely identify the multipliers, and in the case of constant collision frequencies, several works \cite{HHM,bobylev2018general,klingenberg2017consistent} have explored approaches that use the remaining degrees of freedom to satisfy additional properties and/or match transport coefficients.  
For velocity-dependent collision frequencies, the constraints in \eqref{eq:moment_inter} cannot, in general, be computed analytically.  
However, similar to the intra-species setting, the target functions in the inter-species setting can be formulated as the solution of the weighted entropy minimization problem
\begin{align}
\label{eq:entropy_min_cross_terms}
\min_{g_{1}, g_{2} \in \chi_{12}}  \int \nu_{12} h(g_1) \dv + \int \nu_{21} h(g_2) \dv , 
\end{align}
where $h$ is defined in \eqref{eq:h} and
\begin{align}
\label{eq:chi12_def}
\begin{split}
\chi_{12}= \Bigg\lbrace (g_{1}, g_{2})~\Big|& ~ g_{1}, g_{2} >0,\, 
\nu_{12} (1+ |\mbv|^2) g_{1},\,
\nu_{21} (1+ |\mbv|^2) g_{2}  \in L^1(\mathbb{R}^3),\\ 
&\int m_1 \nu_{12}  (g_{1} - f_1)\dv =0, \quad 
\int  m_2 \nu_{21} (g_{2} - f_2)\dv =0,  \\
&\int m_1 \nu_{12} 
\begin{pmatrix}
 \mbv \\  |\mbv|^2
\end{pmatrix} (g_{1} - f_1)\dv  + 
\int m_2 \nu_{21} \begin{pmatrix}
 \mbv \\ |\mbv|^2
\end{pmatrix} (g_{2} - f_2) \dv = 0 \Bigg\rbrace.
\end{split}
\end{align}
The solution to \eqref{eq:entropy_min_cross_terms} is the pair of target functions given in \eqref{eq:shape_2}.  Moreover, if $\Lambda_{12} = \{(\mbalpha_{12},\mbalpha_{21}): \mbalpha_{12} = (\alpha^0_{12}, \mbalpha^1, \alpha^2)^\top \in \mathbb{R}\times \mathbb{R}^3 \times \mathbb{R}^-, \mbalpha_{21} = (\alpha^0_{21}, \mbalpha^1, \alpha^2)^\top \in\mathbb{R}\times \mathbb{R}^3 \times \mathbb{R}^-\}$, then the multipliers $\mblambda_{12},\mblambda_{21}$ satisfy the dual problem
\begin{align}\label{eq:dual_2}
\begin{split}
    (\mblambda_{12},\mblambda_{21})= \operatorname*{argmin}_{(\mbalpha_{12},\mbalpha_{21}) \in \Lambda_{12}} \Bigg\{ &-\int  (\nu_{12}e^{\mbalpha_{12} \cdot \mba_{12}(\mbv)} + \nu_{21}e^{\mbalpha_{21} \cdot \mba_{21}(\mbv)})\dv \\
&+\alpha_{12}^0 \int m_1 \nu_{12} f_1\dv + \alpha_{21}^0 \int m_2 \nu_{21} f_2\dv \\
& + \mbalpha^1 \cdot  \int    \mbv (m_1 \nu_{12} f_1+m_2 \nu_{21} f_2) \dv \\
&+ \alpha^2  \int    |\mbv|^2 (m_1 \nu_{12} f_1+m_2 \nu_{21} f_2) \dv \Bigg\}.
\end{split}
\end{align}

The existence and uniqueness of solutions $\mblambda_1 \in \Lambda$, $\mblambda_2 \in \Lambda$, and $(\mblambda_{12},\mblambda_{21}) \in \Lambda_{12}$ to the dual problems in \eqref{eq:dual_1} and \eqref{eq:dual_2} are proven in \cite{HHKPW}.   As a consequence, solutions to \eqref{eq:BGK} satisfy the appropriate conservation laws, dissipate the total entropy density 
\begin{equation}
    H(f_1,f_2) = \int [h(f_1)+h(f_2)]\dv,
\end{equation}
and verify an H-Theorem. 
Specifically, we have the following
\begin{theorem}[\cite{HHKPW}]
Let $f_1 \geq 0 $ and $f_2 \geq 0 $, with neither identically zero, solve \eqref{eq:BGK} with target functions defined by \eqref{min} and \eqref{eq:entropy_min_cross_terms}.
Then the following conservation laws hold:
\begin{subequations}\label{eq:cons_macros}
\begin{gather} 
    \partial_t \rho_1 + \nabla_x \cdot  \int m_1 \mbv f_1 \dv = 0, \\
    \partial_t \rho_2 + \nabla_x \cdot  \int m_2 \mbv f_2 \dv = 0, \\
    \partial_t  (\rho_1 \mbu_1  + \rho_2 \mbu_2)
    + \nabla_x \cdot \left(
    \int \mbv \otimes \mbv  (m_1 f_1 + m_2 f_2) \dv 
    \right) = 0,\\
     \partial_t  \left(\frac{\rho_1 |\mbu_1|^2}{2}  + \frac{3\rho_1 T_1}{2m_1} + 
     \frac{\rho_2 |\mbu_2|^2}{2}  + \frac{3\rho_2 T_2}{2m_2} \right) 
    + \nabla_x \cdot \left(
    \int  \frac12 \mbv |\mbv|^2  (m_1 f_1 + m_2 f_2)\dv 
    \right) =0,
\end{gather} 
\end{subequations}
where $\rho_i$, $\mbu_i$, and $T_i$ are defined in \eqref{eq:conserved} and \eqref{eq:u_T}.  Moreover, 
\begin{align}
    \partial_t H(f_1,f_2) + \nabla_\mbx \cdot \left(\int \mbv [h(f_1) + h(f_2)]\dv \right) \leq 0,
\end{align}
with equality if and only if $f_1$ and $f_2$ are two Maxwellian distributions with the same mean velocities $\mbu_{\rm eq}(\mbx)$ and temperatures $T_{\rm eq}(\mbx)$.
\end{theorem}

Before moving to the numerical implementation of the velocity-dependent BGK model, we present some additional definitions and a key assumption on the collision frequencies. 

\begin{definition}
The mixture mean velocity $\mbu_\mix$ and the mixture temperature $T_\mix$ are given by
\begin{align} \label{eq:u_mixture}
 \mbu_\mix = \frac{\rho_1 \mbu_1+ \rho_2 \mbu_2}{\rho_1+\rho_2 }
\end{align}  
and
\begin{subequations}
 \label{eq:T_mixture}
\begin{align}
 \label{eq:T_mixture_a}
    T_\mix &= \frac{n_1 T_1 + n_2 T_2}{n_1+n_2} + \frac{\rho_1 (|\mbu_1|^2 - |\mbu_\mix|^2) + \rho_2(|\mbu_2|^2 - |\mbu_\mix|^2)}{{3(n_1 + n_2)}}  \\
 \label{eq:T_mixture_b}
    &= \frac{n_1 T_1 + n_2 T_2}{n_1+n_2} + \frac{1}{3}\frac{\rho_1\rho_2}{\rho_1+\rho_2}  \frac{|\mbu_1-\mbu_2|^2}{n_1 + n_2}.
\end{align} 
\end{subequations}
\end{definition}

\begin{proposition} \label{prop:mixture_quantities}
In the spatially homogeneous setting, $\mbu_\mix$ and $T_\mix$ are constant in time.

\begin{proof}
In the homogeneous setting, $\rho_1$, $\rho_2$, and $\rho_1 \mbu_1+ \rho_2 \mbu_2$ are all constant in time.  Hence the ratio in \eqref{eq:u_mixture} that defines $\mbu_\mix$ is also constant in time.  To show that $T_\mix$ is constant in time, we use \eqref{eq:T_mixture_a} to write
\begin{equation}
    \label{eq:T_mix_form}
    \frac{3}{2}(n_1 + n_2) T_\mix = \rm{I} - \rm{II},
\end{equation}
where
\begin{align}
{\rm{I}} = \frac{1}{2} \rho_1 |\mbu_1|^2 + \frac{3}{2} n_1 T_1 
        + \frac{1}{2} \rho_2 |\mbu_2|^2 + \frac{3}{2} n_2 T_2
\end{align}
is the total energy and 
\begin{align} 
{\rm{II}} 
    =  \frac{1}{2} (\rho_1 + \rho_2 ) |\mbu_\mix|^2 
    =  \frac{1}{2} \frac{(\rho_1 \mbu_1+ \rho_2 \mbu_2)^2 }{\rho_1 + \rho_2}  .
\end{align}
In the homogeneous setting, both I and II are constant in time, as are $n_1$ and $n_2$.  Thus the formula in \eqref{eq:T_mix_form} implies $T_\mix$ is also constant in time.
\end{proof}
\end{proposition}

For the remainder of the paper we make the following assumption on the collision frequencies $\nu_{ij}$:
\begin{assumption} \label{assumption:nu}
The space and time dependency of the collision frequencies $\nu_{ij}$ arises only via dependence on the mass densities $\rho_i$, the mixture mean velocity $\mbu_\mix$, and the mixture temperature $T_\mix$. Because the collision operators in \eqref{eq:BGK} conserve these quantities, the collision frequencies $\nu_{ij}$ are independent of time in the space homogeneous setting.
\end{assumption}

This assumption is common for standard collision rates in the literature and follows from cross section definitions; see for example \cite{KrallTriv, HHM}. 
In this paper, we add a dependence on the microscopic velocity $\mbv$ to $\nu_{ij}$, which arises naturally from the derivation of the BGK operator from the Boltzmann equation \cite{Struchtrup97, HHM}. 
This dependence is neglected for computational convenience in the standard BGK model, and may have profound effects on the relaxation process and the resulting hydrodynamic behavior.  
In particular, transport coefficients derived via the Chapman-Enskog expansion (e.g. thermal conductivity) are sensitive to the dynamics of the tails of the kinetic distribution \cite{CC}.

\section{Time discretization} 
\label{sec:time}

Let $i,j=1,2$ and $i\neq j$. We write \eqref{eq:BGK} as 
\begin{align}
    \partial_t f_i + \mathcal{T}(f_i) = \mathcal{R}_{i}(f_i,f_j)
\end{align}
with the relaxation operator  
\begin{align}
\label{eq:R-def}
\mathcal{R}_{i}(f_i,f_j) = \nu_{ii} \left(A_{ii}-f_i\right) + \nu_{ij} \left(A_{ij}-f_i\right)
\end{align}
and the transport operator
\begin{align}
\mathcal{T}(f_i) = \mbv\cdot \nabla_\mbx f_i. 
\end{align}
When the collision frequencies are large, the operator $\mathcal{R}_{i}$ becomes stiff and an implicit treatment is preferred.  With this fact in mind, we pursue implicit-explicit (IMEX) schemes where $\mathcal{T}$ is treated explicitly and $\mathcal{R}_{i}$ is treated implicitly.  

When the collision frequencies are constant, the inversion of $\mathcal{R}_{i}$ is not difficult. In the single-species case, the inversion is trivial because the target function does not evolve during the evolution of the space homogeneous system \cite{coron1991numerical,pieraccini2007implicit}.  Thus the inversion of   $\mathcal{R}_{i}$ can be reduced to a linear solve.  In the multi-species case, the target functions of the space homogeneous system do evolve, but they can be computed with an iterative solver for the mean velocities and temperatures. 

For velocity-dependent cross-sections, the problem is much more delicate.  This is because the target function parameters are no longer related to the moments in an analytical way, even in the one-species case.  One approach that is presented in \cite{mieussens2000discrete} and \cite{MieussensStruchtrup2004} is to linearize the attractor around the ansatz at the current value, in order to handle difficulties with its evaluation at the next time step.  The result is an efficient scheme for simulating steady-state solutions, but as noted in \cite{mieussens2000discrete}, this approach lacks conservation and entropy properties at the discrete level.

The schemes presented below preserve conservation properties, and the first-order version inherits additional desirable properties from the continuum model.  These properties are enforced by evaluating target functions at the next time step exactly (up to numerical tolerances) using a minimization procedure that mimics the theoretical formulations in \eqref{min} and \eqref{eq:entropy_min_cross_terms}.  The approach works for multi-species BGK equations equipped with a broad class of collision frequencies.   However, it does rely on Assumption \ref{assumption:nu}.  For example,
 given $t_\ell = \ell \Delta t$ for $\ell\in \mathbb{N}_0$ a simple update of $f_i^\ell \approx f_i(\mbx,\mbv,t_\ell)$ from $t_{\ell}$ to $t_{\ell+1}$ uses the approximation 
\begin{align}
\mathcal{R}_{i}(f_i^{\ell+1},f_j^{\ell+1}) \approx \nu_{ii}^\ell \left(A_{ii}^{\ell+1}-f_i^{\ell+1}\right) + \nu_{ij}^\ell \left(A_{ij}^{\ell+1}-f_i^{\ell+1}\right),
\end{align}
where $A_{ii}^{\ell+1}$ and $A_{ij}^{\ell+1}$ are discrete target functions that, as described in Section \ref{sec:generaltime}, depend on $f_i^{\ell+1}, f_j^{\ell+1}, \nu_{ii}^\ell$ and $\nu_{ij}^\ell$ via the solution of a convex minimization problem. 
Under Assumption \ref{assumption:nu}, the collision frequencies depend on quantities that are unchanged by the collisional process.  Thus their evaluation at time step $t_\ell$ is justified, since
\begin{align}\label{eq:nu_const}
    \nu_{ij}^{\ell+1} = \nu_{ij}(\rho_i^{\ell+1}, \rho_j^{\ell+1}, \mbu_\mix^{\ell+1}, T_\mix^{\ell+1}) = \nu_{ij}(\rho_i^{\ell}, \rho_j^{\ell}, \mbu_\mix^{\ell}, T_\mix^{\ell}) = \nu_{ij}^\ell.
\end{align}
However, in more general settings, lagging the collision frequencies in this way may cause a drop in temporal order for an otherwise high-order scheme \cite{LowrieLagging}.

\subsection{First-order splitting} \label{subsec:firstordersplit}
We split  \eqref{eq:BGK} into a relaxation step and the transport step. 
\paragraph{Relaxation} We execute the relaxation step in each spatial cell using a backward Euler method
\begin{align} \label{eq:relax}
\frac{f_i^{\ell'} -f_i^\ell}{\Delta t} = \mathcal{R}_{i}(f_i^{\ell'},f_j^{\ell'}),
\end{align}
which can be rewritten to express $f_i^{\ell'}$ as the convex combination
\begin{align} \label{eq:update_split}
    f_i^{\ell'} = c_i^\ell f_i^{\ell} + c_i^\ell \Delta t (\nu_{ii}^{\ell} A_{ii}^{\ell'} + \nu_{ij}^\ell A_{ij}^{\ell'})
\end{align}
with 
\begin{align} 
    c_i^\ell = \frac{1}{1+\Delta t  (\nu_{ii}^{\ell} + \nu_{ij}^{\ell})}.
\end{align}
If $A_{ii}^{\ell'}$ and $A_{ij}^{\ell'}$ can be expressed as functions of $f_i^\ell$, then \eqref{eq:update_split} provides an explicit update formula for $f_i^{\ell'}$. In Section \ref{sec:generaltime} we show how to determine $A_{ii}^{\ell'}$ and $A_{ij}^{\ell'}$ while preserving the conservation properties \eqref{eq:moment_intra} and \eqref{eq:moment_inter} at the discrete level.
\paragraph{Transport} We solve the transport in $x$ for $f_i^{\ell+1}$ by a forward Euler method with initial data $f_i^{\ell'}$:
\begin{align} \label{eq:transport_x}
\frac{f_i^{\ell+1} -f_i^{\ell'}}{\Delta t}+ \mathcal{T}(f_i^{\ell'}) = 0.
\end{align}
Details on the numerical approximation of $\mathcal{T}$ are given in Section \ref{sec:space}.

\subsection{Second-order IMEX Runge-Kutta} \label{subsec:secondorderIMEX}
For a second-order method, we use the following IMEX Butcher tableaux \cite{ARS97} \begin{align}
\renewcommand\arraystretch{1.2}
\begin{array}
{c|ccc}
0\\
\gamma &0& \gamma \\
1 &0&1-\gamma&\gamma \\
\hline
& 0&1-\gamma&\gamma
\end{array}
\hspace*{2cm}
\renewcommand\arraystretch{1.2}
\begin{array}
{c|ccc}
0\\
\gamma & \gamma \\
1 &\delta &1-\delta&0 \\
\hline
& \delta &1-\delta&0
\end{array}
\end{align}
with 
\begin{equation}
    \gamma = 1- \frac{\sqrt{2}}{2} 
    \quad \text{and} \quad
    \delta = 1-\frac{1}{2 \gamma}.
\end{equation}
The left table is used for the relaxation step, and the right table is used for the transport step. 
This IMEX Runge-Kutta scheme is L-stable and globally stiffly accurate (GSA).\footnote{The GSA property means that the numerical solution $f_i^{\ell+1}$ coincides with the last stage value of the method, which is important for accuracy of the method when the collision frequencies become large.}

Applying the method to \eqref{eq:BGK} results in the following scheme:  
\begin{subequations} \label{eq:update_imex}
\begin{align} 
    f_i^{(1)} &= f_i^\ell - \gamma \Delta t \,\mathcal{T}(f_i^\ell) + \gamma \Delta t\, \mathcal{R}_{i}(f_i^{(1)},f_j^{(1)}), 
    \label{eq:update_imex_1}\\
    f_i^{(2)} &= f_i^\ell - \delta \Delta t \,\mathcal{T}(f_i^\ell) - (1-\delta)\Delta t\, \mathcal{T}(f_i^{(1)}) \nonumber \\
    &\qquad \qquad \qquad+ (1-\gamma) \Delta t\,  \mathcal{R}_{i}(f_i^{(1)},f_j^{(1)}) + \gamma \Delta t \mathcal{R}_{i}(f_i^{(2)},f_j^{(2)}),
    \label{eq:update_imex_2} \\
    f_i^{\ell+1} &= f_i^{(2)}. 
\end{align}
\end{subequations}
Using the constants
\begin{align} 
    c_i^{(r)} = \frac{1}{1+\gamma \Delta t  (\nu_{ii}^{(r)} + \nu_{ij}^{(r)})},
\end{align}
we can rewrite \eqref{eq:update_imex_1} and \eqref{eq:update_imex_2} as convex combination of three terms
\begin{subequations}  \label{eq:update_ARS}
\begin{alignat}{3}
f_i^{(1)} &= c_i^{(1)}  G_i^{(1)}  &&+ c_i^{(1)}  \gamma \Delta t \,\nu_{ii}^{(1)} A_{ii}^{(1)} &&+  c_i^{(1)}  \gamma \Delta t \,\nu_{ij}^{(1)} A_{ij}^{(1)}  \\
f_i^{(2)} &= c_i^{(2)} G_i^{(2)} &&+ c_i^{(2)} \gamma \Delta t \, \nu_{ii}^{(2)} A_{ii}^{(2)} &&+ c_i^{(2)} \gamma \Delta t \, \nu_{ij}^{(2)} A_{ij}^{(2)} ,
\end{alignat}
\end{subequations}
where the quantities
\begin{subequations} 
\label{eq:Gs}
\begin{align}
G_i^{(1)}  &= f_i^\ell - \Delta t \, \gamma\, \mathcal{T}(f_i^\ell) \\
G_i^{(2)} &=  f_i^\ell - \Delta t \, \delta\, \mathcal{T}(f_i^\ell)-\Delta t \,(1-\delta) \mathcal{T}(f_i^{(1)}) + \Delta t \,(1-\gamma) \mathcal{R}_{i}(f_i^{(1)},f_j^{(1)})
\end{align}
\end{subequations}
depend on known data. The collision frequencies $\nu_{ii}^{(r)}$, $\nu_{ij}^{(r)}$ and constants $c_i^{(r)}$ are evaluated at the intermediate steps $G_i^{(r)}$.  This option maintains second-order accuracy as long as Assumption \ref{assumption:nu} applies.

The main computational challenge in each stage of \eqref{eq:update_ARS} is determining the parameters of the target functions. In the following section, we explain how to manage this.

\subsection{General implicit solver} \label{sec:generaltime}
We write the implicit updates in \eqref{eq:update_split} and \eqref{eq:update_ARS} above in a generic steady state form
\begin{align} \label{eq:update_general}
\psi_i  = c_i G_i  + c_i \gamma\Delta t  (\nu_{ii} B_{ii}  + \nu_{ij} B_{ij} ) 
\end{align}
where $B_{ii}$ and $B_{ij}$ are the unique target functions associated to $\psi_i$, 
\begin{equation}
    c_i = \frac{1}{1+\gamma\Delta t  (\nu_{ii} + \nu_{ij} )},
\end{equation}
and $G_i$ is a known function.  
The goal now is to express $B_{ii} $ and $B_{ij} $ as functions of $G_i$ and $G_j$ so that \eqref{eq:update_general} provides an explicit update formula for $\psi_i$.  Applying the conservation properties \eqref{eq:moment_intra} and \eqref{eq:moment_inter} to \eqref{eq:update_general} gives
\begin{align}
\begin{split}
\int \nu_{11}  B_{11}  \,\mba_1(\mbv) \dv +
\int \nu_{22}  B_{22}  \,\mba_2(\mbv) \dv +
\int \nu_{12}  B_{12}  \,\mba_1(\mbv) \dv + 
\int \nu_{21}  B_{21}  \,\mba_2(\mbv)  \dv \end{split} \nonumber \\ 
\begin{split} \overset{\eqref{eq:moment_intra},\eqref{eq:moment_inter}}{=}
\int \nu_{11}  \psi_1  \,\mba_1(\mbv) \dv +
\int \nu_{22}  \psi_2  \,\mba_2(\mbv) \dv +
\int \nu_{12}  \psi_1  \,\mba_1(\mbv) \dv +
\int \nu_{21}  \psi_2  \,\mba_2(\mbv) \dv 
\end{split} \nonumber \\ 
\begin{split}
&\hspace{0.25cm}\overset{(\ref{eq:update_general})}{=}
 \int  \nu_{11}   c_1  \left[ G_1  +  \Delta t \, \gamma \nu_{11}  B_{11}  + \Delta t \, \gamma \nu_{12}  B_{12}  \right] \,\mba_1(\mbv) \dv  \\ 
&\hspace*{1.5cm}
+ \int \nu_{22}  c_2  \left[ G_2  + \Delta t  \, \gamma   \nu_{22}  B_{22}  + \Delta t \, \gamma \nu_{21}  B_{21}  \right]   \,\mba_2(\mbv) \dv \\
&\hspace*{1.5cm}
+
\int \nu_{12}  c_1  \left[G_1  + \Delta t \, \gamma \nu_{11}  B_{11}  + \Delta t  \, \gamma \nu_{12}  B_{12} \right] \mba_1(\mbv) \dv \\ 
&\hspace*{1.5cm}
+ \int \nu_{21}  c_2  \left[G_2  + \Delta t \, \gamma \nu_{22}  B_{22}  + \Delta t \, \gamma \nu_{21}  B_{21}  \right]\mba_2(\mbv)   \dv. \end{split}
\end{align}
After sorting terms, we arrive at the following moment equations
\begin{align} \label{eq:moment_A_general}
\begin{split}
\int c_1  \left(\nu_{11} B_{11}+\nu_{12} B_{12}\right) \mba_1(\mbv)  \dv
+  \int c_2   \left(\nu_{21} B_{21}+\nu_{22} B_{22}\right) \mba_2(\mbv) \dv \\ =
\int  c_1  \left(\nu_{11} +\nu_{12} \right) G_1  \mba_1(\mbv)
+ \int  c_2  \left(\nu_{22} +\nu_{21} \right) G_{2}  \mba_2(\mbv) \dv,
\end{split}
\end{align}
which provide a set of constraints to determine  $B_{ii}$ and $B_{ij}$ from the given data $G_i$ and $G_j$. 

The constraints in \eqref{eq:moment_A_general} represent first-order optimality conditions associated to the minimization of the convex function 
\begin{align}\label{eq:potentialfunction}
   \varphi_{\rm tot}(\mbalpha_1,\mbalpha_2,\mbalpha) 
   =
-\int \left[ c_1 \nu_{11} B_{11} +  c_2 \nu_{22} B_{22} +  c_1 \nu_{12} B_{12} + c_2 \nu_{21} B_{21} \right] \dv 
+ \mbmu_1 \cdot \mbalpha_1 + \mbmu_2 \cdot \mbalpha_2 + \mbmu \cdot \mbalpha 
\end{align} 
where $\mbalpha_i = (\alpha_i^0,\mbalpha_i^1,\alpha_i^2)^\top \in\mathbb{R}\times \mathbb{R}^3 \times \mathbb{R}^-$;
\begin{align} \label{eq:input_mu_i}
\mbmu_i = 
\begin{pmatrix}
\mu_{i}^0 \\ \mbmu_{i}^1 \\ \mu_{i}^2
\end{pmatrix} =
\int c_i \nu_{ii}  G_i \mba_i(\mbv) \dv
\end{align}
for $i=1,2$; $\mbalpha = (\alpha_{12}^0,\alpha_{21}^0,\mbalpha^1,\alpha^2)^\top \in \mathbb{R}\times \mathbb{R}\times \mathbb{R}^3 \times \mathbb{R}^-$; and
\begin{align} \label{eq:input_mu}
\mbmu = 
\begin{pmatrix}
\mu_{12}^0 \\ \mu_{21}^0 \\ \mbmu^1 \\ \mu^2
\end{pmatrix} &= \int \left[ 
\begin{pmatrix}
1 \\ 0 \\ \mbv \\ \vert \mbv \vert^2
\end{pmatrix}
 m_1 c_1 \nu_{12}  G_1 +
\begin{pmatrix}
0 \\ 1 \\ \mbv \\ \vert \mbv \vert^2
\end{pmatrix}
m_2 c_2 \nu_{21}  G_2  \right] \dv.
\end{align}
The minimization problem can be decoupled as follows:
\begin{proposition}
The components of the minimizer of (\ref{eq:potentialfunction}) can be found by minimizing the following three convex functions independently:
  \begin{align}
     \varphi_i(\mbalpha_i) &= -\int c_i \nu_{ii} B_{ii} \dv + \mbmu_i\cdot \mbalpha_i \quad \text{for} \quad i=1,2 \quad \text{and} \label{eq:potentialfunction_single} \\
     \varphi(\mbalpha) &= -\int \left[  c_1 \nu_{12} B_{12} +  c_2 \nu_{21} B_{21} \right] \dv  + \mbmu \cdot \mbalpha \label{eq:potentialfunction_mixed}
 \end{align}
and the minimum of \eqref{eq:potentialfunction} is the sum of their minima.
\end{proposition}

\begin{proof}
The statement is trivial since \eqref{eq:potentialfunction} can be written as sum of the three potential functions, whose arguments are independent; that is,  $\varphi_{\rm tot}(\mbalpha_1,\mbalpha_2,\mbalpha) = \varphi_1(\mbalpha_1) + \varphi_2(\mbalpha_2) + \varphi_1(\mbalpha)$.
\end{proof}

The minimization problems in \eqref{eq:potentialfunction_single} and \eqref{eq:potentialfunction_mixed} are numerical analogs of  \eqref{eq:dual_1} and \eqref{eq:dual_2}, respectively. Indeed, the temporal discretization simply introduces the additional weights $c_i^\ell \to 1$  as $\Delta t \to 0$.  Importantly, the existence and uniqueness of solutions to \eqref{eq:potentialfunction_single} and \eqref{eq:potentialfunction_mixed} are guaranteed by the theory in \cite{HHKPW}.  Essentially one need only replace the collision frequencies $\nu_{ij}$ by 
\begin{equation}
    \nu^{\ast}_{ij} = c_i \nu_{ij} = \frac{\nu_{ij}}{1+\gamma\Delta t  (\nu_{ii} + \nu_{ij})} 
\end{equation}
and then verify that $\nu^{\ast}_{ij}$ satisfies the conditions used in \cite{HHKPW}. These conditions are mild integrability conditions that, because $0 < c_i < 1$, are easily satisfied by $\nu_{ij}^{\ast}$ whenever they are satisfied by $\nu_{ij}$.

The minimum of each potential function in \eqref{eq:potentialfunction_single} and \eqref{eq:potentialfunction_mixed} is found using Newton's method for convex optimization. 
The details of this implementation are given in Section 6.

\section{Space discretization}
\label{sec:space}

For the simulations in this paper, we assume a slab geometry for which $\partial_{x^2} f_i = \partial_{x^3} f_i =0$.  Thus while the  (microscopic) velocity space remains three-dimensional $(\mbv = (v^1,v^2,v^3))$, the physical space dimension can be reduced to one dimension; and in a slight abuse of notation, we set $x = x^1$.  We divide the spatial domain $[x_{\min},x_{\max}]$ into uniform cells $I_k = [x_k-\Delta x/2, x_k+\Delta x/2]$ for $k \in \{0,\dots,K\}$. 

We employ a second-order finite volume framework that tracks approximate cell-averaged quantities 
\begin{equation}
    f_{i,k}^\ell \approx \frac{1}{\Delta x}\int_{I_k} f_i(x,\mbv,t^{\ell}) \dx.
\end{equation}
To approximate the relaxation operator, we use the second-order approximation
\begin{equation}
    \mathcal{R}_{i,k}^\ell 
        = \mathcal{R}_{i}(f_{i,k}^\ell,f_{j,k}^\ell) 
        \approx \frac{1}{\Delta x} \int_{I_k} \mathcal{R}\left(\,f_i(x,\mbv,t^{\ell}),f_j(x,\mbv,t^{\ell}) \,\right) \dx.
\end{equation}
Meanwhile, the transport operator $\mathcal{T}$ is discretized with numerical fluxes $\mathscr{F}_{k+\frac{1}{2}}$ by 
\begin{align} \label{eq:TO-discrete}
\mathcal{T}\left( g \right) \approx \mathcal{T}_{k}(g) = \frac{1}{\Delta x} \left( \mathscr{F}_{k+\frac{1}{2}}(g) - \mathscr{F}_{k-\frac{1}{2}}(g) \right)
\end{align}
for any grid function $g = \{g_{k}\}$.  We use, following \cite{MieussensStruchtrup2004}, 
\begin{align} \label{eq:flux}
    \mathscr{F}_{k+\frac{1}{2}}(g) = \frac{v^1}{2} \left( g_{k+1} + g_{k} \right) - \frac{\vert v^1 \vert }{2} \left( g_{k+1} - g_{k} - \phi_{k+\frac{1}{2}}(g) \right)
\end{align}
where $\phi_{k+\frac{1}{2}}$ is a flux limiter. The choice $\phi_{k+\frac{1}{2}}=0$ leads to a first-order approximation (the well-known upwind fluxes). A second-order method is provided by letting
\begin{align}
\label{eq:second-order-limiter}
    \phi_{k+\frac{1}{2}}(g) = \operatorname{minmod}\left( (g_{k}-g_{k-1}), (g_{k+1}-g_{k}), (g_{k+2}-g_{k+1})\right)
\end{align}
where
\begin{align}
    \operatorname{minmod}(a,b,c) = 
    \begin{cases}
     s \min(|a|,|b|,|c|), \quad & \mathrm{sign}(a) = \mathrm{sign}(b)=\mathrm{sign}(c) =:s,\\
     0, \quad &\text{otherwise}.
    \end{cases}
\end{align}
For a simple forward Euler update of \eqref{eq:transport_x}, i.e., 
    \begin{align} \label{eq:update_transport}
    f_{i,k}^{\ell+1} = f_{i,k}^{\ell} - \frac{\Delta t}{\Delta x} \left(\mathscr{F}_{k+\frac{1}{2}}(f_i^{\ell}) - \mathscr{F}_{k-\frac{1}{2}}(f_i^{\ell})\right) ,
    \end{align}
the positivity of $f_i$ is guaranteed by enforcing the CFL condition
\begin{align} \label{eq:CFL}
    \Delta t < \alpha \frac{\Delta x}{\max |v^1|}
\end{align}
with $\alpha=1$ for the first-order flux and $\alpha = \frac{2}{3}$ for the second-order flux.  (See Proposition \ref{prop:positivity_1st_order}.)

\section{Properties of the semi-discrete scheme} \label{sec:discreteproperties}

In this section, we review the positivity, conservation properties, and the entropy behavior of the semi-discrete scheme. A convergence study of the schemes is addressed in the Appendix.

\subsection{Positivity of distribution functions}

The first-order time stepping scheme in Section \ref{subsec:firstordersplit} preserves positivity for both first- and second-order numerical fluxes in space; see Proposition \ref{prop:positivity_1st_order}. 
Additionally, we discuss the positivity for the second-order scheme from Section \ref{subsec:secondorderIMEX} in Proposition \ref{prop:positivity_2nd_order}, and give a sufficient criterion for the space homogeneous case.

\begin{proposition}\label{prop:positivity_1st_order}
The first-order time discretization in Section \ref{subsec:firstordersplit} together with the space discretization described in Section \ref{sec:space} is positivity preserving, provided that
 \begin{align}
 \label{eq:pos_cfl}
     \Delta t \leq \alpha \frac{\Delta x}{\max |v^1|},
 \end{align}
 with $\alpha=1$ and $\alpha = \frac{2}{3}$ for the first-order and second-order fluxes, respectively. 
\end{proposition}
\begin{proof}
Let $f_{i,k}^{\ell}\geq 0$. For the relaxation step,
\begin{align}
    f_{i,k}^{\ell'} \overset{\eqref{eq:update_split}}{=} c_{i,k}^{\ell} f_{i,k}^{\ell} + c_{i,k}^{\ell} \Delta t (\nu_{ii,k}^{\ell} A_{ii,k}^{\ell'} + \nu_{ij,k}^{\ell} A_{ij,k}^{\ell'}) \geq 0
\end{align}
because $c_{i,k}^{\ell}, \nu_{ii,k}^{\ell}, \nu_{ij,k}^{\ell}, A_{ii,k}^{\ell'}, A_{ij,k}^{\ell'}\geq 0$.
For the transport step \eqref{eq:update_transport}, we have with the first-order fluxes
\begin{align}
    f_{i,k}^{\ell+1} = \left( 1-\frac{\Delta t}{\Delta x}|v^1|\right) f_{i,k}^{\ell'} + \frac{\Delta t}{\Delta x}|v^1| f_{i,k-\operatorname{sign}(v^1)}^{\ell'} \geq 0,
\end{align}
where the last inequality holds in each cell provided that the given CFL condition in \eqref{eq:pos_cfl} holds with $\alpha=1$.

For the second-order fluxes, define $\sigma:=\operatorname{sign}(f_{i,k}^\ell-f_{i,k-1}^\ell)$. Then one can show that
\begin{align}
    \phi_{k+\frac{1}{2}}(f_i^\ell) \geq 
    \begin{cases}
        0 \quad &\text{if} \quad \sigma=+1 \\
        f_{i,k+1}^\ell-f_{i,k}^\ell \quad &\text{if} \quad \sigma=-1
    \end{cases},\\
    -\phi_{k-\frac{1}{2}}(f_i^\ell) \geq 
    \begin{cases}
        f_{i,k-1}^\ell-f_{i,k}^\ell \quad &\text{if} \quad \sigma=+1\\
        0 \quad &\text{if} \quad \sigma=-1 
    \end{cases}.
\end{align}
Hence
\begin{equation}
\begin{aligned}
    f_{i,k}^{\ell+1} 
    &\overset{\eqref{eq:update_transport}}{=} \left( 1-\frac{\Delta t}{\Delta x}|v^1|\right) f_{i,k}^{\ell} + \frac{\Delta t}{\Delta x}|v^1| f_{i,k-\operatorname{sign}(v^1)}^{\ell} + \frac{\Delta t}{\Delta x} \frac{|v^1|}{2} (\phi_{k+\frac{1}{2}}(f_i^\ell)-\phi_{k-\frac{1}{2}}(f_i^\ell)) \\
    &\geq \left( 1-\frac{\Delta t}{\Delta x}|v^1|\right) f_{i,k}^{\ell} + \frac{\Delta t}{\Delta x}|v^1| f_{i,k-\operatorname{sign}(v^1)}^{\ell} + \frac{\Delta t}{\Delta x} \frac{|v^1|}{2}
    \begin{cases}
        (f_{i,k-1}^\ell-f_{i,k}^\ell) \; &\text{if} \quad \sigma=+1 \\
        (f_{i,k+1}^\ell-f_{i,k}^\ell) \; &\text{if} \quad \sigma=-1
    \end{cases} \\
    &= \left( 1-\frac{3}{2}\frac{\Delta t}{\Delta x}|v^1|\right) f_{i,k}^{\ell} + \frac{\Delta t}{\Delta x}|v^1| f_{i,k-\operatorname{sign}(v^1)}^{\ell} + \frac{\Delta t}{\Delta x} \frac{|v^1|}{2}
    \begin{cases}
        f_{i,k-1}^\ell \; &\text{if} \quad \sigma=+1 \\
        f_{i,k+1}^\ell \; &\text{if} \quad \sigma=-1
    \end{cases} \\
    &\geq 0,
\end{aligned}
\end{equation}
provided that the CFL condition in \eqref{eq:pos_cfl} holds with $\alpha= \frac{2}{3}$.
\end{proof}

It is more difficult to guarantee positivity with second-order time-stepping. Unconditionally strong stability preserving (SSP) implicit Runge-Kutta schemes, which preserve any convex property, e.g. positivity, are at most first-order accurate \cite{GottliebShuTadmor2001}. Modified IMEX Runge-Kutta schemes that preserve positivity for the classical single-species BGK equation have been recently developed in \cite{HSZ}.  However, to our knowledge, these schemes cannot be applied directly to BGK models with velocity-dependent collision frequencies. 

Nevertheless, we derive some sufficient conditions on $ \Delta t$ for positivity preservation in the second-order scheme presented in Section \ref{subsec:secondorderIMEX}.

\begin{proposition} \label{prop:positivity_2nd_order}
For the space homogeneous case, the second-order IMEX scheme presented in Section \ref{subsec:secondorderIMEX} is positivity preserving provided that
\begin{align} \label{eq:timerestriction-IMEX}
\Delta t \leq \frac{1}{(1-2\gamma)  (\nu_{ii}^{(1)}+\nu_{ij}^{(1)})}
\end{align}
 for $i,j=1,2.$
\end{proposition}
\begin{proof}
The positivity of $f_i^{(1)}$ follows directly from its definition without any restriction on the time step. 
For the positivity of $f_i^{\ell+1}=f_i^{(2)}$ we require  $G_i^{(2)} \geq 0$. 
Using the definition of $f_i^{(1)}$, we obtain
\begin{align}
 0\leq G_i^{(2)} &= f_i^{\ell} + \Delta t (1-\gamma) \left[ \nu_{ii}^{(1)}A_{ii}^{(1)} + \nu_{ij}^{(1)}A_{ij}^{(1)} - (\nu_{ij}^{(1)}+\nu_{ij}^{(1)}) f_i^{(1)}\right] \\
 & \, = f_i^{\ell} \left[1-\Delta t (1-\gamma) c_i^{(1)} (\nu_{ii}^{(1)}+\nu_{ij}^{(1)}) \right] + \Delta t (1-\gamma) c_i^{(1)} \left[\nu_{ii}^{(1)}A_{ii}^{(1)} + \nu_{ij}^{(1)}A_{ij}^{(1)} \right].
\end{align}
Then, the most obvious sufficient condition for positivity reads
\begin{align} 
\label{eq:cfl-order-two}
1-\Delta t (1-\gamma) c_i^{(1)} (\nu_{ii}^{(1)}+\nu_{ij}^{(1)}) \geq 0 \quad \Longleftrightarrow \quad 
\Delta t \leq \frac{1}{(1-2\gamma)  (\nu_{ii}^{(1)}+\nu_{ij}^{(1)})}.
\end{align}
\end{proof}

The time step condition \eqref{eq:timerestriction-IMEX} can be restrictive if $\nu_{ij}^{(1)}$ become large. For this reason, one may instead enforce the milder (but still sufficient) local conditions
\begin{align} \label{eq:timestep-restriction-local-IMEX-1}
\Delta t \leq \frac{f_i^{\ell}}{(1-2\gamma)  (\nu_{ii}^{(1)}+\nu_{ij}^{(1)})f_i^{\ell} - (1-\gamma)( \nu_{ii}^{(1)}A_{ii}^{(1)} + \nu_{ij}^{(1)}A_{ij}^{(1)})}
\end{align}
and
\begin{align}\label{eq:timestep-restriction-local-IMEX-2}
\Delta t \leq \frac{f_i^{\ell}}{(1-\gamma)  \left[(\nu_{ii}^{(1)}+\nu_{ij}^{(1)})f_i^{(1)} - (\nu_{ii}^{(1)}A_{ii}^{(1)} + \nu_{ij}^{(1)}A_{ij}^{(1)}) \right]}.
\end{align}
When the frequencies are large, the difference between each numerical kinetic distribution and its corresponding target function is to scale with the inverse of the frequency, in which case \eqref{eq:timestep-restriction-local-IMEX-1} and \eqref{eq:timestep-restriction-local-IMEX-2} are not restrictive.

In our numerical tests, the time step $\Delta t$ is set according to the CFL condition \eqref{eq:CFL} by default. If positivity is violated, we reduce the time step size according to (\ref{eq:timerestriction-IMEX}). Thus we guarantee positivity while maintaining large time steps whenever possible.  One could instead use the less restrictive local conditions in \eqref{eq:timestep-restriction-local-IMEX-1} and  \eqref{eq:timestep-restriction-local-IMEX-2}, which requires additional iterations over the grid to find a global value for the time step.  However, in practice, violations of positivity are rare and thus we use \eqref{eq:timerestriction-IMEX} for simplicity.

\subsection{Conservation of mass, total momentum and total energy} \label{subsec:discreteconservation}

In this section, we address the conservation of mass, total momentum, and total energy for the semi-discrete scheme (before velocity discretization).

\begin{proposition} \label{theo:conservation_relaxation}
The relaxation step in the first-order splitting scheme presented in Section \ref{subsec:firstordersplit} satisfies the conservation laws
\begin{gather}
   \int m_1 f_1^{\ell'} d \mbv = \int m_1 f_1^{\ell} \dv , \quad \int  m_2 f_2^{\ell'} \dv = \int  m_2 f_2^{\ell} \dv, \\
    \int  \left( m_1 \mbv f_1^{\ell'} +  m_2 \mbv f_2^{\ell'} \right) \dv
    = \int \left(  m_1 \mbv f_1^{\ell} + m_2 \mbv f_2^{\ell}   \right) \dv, \label{eq:cons_momentum_discrete}\\
        \int  \left( m_1 |\mbv|^2 f_1^{\ell'} +  m_2 |\mbv|^2 f_2^{\ell'} \right) \dv
    = \int  \left(  m_1 |\mbv|^2 f_1^{\ell} + m_2 |\mbv|^2 f_2^{\ell}   \right) \dv. \label{eq:cons_energy_discrete}
\end{gather}
\end{proposition}
\begin{proof}
We multiply the relaxation step \eqref{eq:update_split} by $\mba_i$, sum over $i=1,2$, and integrate with respect to $\mbv$. Sorting terms yields
\begin{align}
\begin{split}
\int & \left( f_1^{\ell'} \mba_1+ f_2^{\ell'} \mba_2\right)  \dv  -
\int \left(   f_1^{\ell} \,\mba_1+  f_2^{\ell} \,\mba_2\right) \dv  \\ 
&\overset{\eqref{eq:update_split}}{=}
 \Delta t   \left[ \int \left(  c_1^{\ell} \nu_{11}^{\ell} A_{11}^{\ell'}\mba_1 +  c_2^{\ell} \nu_{22}^{\ell} A_{22}^{\ell'} \mba_2 +  c_1^{\ell} \nu_{12}^{\ell} A_{12}^{\ell'} \mba_1+  c_2^{\ell} \nu_{21}^{\ell} A_{21}^{\ell'} \mba_2\right) \dv \right. \\ 
&\hspace*{1.2cm}\left. -
\int \left[  c_1^{\ell} \left(\nu_{11}^{\ell}+\nu_{12}^{\ell}\right) f_1^{\ell} \,\mba_1 +  c_2^{\ell} \left(\nu_{22}^{\ell}+\nu_{21}^{\ell}\right) f_{2}^{\ell}\,\mba_2 \right] \dv   \right]. 
\end{split}
\end{align} 
The right-hand side above corresponds to the first-order optimality conditions in \eqref{eq:moment_A_general}. By minimizing the corresponding functions in \eqref{eq:potentialfunction_single} and \eqref{eq:potentialfunction_mixed}, we guarantee that this term is identically zero, which in turn proves the conservation statement \eqref{eq:cons_momentum_discrete} and \eqref{eq:cons_energy_discrete}.  For the masses we execute the above procedure for each species individually and obtain
\begin{align}
    \int f_i^{\ell'} m_i \dv - \int f_i^{\ell} m_i \dv  &\overset{\eqref{eq:update_split}}{=}
 \Delta t   \int \left(  c_i^{\ell} \nu_{ii}^{\ell} A_{ii}^{\ell'}m_i +  c_i^{\ell} \nu_{ij}^{\ell} A_{ij}^{\ell'} m_i \right) \dv 
- \int   c_i^{\ell} \left(\nu_{ii}^{\ell}+\nu_{ij}^{\ell}\right) f_i^{\ell} \,m_i  \dv \\
&= \Delta t \left[ \partial_{\lambda_i^0} \varphi_i(\mblambda_i) +  \partial_{\lambda_{ij}^0} \varphi(\mblambda)  \right] 
\end{align}
which vanishes due to first-order optimality conditions on $\varphi$ and $\varphi_i$ being defined in \eqref{eq:potentialfunction_single} and \eqref{eq:potentialfunction_mixed}.
\end{proof}

\begin{proposition} \label{theo:conservation_transport}
For each $i=1,2$, the transport step in the first-order splitting scheme in Section \ref{subsec:firstordersplit}, combined with the space discretization presented in Section \ref{sec:space} satisfies the conservation laws
\begin{align}
    \sum_{k=0}^K \int \mba_i f^{\ell+1}_{i,k} \dv \Delta x = \sum_{k=0}^K \int \mba_i f^{\ell'}_{i,k} \dv \Delta x
\end{align}
for periodic or zero boundary conditions. 
\end{proposition}
\begin{proof}
For $i=1,2$, we multiply the transport step \eqref{eq:transport_x} by $\mba_i$, integrate with respect to $\mbv$ and sum over all cell averages in $x$. The result is
\begin{align}
\sum_{k=0}^K \int \mba_i f^{\ell+1}_{i,k} \dv \Delta x &\overset{\eqref{eq:transport_x}}{=}  
\sum_{k=0}^K \int \mba_i f^{\ell'}_{i,k} \dv \Delta x
- \sum_{k=0}^{K} \int \frac{\Delta t}{\Delta x} \left( \mathscr{F}_{k+\frac{1}{2}}(f_i^{\ell'}) - \mathscr{F}_{k-\frac{1}{2}}(f_i^{\ell'}) \right) \mba_i\dv \Delta x\\
&=  \sum_{k=0}^K \int \mba_i f^{\ell'}_{i,k} \dv \Delta x - \Delta t \, \mathbf{\Omega }
\end{align}
where the remnant of the telescoping sum 
\begin{align}
    \mathbf{\Omega } = \int \mba_i  \mathscr{F}_{K+\frac{1}{2}}(f_i^{\ell'}) \dv - \int \mba_i \mathscr{F}_{-\frac{1}{2}}(f_i^{\ell'}) \dv 
\end{align}
vanishes for periodic or zero boundary conditions, e.g. $\mathscr{F}_{K+\frac{1}{2}}(f_i^{\ell}) = \mathscr{F}_{-\frac{1}{2}}(f_i^{\ell})$  and $\mathscr{F}_{K+\frac{1}{2}}(f_i^{\ell}) = \mathscr{F}_{-\frac{1}{2}}(f_i^{\ell}) = 0$, respectively.
\end{proof}

The second-order time-stepping scheme in Section \ref{subsec:secondorderIMEX} can be broken into relaxation and transport parts, each of which preserves the conservation of mass, total momentum, and total energy.  As a result, we have the following.
\begin{corollary}
For periodic or zero boundary conditions, any combination of temporal and space discretization presented in Sections \ref{sec:time} and \ref{sec:space}, respectively, conserves mass, total momentum and total energy. 
\end{corollary}

\subsection{Entropy inequality}

We discuss the entropy behavior for the first-order scheme in Section \ref{subsec:firstordersplit}. Both the relaxation and the transport step dissipate entropy; see Propositions \ref{theo:entropy_inequality_relaxation} and \ref{theo:entropy_inequality_transport}. Additionally, we show in Proposition \ref{theo:entropy_equilibrium} that the minimal entropy is reached for the relaxation step if the distribution functions coincide with the corresponding target functions. These results rely on dissipative properties of the backward Euler method and a clever use of the conservation properties at the semi-discrete level.

\begin{proposition}\label{theo:entropy_inequality_relaxation}
 Let $h(f) = f \log f -f$. The relaxation step in the first-order splitting scheme in Section \ref{subsec:firstordersplit} fulfills the discrete entropy inequality
 \begin{align}
 \label{eq:entropy_inequality}
     \int  h(f_1^{\ell'}) + h(f_2^{\ell'})\dv \leq  \int h(f_1^{\ell}) + h(f_2^{\ell}) \dv.
 \end{align}
\end{proposition}
\begin{proof}
By convexity
\begin{equation}
h(f_i^{\ell}) \geq h(f_i^{\ell'}) + h'(f_i^{\ell'})(f_i^{\ell}- f_i^{\ell'}).
\end{equation}
The implicit step \eqref{eq:update_split} is 
\begin{equation}
\label{eq:implicit_step_restated}
f_i^{\ell'} - f_i^{\ell} = \Delta t \nu_{ii}^{\ell}(A_{ii}^{\ell'} -f_i^{\ell'}) + \Delta t \nu_{ij}^{\ell} (A_{ij}^{\ell'} -f_i^{\ell'}).
\end{equation}
Using \eqref{eq:implicit_step_restated} and the convexity of $h$ gives
\begin{equation}
\label{eq:h_convex}
\begin{aligned}
h(f_i^{\ell'}) - h(f_i^{\ell}) &\leq h'(f_i^{\ell'})(f_i^{\ell'}- f_i^{\ell}) \\
    &\overset{\eqref{eq:implicit_step_restated}}{=}  \Delta t \,\nu_{ii}^{\ell} h'(f_i^{\ell'})(A_{ii}^{\ell'} -f_i^{\ell'}) + \Delta t \nu_{ij}^{\ell} h'(f_i^{\ell'})(A_{ij}^{\ell'} -f_i^{\ell'}) \\
    &= \Delta t \,\nu_{ii}^{\ell} \left[ h'\left(\frac{f_i^{\ell'}}{A_{ii}^{\ell'}}\right)(A_{ii}^{\ell'} -f_i^{\ell'})
     + h'(A_{ii}^{\ell'})(A_{ii}^{\ell'} -f_i^{\ell'}) \right] \\
    &\hspace{0.5cm}+\Delta t \,\nu_{ij}^{\ell} \left[ h'\left(\frac{f_i^{\ell'}}{A_{ij}^{\ell'}}\right)(A_{ij}^{\ell'} -f_i^{\ell'})
    + h'(A_{ij}^{\ell'})(A_{ij}^{\ell'} -f_i^{\ell'}) \right],
\end{aligned}
\end{equation}
where in the last line above, we have added and subtracted the same quantity.  After integration in $\mbv$, some terms in \eqref{eq:h_convex} disappear.  Specifically, because 
\begin{align}
h'(A_{ii}^{\ell'}) = \log(A_{ii}^{\ell'}) = \mblambda_i \cdot \mba(\mbv)
\end{align}
it follows that
\begin{equation} \label{eq:zero-intra}
\int \nu_{ii}^{\ell} \, h'(A_{ii}^{\ell'})(A_{ii}^{\ell'} -f_i^{\ell'}) \dv 
    =\mblambda_i \cdot \int \nu_{ii}^{\ell} \,  \mba(\mbv)(A_{ii}^{\ell'} -f_i^{\ell'}) \dv 
    =0. 
\end{equation}
Analogously for the inter-species terms,
\begin{equation} \label{eq:zero-inter}
\begin{aligned}
    \int \nu_{12}^{\ell} \, &h'(A_{12}^{\ell'})(A_{12}^{\ell'} -f_1^{\ell'}) \dv + \int \nu_{21}^{\ell} \, h'(A_{21}^{\ell'})(A_{21}^{\ell'} -f_2^{\ell'}) \dv \\
    &=\lambda_{12}^0 \int \nu_{12}^{\ell} \,  (A_{12}^{\ell'} -f_1^{\ell'}) \dv + \lambda_{21}^0 \int \nu_{21}^{\ell} \,  (A_{21}^{\ell'} -f_2^{\ell'}) \dv \\
    &\quad + \begin{pmatrix} \mblambda^1 \\ \lambda^2 \end{pmatrix} \cdot \int (\nu_{12}^{\ell} \,  (A_{12}^{\ell'} -f_1^{\ell'}) + \nu_{21}^{\ell} \,  (A_{21}^{\ell'} -f_2^{\ell'}) ) \begin{pmatrix} \mbv \\ |\mbv|^2 \end{pmatrix} \dv  \\
    &= 0. 
\end{aligned}
\end{equation}
The integrals \eqref{eq:zero-intra} and \eqref{eq:zero-inter} vanish as the conservation properties are satisfied at the semi-discrete level as well by construction of the scheme.
Thus after integrating \eqref{eq:h_convex} in $\mbv$,
\begin{equation}
\label{eq:time_discrete_entropy_inequality}
\begin{aligned}
\int & h(f_1^{\ell'}) \dv - \int h(f_1^{\ell}) \dv + \int h(f_2^{\ell'}) \dv - \int h(f_2^{\ell}) \dv \\
    &\leq \Delta t \,\nu_{11}^{\ell} \int h'\left(\frac{f_1^{\ell'}}{A_{11}^{\ell'}}\right)(A_{11}^{\ell'} -f_1^{\ell'}) \dv + \Delta t \,\nu_{22}^{\ell} \int h'\left(\frac{f_2^{\ell'}}{A_{22}^{\ell'}}\right)(A_{22}^{\ell'} -f_2^{\ell'}) \dv \\
    &\quad + \Delta t \,\nu_{12}^{\ell} \int h'\left(\frac{f_1^{\ell'}}{A_{12}^{\ell'}}\right)(A_{12}^{\ell'} -f_1^{\ell'}) \dv + \Delta t\, \nu_{21}^{\ell} \int h'\left(\frac{f_2^{\ell'}}{A_{21}^{\ell'}}\right)(A_{21}^{\ell'} -f_2^{\ell'}) \dv \\
    &\leq 0
\end{aligned}
\end{equation}
because $\log\left(\frac{x}{y}\right)(y-x)\leq 0$ for all $x,y\in\mathbb{R}^+$.  
\end{proof}

\begin{proposition} \label{theo:entropy_equilibrium}
 The inequality in Proposition \ref{theo:entropy_inequality_relaxation} is an equality if and only if $f_1^{\ell} = A_{12}^{\ell}$ and  $f_2^{\ell} = A_{21}^{\ell}$.  In such cases $f_1^{\ell'} = A_{12}^{\ell'}$ and  $f_2^{\ell'} = A_{21}^{\ell'}$.
\end{proposition} 

\begin{proof}
Suppose first that $f_1^{\ell} = A_{12}^{\ell}$ and  $f_2^{\ell} = A_{21}^{\ell}$. Then according to \cite[Theorem 2]{HHKPW}, 
\begin{equation}
    h(f_1^{\ell}) + h(f_2^{\ell}) \leq h(g_1) + h(g_2) 
\end{equation}
for any measurable positive functions $g_1$ and $g_2$ such that
\begin{gather}
   \int m_1 g_1 d \mbv = \int m_1 f_1^{\ell} \dv , \quad \int  m_2 g_2 \dv = \int  m_2 f_2^{\ell} \dv, \\
    \int  \left( m_1 \mbv g_1 +  m_2 \mbv g_2 \right) \dv
    = \int \left(  m_1 \mbv f_1^{\ell} + m_2 \mbv f_2^{\ell}   \right) \dv,\\
        \int  \left( m_1 |\mbv|^2 g_1 +  m_2 |\mbv|^2 g_2 \right) \dv
    = \int  \left(  m_1 |\mbv|^2 f_1^{\ell} + m_2 |\mbv|^2 f_2^{\ell}   \right) \dv.
\end{gather}
These conditions are exactly those satisfied by $f_1^{\ell'}$ and $f_1^{\ell'}$ (cf. Theorem \ref{theo:conservation_relaxation}).  Hence
\begin{equation}
    h(f_1^{\ell}) + h(f_2^{\ell}) \leq h(f_1^{\ell'}) + h(f_1^{\ell'}) 
\end{equation}
which shows that \eqref{eq:entropy_inequality} is an equality.  To show the converse statement, suppose that \eqref{eq:entropy_inequality} holds as an equality.   Then according to \eqref{eq:time_discrete_entropy_inequality}
$f_1^{\ell'} = A_{11}^{\ell'} = A_{12}^{\ell'}$ and $f_2^{\ell'} = A_{21}^{\ell'} = A_{22}^{\ell'}$.  Therefore, by definition of $\mathcal{R}_i$ in \eqref{eq:R-def}, $\mathcal{R}_{i}(f_i^{\ell'},f_j^{\ell'})=0$, which when plugged into \eqref{eq:relax}, gives  $f_1^{\ell} = f_1^{\ell'}$ and $f_2^{\ell} = f_2^{\ell'}$.
\end{proof}

\begin{proposition}\label{theo:entropy_inequality_transport}
 Let $h(f) = f \log f -f$. The transport step in the first-order splitting scheme in Section \ref{subsec:firstordersplit} combined with the first-order spatial discretization in Section \ref{sec:space} fulfills the discrete entropy inequality
 \begin{align}
 \label{eq:transport_entropy_inquality}
    \sum_{k=0}^K \left\{ \int  h(f_{1,k}^{\ell+1}) + h(f_{2,k}^{\ell+1})\dv \right\} \Delta x
        \leq \sum_{k=0}^K  \left\{   \int h(f_{1,k}^{\ell'}) + h(f_{2,k}^{\ell'}) \dv \right\} \Delta x
 \end{align}
 for periodic or zero boundary conditions, provided that 
  \begin{align}
     \Delta t \leq \frac{\Delta x}{\max |v^1|}.
 \end{align}
\end{proposition}
\begin{proof}
Using the notation $v^+ := \frac{v^1+|v^1|}{2}$ and $v^- := \frac{v^1-|v^1|}{2}$ we write the update formula of \eqref{eq:transport_x} with the first-order numerical fluxes as
\begin{equation}
\begin{aligned}
    f_{i,k}^{\ell+1} &= f_{i,k}^{\ell'} - \frac{\Delta t}{\Delta x}\left(v^+ f_{i,k}^{\ell'} + v^- f_{i,k+1}^{\ell'} - v^+ f_{i,k-1}^{\ell'} - v^- f_{i,k}^{\ell'}\right) \\
    &= \left(1-\frac{\Delta t}{\Delta x}|v^1|\right) f_{i,k}^{\ell'} - \frac{\Delta t}{\Delta x} v^- f_{i,k+1}^{\ell'} + \frac{\Delta t}{\Delta x} v^+ f_{i,k-1}^{\ell'}.
\end{aligned}
\end{equation}
Clearly if the CFL condition is fulfilled, then $f_{i,k}^{\ell+1}$ is a convex linear combination of $f_{i,k}^{\ell'}$, $f_{i,k-1}^{\ell'}$, and $f_{i,k+1}^{\ell'}$. Thus by the convexity of $h$, for each $\mbv$,
\begin{equation}
\label{eq:entropy_transport_v}
\begin{aligned}
    \sum_{k=0}^K h (f_{i,k}^{\ell+1})  \Delta x
        & \leq \sum_k  \left[\left(1-\frac{\Delta t}{\Delta x}|v^1|\right)  h (f_{i,k}^{\ell'}) - \frac{\Delta t}{\Delta x} v^- h (f_{i,k+1}^{\ell'}) + \frac{\Delta t}{\Delta x} v^+ h (f_{i,k-1}^{\ell'})\right] \Delta x \\
        &= \sum_k h (f_{i,k}^{\ell'}) \Delta x + \Delta t \,\Omega 
\end{aligned}
\end{equation}
where the boundary term
\begin{align}
    \Omega = 
     v^-f_{i,0}^{\ell'}\log (f_{i,0}^{\ell'}) - v^-f_{i,K+1}^{\ell'}\log (f_{i,K+1}^{\ell'})  - v^+ f_{i,K}^{\ell'}\log (f_{i,K}^{\ell'})  + v^+ f_{i,-1}^{\ell'}\log (f_{i,-1}^{\ell'}). 
\end{align}
is the only remnant of the telescoping sum and vanishes for periodic or zero boundary conditions.
Thus summation over $i$ and integration of \eqref{eq:entropy_transport_v} with respect to $\mbv$ yields the entropy inequality in \eqref{eq:transport_entropy_inquality}.

\end{proof}
Combining the two results above gives the following: 
\begin{corollary}
  For periodic or zero boundary conditions, the first-order splitting scheme from Section \ref{subsec:firstordersplit} combined with the first-order numerical fluxes in Section \ref{sec:space} fulfills the discrete entropy inequality
 \begin{align}
    \sum_{k=0}^K \left\{ \int  h(f_{1,k}^{\ell+1}) + h(f_{2,k}^{\ell+1})\dv \right\} \Delta x
        \leq \sum_{k=0}^K  \left\{   \int h(f_{1,k}^{\ell'}) + h(f_{2,k}^{\ell'}) \dv \right\} \Delta x
 \end{align}
provided that 
  \begin{align}
     \Delta t \leq \frac{\Delta x}{\max |v^1|}.
 \end{align}
\end{corollary}

\section{Velocity discretization} 
\label{sec:velocity}

In order to obtain a fully-discrete scheme, we finally discretize the velocity variable. We center the discrete velocities
  $\mbv_\mbq = (v_{q_1}^1, v_{q_2}^2,v_{q_3}^3)^\top$, with $\mbq=(q_1,q_2,q_3) \in \mathbb{N}^3_0$, around the mixture mean velocity $\mbu_{\rm mix}$ and restrict them to a finite cube.  That is, for each $p \in \{1,2,3\}$,
\begin{align}
v^p \in [u_{{\rm mix}}^p - 6 v_{{\rm{th}},i} , u_{{\rm mix}}^p + 6v_{{\rm{th}},i} ]
\end{align} 
where $v_{{\rm{th}},i} = \sqrt{T_{\rm mix}/m_i}$ is the thermal velocity of species $i$.
To ensure adequate resolution of the velocity domain, the velocity mesh size is chosen, as in \cite{mieussens2000discrete}, to be $\Delta v = 0.25 v_{{\rm{th}},i}$ in each direction.  

An advantage of the BGK model \cite{HHM} is that it is possible to use different velocity grids for each species/equation since the distributions of different species only interact via their moments. 
This feature is a substantial benefit when the species masses, and hence the reference thermal speeds for each species, differ significantly.

Using the grid described above, all velocity integrals are replaced by discrete sums using the trapezoidal rule, which is known to perform well for smooth, compactly supported functions, since they can be viewed as periodic.  (See, e.g, \cite[Section 5.4, Corollary 1]{atkinson1989introduction}.)  Thus
\begin{align}
    \int (\cdot) \dv \approx \sum_\mbq \omega_\mbq (\cdot)_\mbq (\Delta v)^3 ,
\end{align}
where $\omega_\mbq=\omega_{q_1}\omega_{q_2}\omega_{q_3}$ are the weights and
\begin{align}
\omega_{q_p} = \begin{cases} 1 \quad \text{if } \min(q_p) < q_p < \max(q_p),\\
\frac{1}{2} \quad \text{else}.
\end{cases}
\end{align}
Due to the quadrature approximation, we have to distinguish between discrete and continuous moments, especially when determining the local equilibria $A_{ii}$ and $A_{ij}$.  In fact, the minimization of \eqref{eq:potentialfunction_single} and \eqref{eq:potentialfunction_mixed} is solved using a discrete velocity grid and discrete moments $\bar{\mbmu}_1,\bar{\mbmu}_2,\bar{\mbmu}$ as input. Thus $\mblambda_1, \mblambda_2, \mblambda$ are such that $A_{ii}$ and $A_{ij}$ have the desired discrete moments and the conservation properties from the previous section are fulfilled at the discrete level. (See \cite{mieussens2000discrete} for a similar approach for the standard, singles-species BGK model.)

\begin{theorem} \label{theo:properties_v_discrete}
Propositions \ref{prop:positivity_1st_order}, \ref{prop:positivity_2nd_order}, and \ref{theo:entropy_inequality_relaxation}-\ref{theo:entropy_inequality_transport} all hold true after replacing continuous integrals by their respective quadratures. 
Additionally, the scheme in Section \ref{sec:generaltime} satisfies the following conservation properties for $\ell\geq 0$ 
\begin{align}
    \sum_{k,\mbq}  \omega_\mbq\left(f_{1,k\mbq}^\ell \mba_{1,\mbq} + f_{2,k\mbq}^\ell \mba_{2,\mbq}\right) (\Delta v)^3 \Delta x = \sum_{k,\mbq}  \omega_\mbq \left(f_{1,k\mbq}^0 \mba_{1,\mbq} + f_{2,k\mbq}^0 \mba_{2,\mbq}\right) (\Delta v)^3 \Delta x
\end{align}
with $\mba_{i,\mbq} = m_i(1,\mbv_\mbq,|\mbv_\mbq|^2)^\top$ and $f_{i,k\mbq}^\ell \approx f_{i,k}^\ell(\mbv_\mbq)$.
\end{theorem}

\paragraph{Optimization algorithm} The minimization of \eqref{eq:potentialfunction_single} and \eqref{eq:potentialfunction_mixed} is solved by Newton's method with a backtracking line search \cite[p. 325]{DennisSchnabel1996linesearch}, using the \texttt{SNESNEWTONLS} solver from PETSc \cite{petsc-web-page,petsc-efficient,petsc-user-ref}.
Newton's methods require the evaluation of gradients:
\begin{align} 
    \nabla_{\mbalpha_i} \varphi_i &\approx -\sum_\mbq \omega_\mbq (c_i  \nu_{ii})_\mbq \,B_{ii,\mbq}\, \mba_{i,\mbq} (\Delta v)^3 
    +  \bar{\mbmu}_{i}, \label{eq:gradient-intra}\\
    \nabla_{\mbalpha} \varphi &\approx  -\sum_\mbq\omega_\mbq \left((c_1\nu_{12})_\mbq \, B_{12,\mbq} \,\mba_{12,\mbq} + (c_2\nu_{21})_\mbq \, B_{21,\mbq}\, \mba_{21,\mbq}\right) (\Delta v)^3
    + \bar{\mbmu},
 \end{align}   
and Hessians:
 \begin{align}
    \nabla_{\mbalpha_i}^2 \varphi_i &\approx -\sum_\mbq \omega_\mbq (c_i  \nu_{ii})_\mbq \, B_{ii,\mbq}\, \mba_{i,\mbq} \otimes \mba_{i,\mbq} (\Delta v)^3,  \\
    \nabla_{\mbalpha}^2 \varphi &\approx  -\sum_\mbq \omega_\mbq \left((c_1\nu_{12})_\mbq \, B_{12,\mbq} \,\mba_{12,\mbq}\otimes \mba_{12,\mbq} + (c_2\nu_{21})_\mbq \, B_{21,\mbq}\, \mba_{21,\mbq}\otimes \mba_{21,\mbq} \right) (\Delta v)^3, \label{eq:Hessian-inter}
\end{align}
where 
$\mba_{12,\mbq} = m_1 (1,0,\mbv_\mbq,|\mbv_\mbq|^2)^\top$ and $\mba_{21,\mbq}= m_2 (0,1,\mbv_\mbq,|\mbv_\mbq|^2)^\top$. The input data in \eqref{eq:input_mu_i} and \eqref{eq:input_mu} is computed in a straightforward way:
\begin{align}
    \bar{\mbmu}_{i} \approx \sum_\mbq \omega_\mbq (c_i  \nu_{ii})_\mbq \,G_{i,\mbq}\, \mba_{i,\mbq} (\Delta v)^3,
\end{align}
and analogously for $\bar{\mbmu}$. The Newton method is considered to have converged if one of the standard termination criteria\footnote{For solving $F(x)=0$, standard termination criteria are: i) $||F||< \epsilon$, ii) $||F||< \epsilon ||F(x_0)||$, and iii) $||\Delta x|| < \epsilon||x||$ for the tolerance $\epsilon$.} is less than $10^{-14}$.

The quadrature operations required to evaluate the gradients and Hessians \eqref{eq:gradient-intra}--\eqref{eq:Hessian-inter} of the optimization algorithm constitute the major cost of the overall scheme. As discussed in the introduction, the optimization algorithm is more expensive than solving the standard BGK model but less expensive than computing the Boltzmann collision operator.

\section{Numerical results}
\label{sec:results}

In this section, we perform a range of numerical tests. We first verify the  properties of our scheme and then present
several examples to illustrate the effect of a velocity-dependent collision frequency.

\subsection{Relaxation in a homogeneous setting}

\subsubsection{Illustrative toy problem}\label{test:verification}
\begin{figure}[htb]
\centering
\begin{subfigure}[c]{0.32\textwidth}
\includegraphics[width=\textwidth]{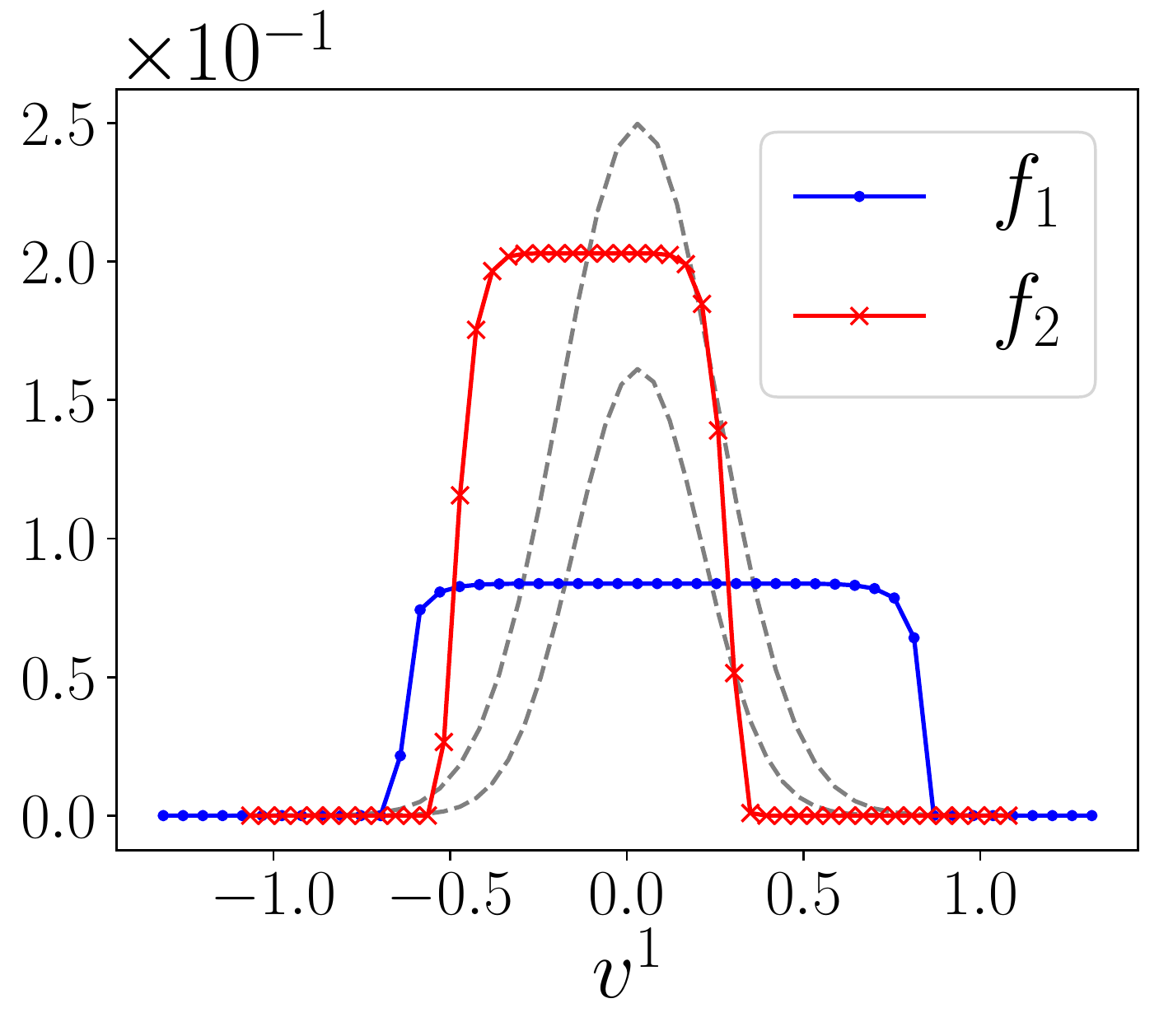}
\subcaption{$t=0$}
\label{fig:h-theorem_ic}
\end{subfigure}
\hfill
\begin{subfigure}[c]{0.32\textwidth}
\includegraphics[width=\textwidth]{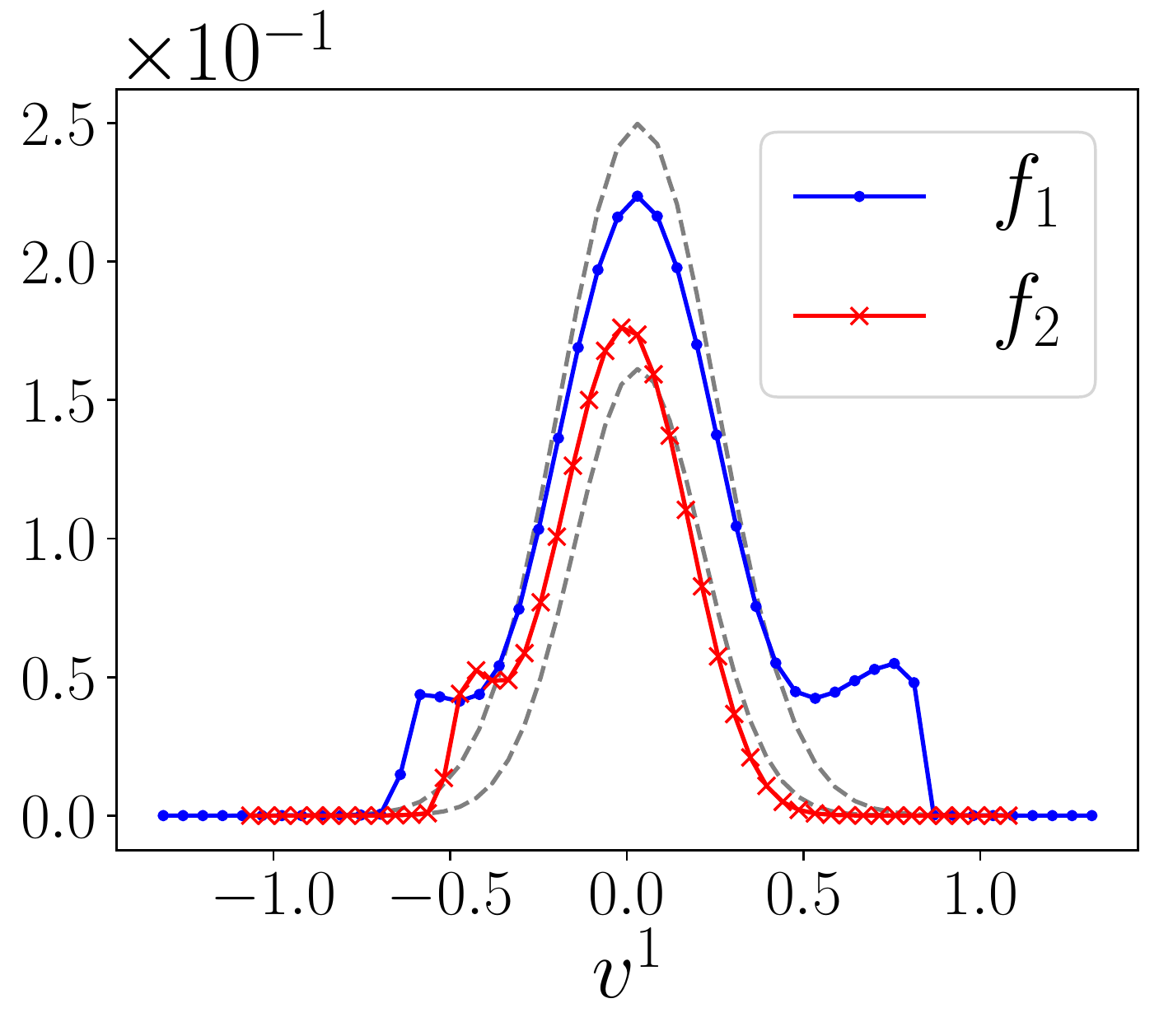}
\subcaption{$t=0.25$}
\end{subfigure}
\hfill
\begin{subfigure}[c]{0.32\textwidth}
\includegraphics[width=\textwidth]{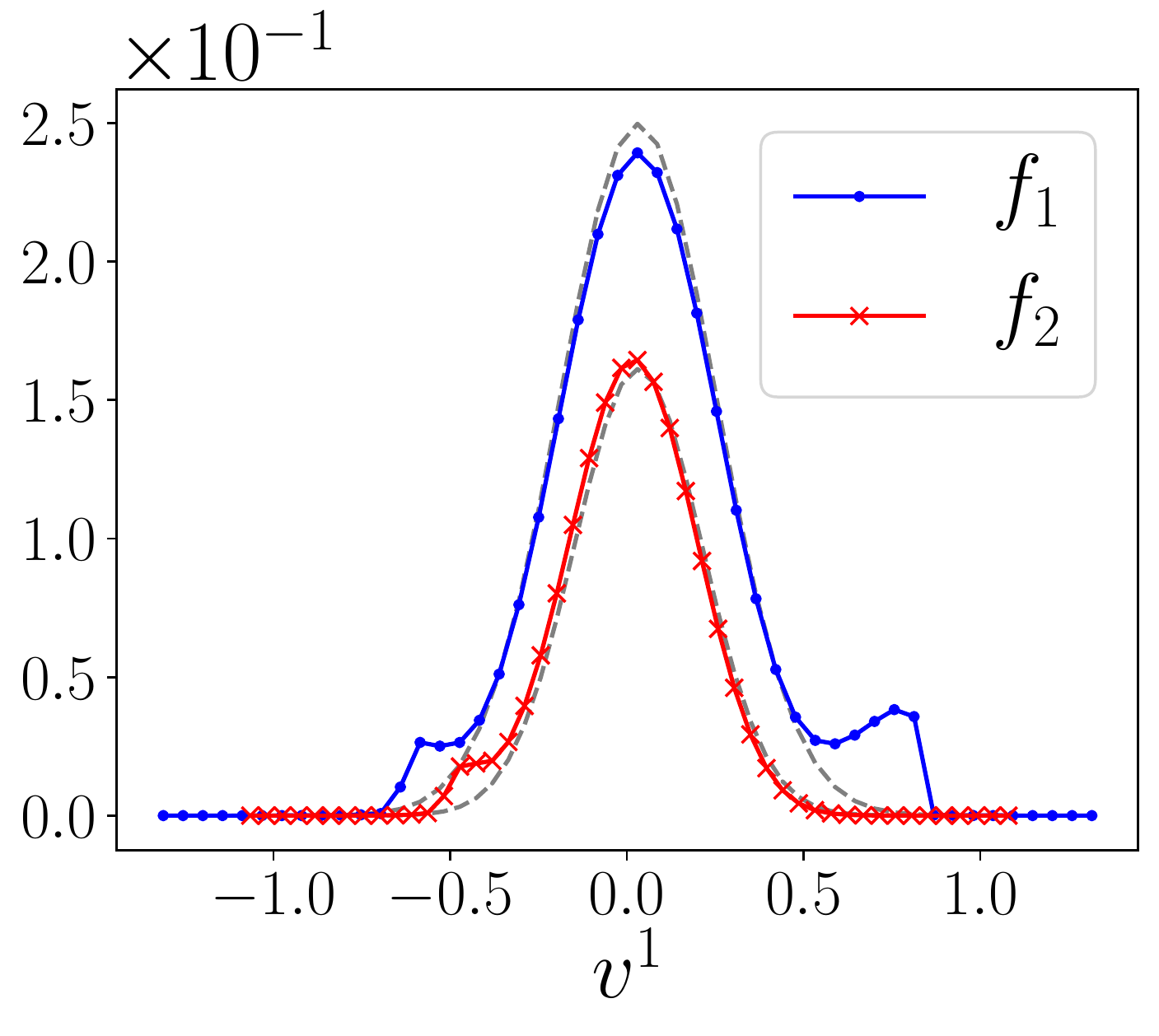}
\subcaption{$t=0.50$}
\end{subfigure}
\hfill
\begin{subfigure}[c]{0.32\textwidth}
\includegraphics[width=\textwidth]{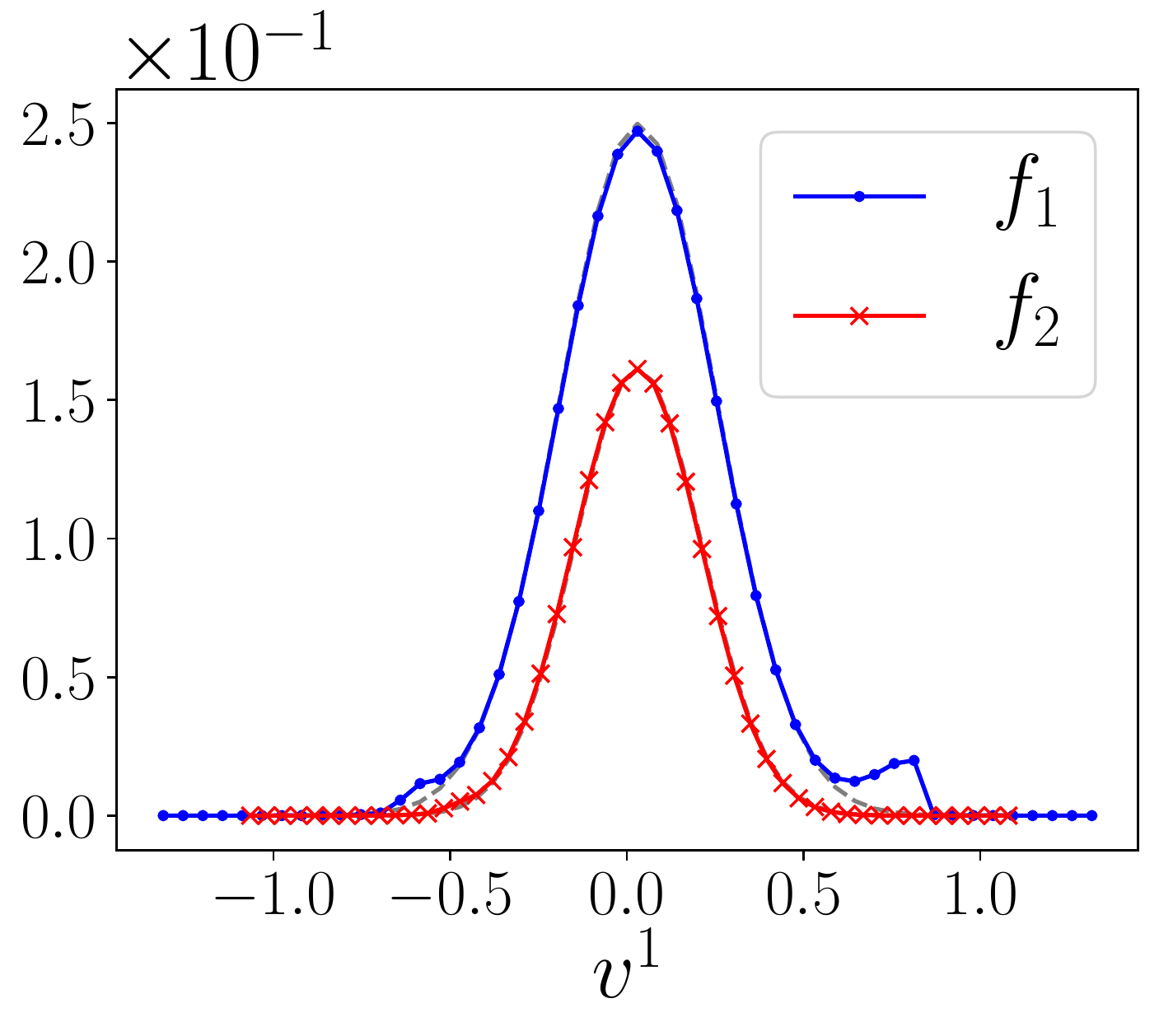}
\subcaption{$t=1.00$}
\end{subfigure}
\hfill
\begin{subfigure}[c]{0.32\textwidth}
\includegraphics[width=\textwidth]{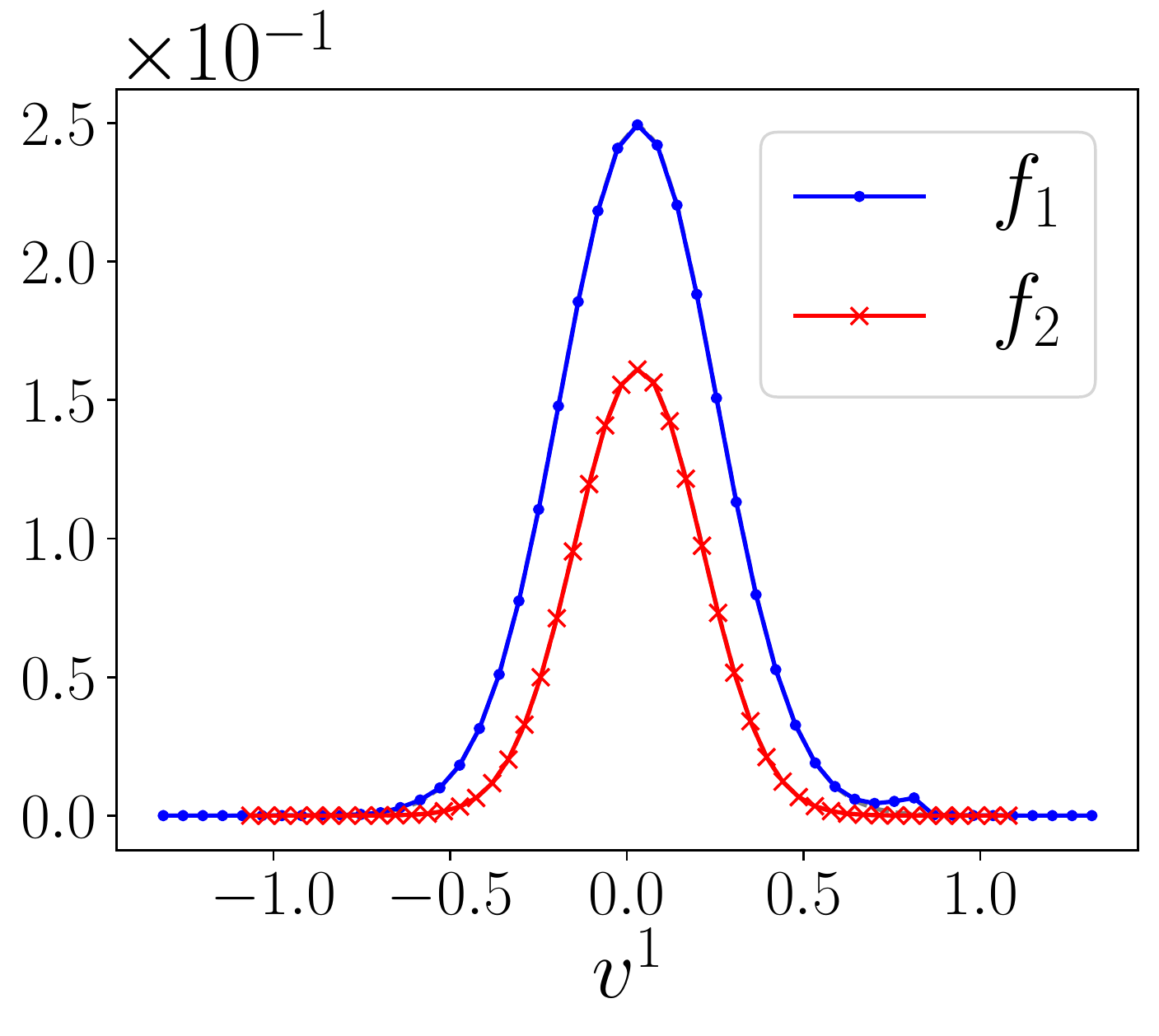}
\subcaption{$t=2.00$}
\end{subfigure}
\hfill
\begin{subfigure}[c]{0.32\textwidth}
\includegraphics[width=\textwidth]{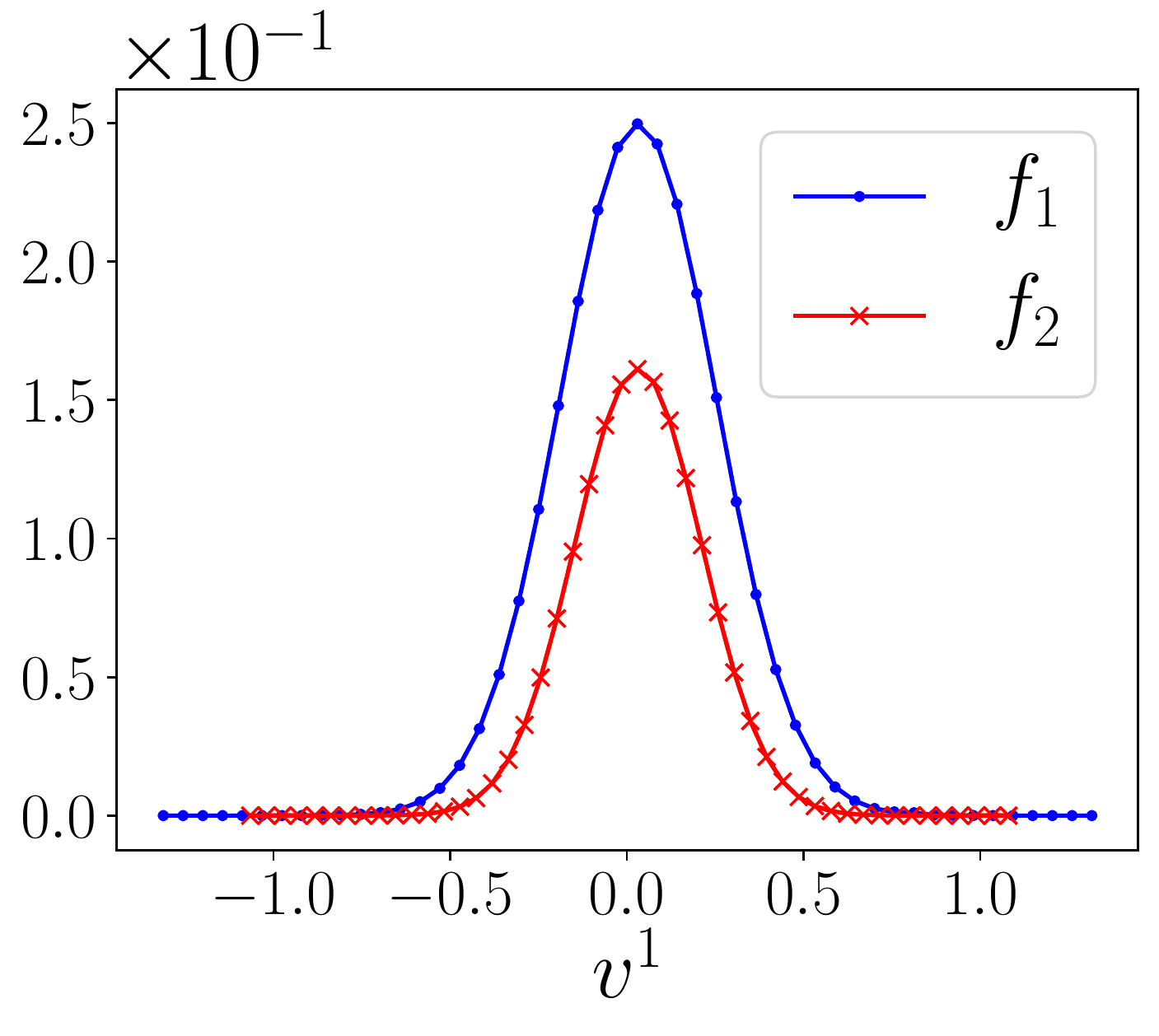}
\subcaption{$t=4.00$}
\end{subfigure}
    \caption{Relaxation of the distribution functions to Maxwellians for the test case in Section \ref{test:verification}. We fix $v^2=v^3=0$ and plot $f_i(v^1,v^2=0,v^3=0,t)$ at times $t$.  At time progress, the two distribution functions converge to Maxwellians centered around a common mean velocity with a width according to their common temperature divided by the respective mass.  For reference, these Maxwellians are shown by dotted gray lines.  }
    \label{fig:h-theorem}
\end{figure}

The purpose of this experiment is to illustrate basic properties of the BGK model.  We solve the spatially homogeneous version of \eqref{eq:BGK} for species with masses $m_1=1$ and $m_2=1.5$.
The initial distribution functions (see Figure \ref{fig:h-theorem_ic}) are given by
 \begin{align}
     f_i(\mbv,t=0) = 0.1\cdot m_i^{27} \cdot \exp\left(-\frac{0.01}{(0.75/m_i)^{10}-|\mbv-\mbu_i(0)|_1^{10}}\right),
 \end{align}
 with $\mbu_1(0)=(0.1,0,0)^\top$ and $\mbu_2(0)=(-0.1,0,0)^\top$.  The parameter choices here are not physical; rather they are chosen to yield an initial distribution with a particular shape that makes the relaxation easier to visualize.  With this initialization, the mixture mean velocity and mixture temperature have numerical values
 \begin{align}
     \mbu_\mix = 0.0322 \quad \text{and} \quad T_\mix=0.0487.
 \end{align}
 According to Proposition \ref{prop:mixture_quantities}, these values stay constant in time.   
 The collision frequencies take the form  
\begin{equation}
 \nu_{ij}(x,\mbv,t) = \frac{10\, n_j}{\delta_{ij} +\vert \mbv-\mbu_\mix \vert^{3} },
\end{equation}
with the regularization parameter $\delta_{ij} = 0.1 \cdot (\Delta v_{ij})^3$ where $\Delta v_{ij} = \frac{1}{4} \sqrt{T_\mix/(2\mu_{ij})}$ and $\mu_{ij} = m_im_j/(m_i+m_j)$.

The simulation is run using a velocity grid with $48^3$ nodes.
and the first-order temporal splitting scheme from Section \ref{subsec:firstordersplit} with time step $\Delta t = 0.01$. As demonstrated in Section \ref{sec:discreteproperties}, this scheme maintains positivity, conservation, and entropy dissipation properties of the BGK model.

In Figure 1, we plot the kinetic distributions at several different times and observe convergence to their respective equilibria. 
It is easy to see that the convergence to equilibrium is much faster in the center than near the tails of the distribution functions.
This is a consequence of the fact that the velocity-dependent collision frequency amplifies the relaxation process for small relative velocities. In Figure 2, we show convergence of the bulk velocities and temperatures to their equilibrium values, given by the mixture values in \eqref{eq:u_mixture} and \eqref{eq:T_mixture}. In Figure 3, we show the evolution of the entropy and the entropy dissipation.  As expected, the entropy decays monotonically.  In Figure 4, we demonstrate conservation properties.

\begin{figure}[htb]
\centering
\begin{subfigure}[c]{0.45\textwidth}
\includegraphics[width=\textwidth]{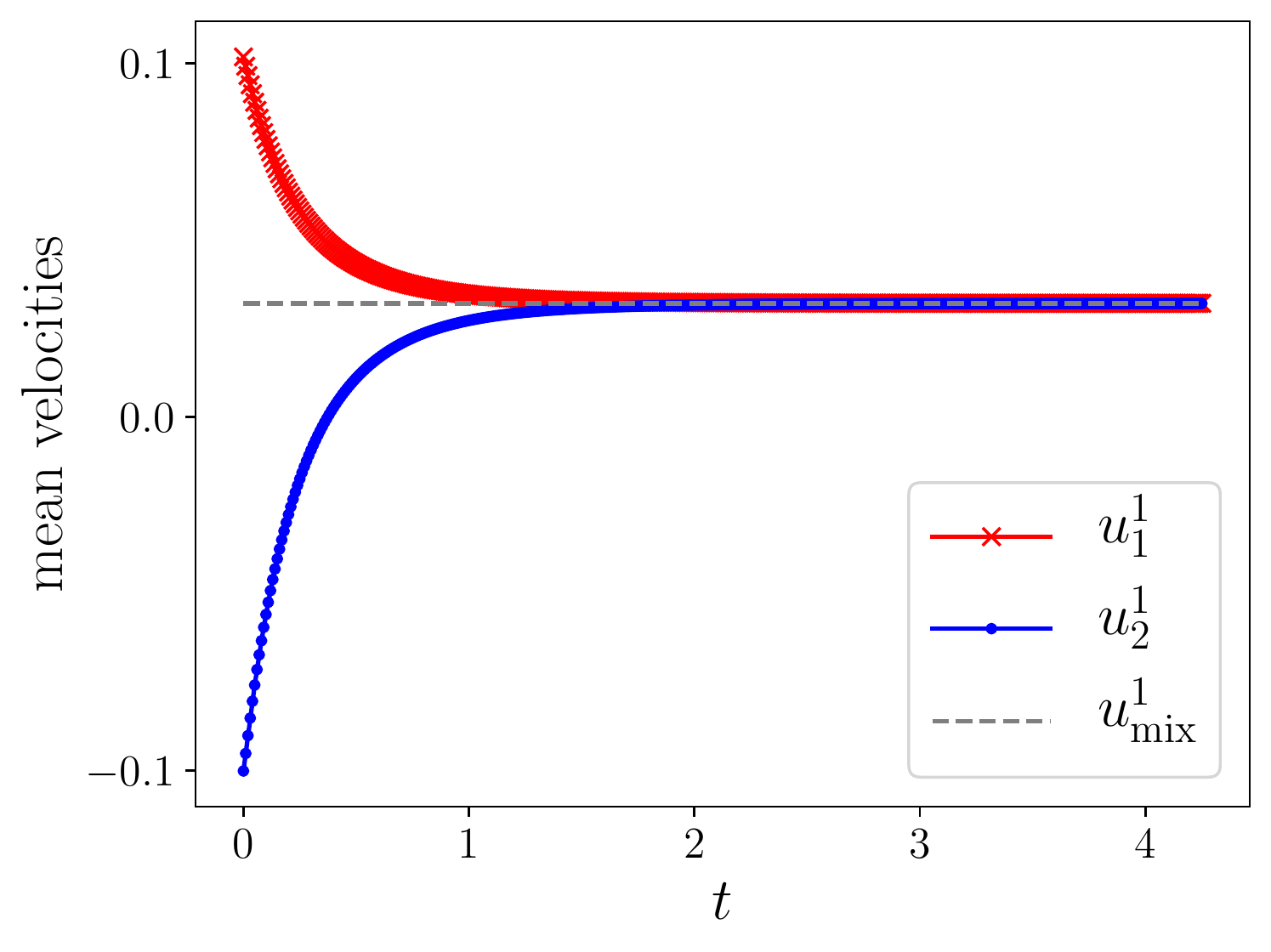}
\subcaption{convergence of mean velocities}
\end{subfigure}
\hfill
\begin{subfigure}[c]{0.45\textwidth}
\includegraphics[width=\textwidth]{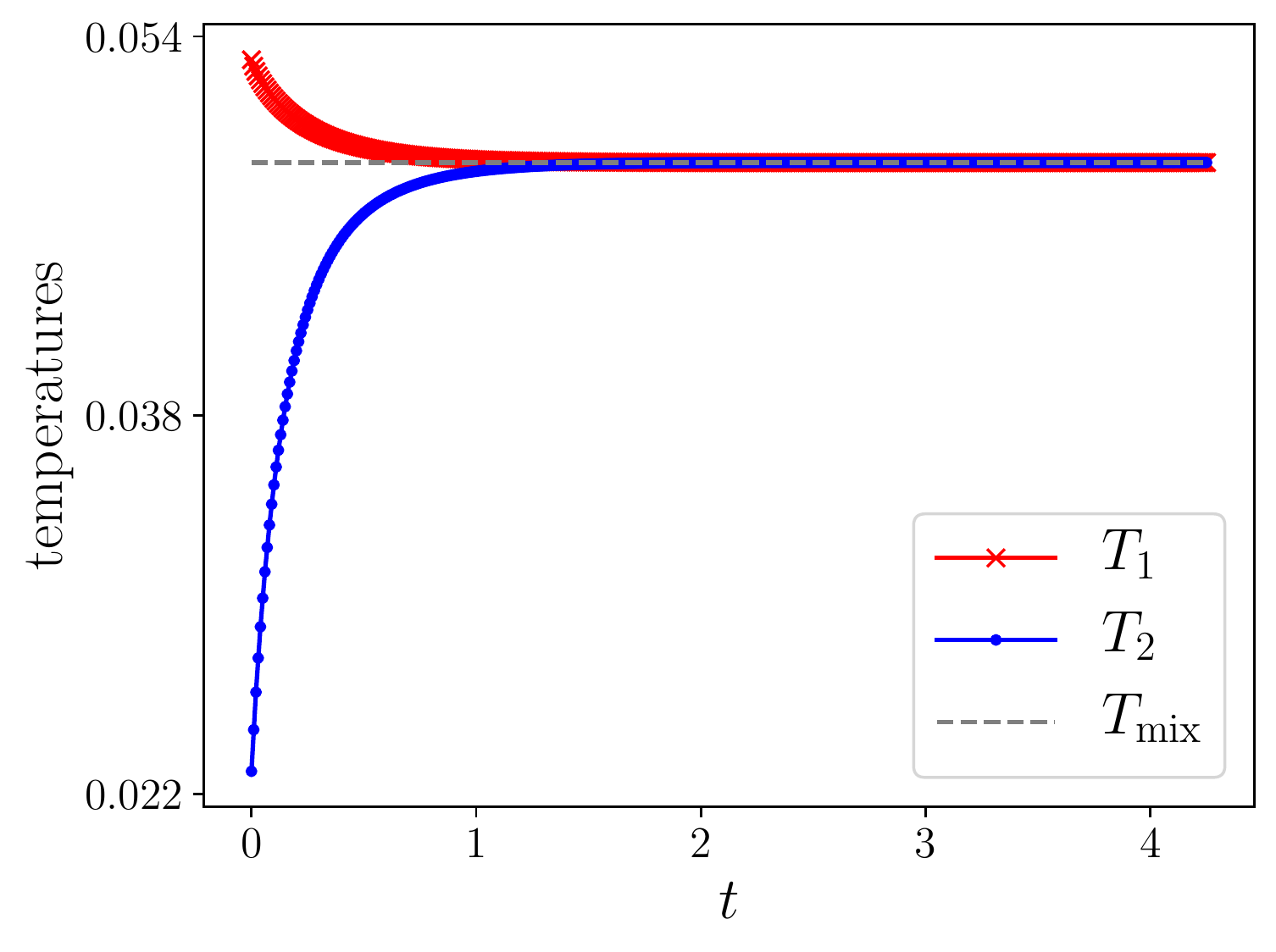}
\subcaption{convergence of temperatures}
\end{subfigure}
    \caption{Convergence of mean velocities and temperatures for the test case in Section \ref{test:verification}. In each plot, the dotted line denotes the mixture values, given in \eqref{eq:u_mixture} and \eqref{eq:T_mixture}. }
\end{figure}

\begin{figure}[htb]
\centering
\begin{subfigure}[c]{0.45\textwidth}
\includegraphics[width=\textwidth]{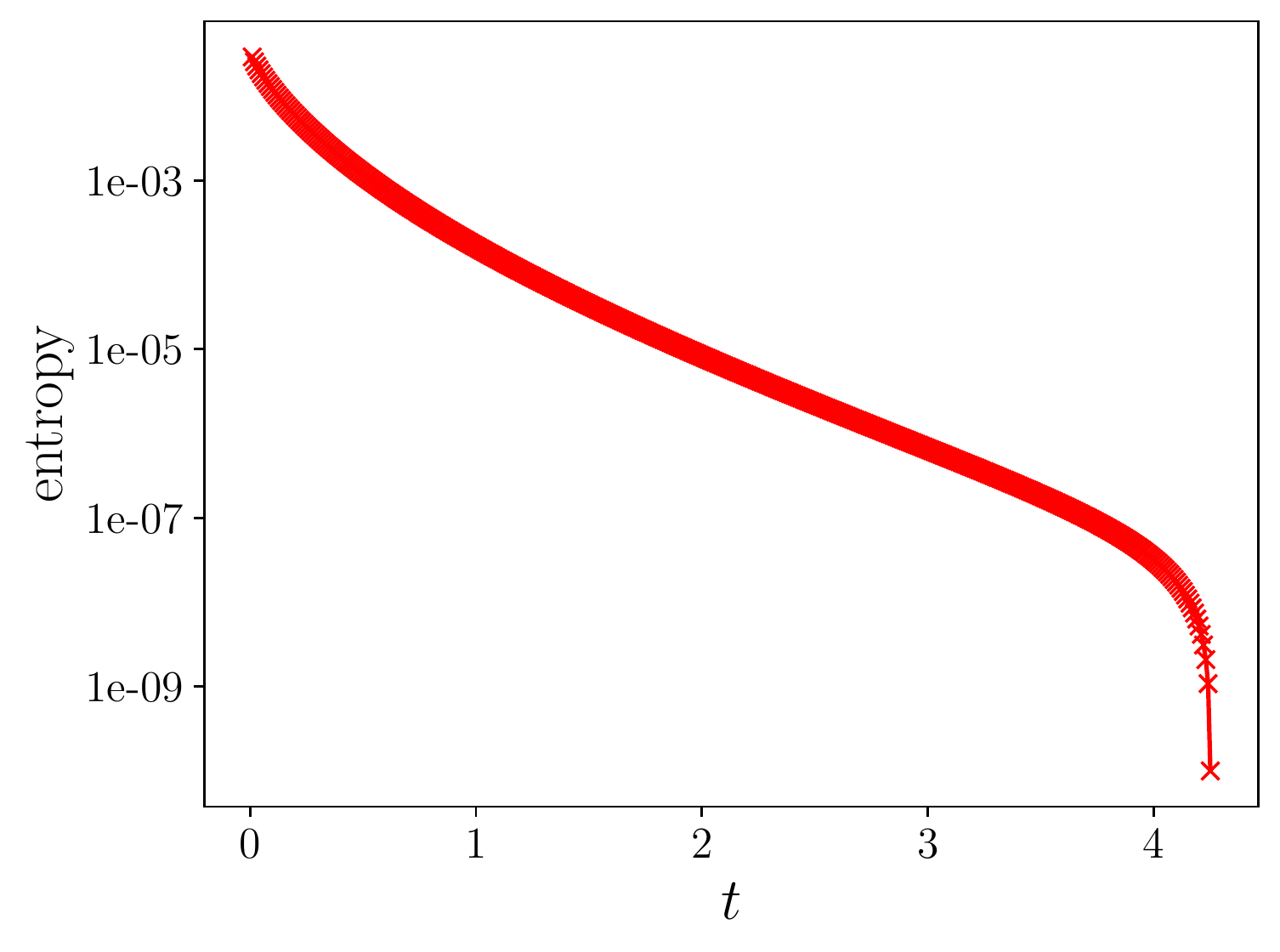}
\subcaption{entropy}
\end{subfigure}
\hfill
\begin{subfigure}[c]{0.45\textwidth}
\includegraphics[width=\textwidth]{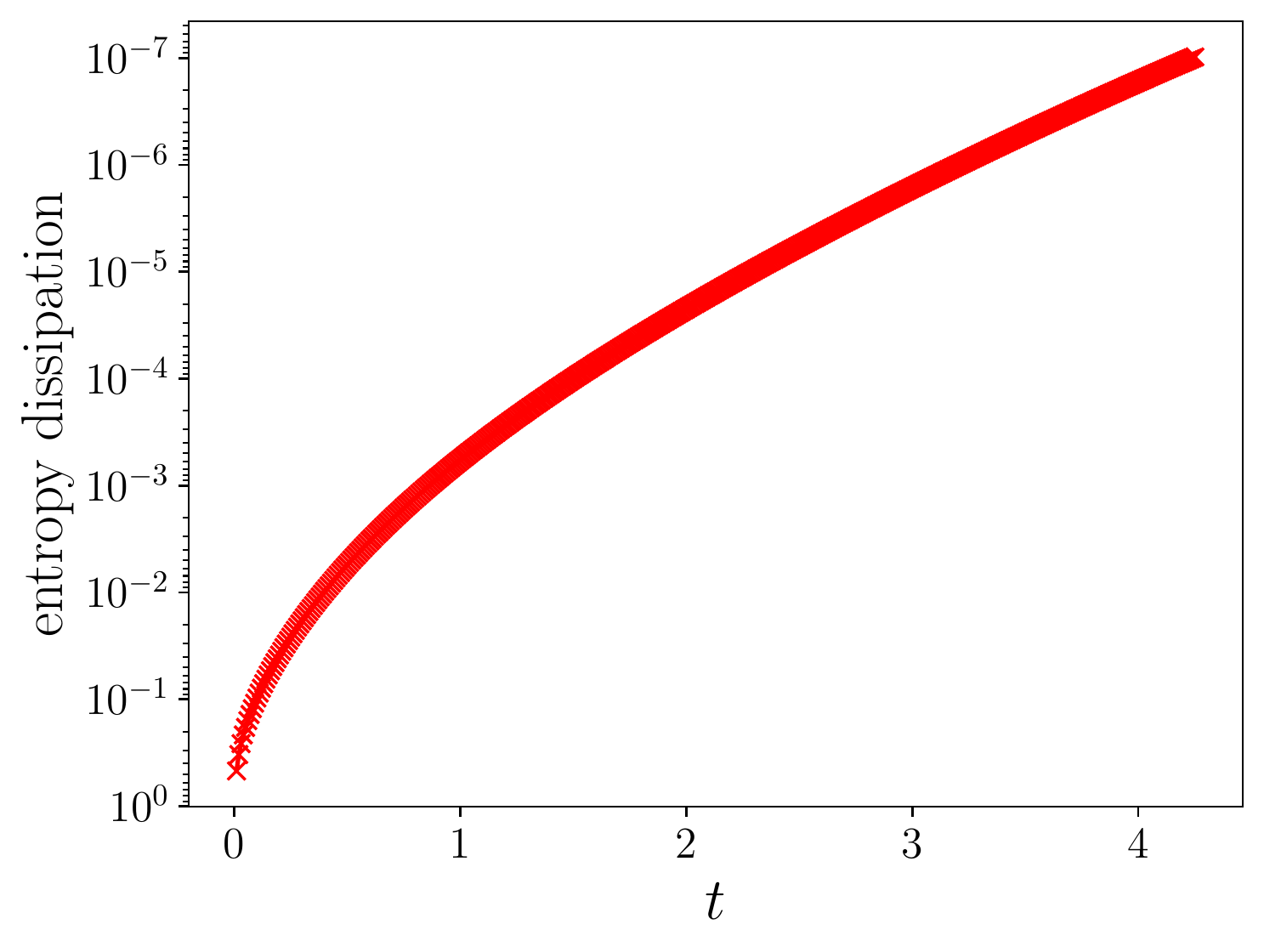}
\subcaption{entropy dissipation }
\end{subfigure}
    \caption{Entropy and entropy dissipation for the test case in Section \ref{test:verification}.  As predicted by the theory, the entropy decays monotonically.
    }
    \label{fig:entropy}
\end{figure}

\begin{figure}[htb]
\centering
\begin{subfigure}[c]{0.32\textwidth}
\includegraphics[width=\textwidth]{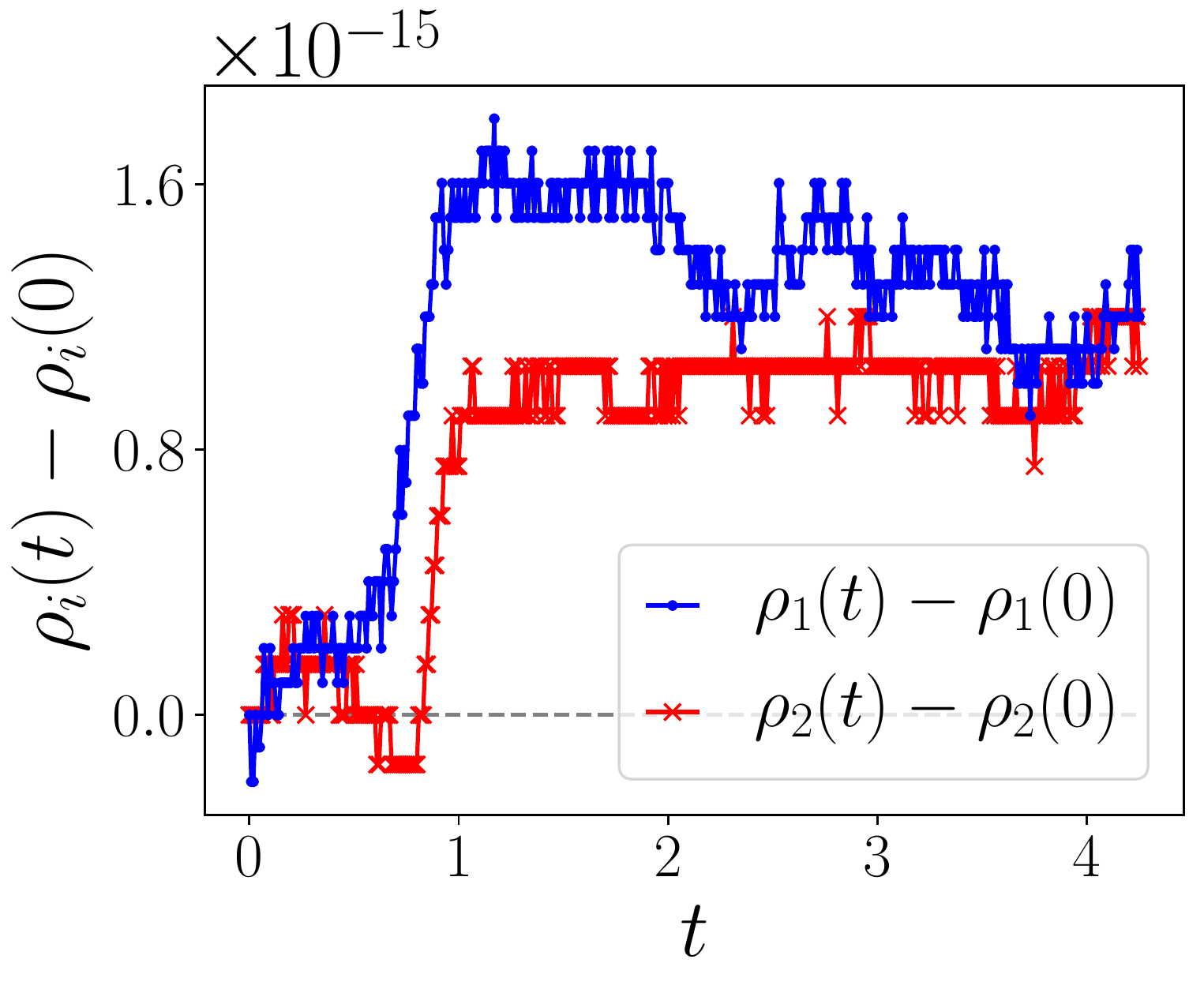}
\subcaption{species mass}
\end{subfigure}
\hfill
\begin{subfigure}[c]{0.32\textwidth}
\includegraphics[width=\textwidth]{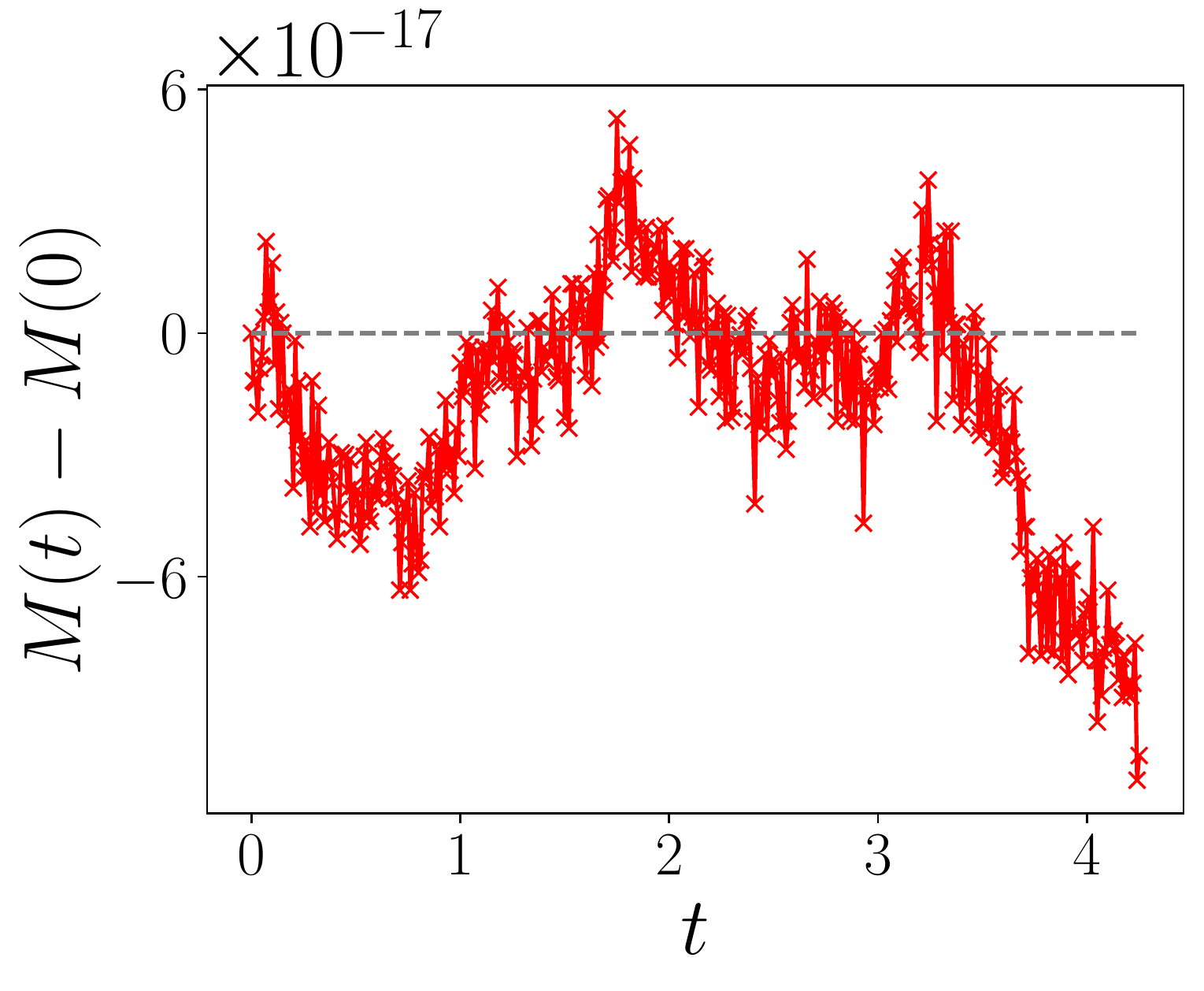}
\subcaption{total momentum}
\end{subfigure}
\hfill
\begin{subfigure}[c]{0.32\textwidth}
\includegraphics[width=\textwidth]{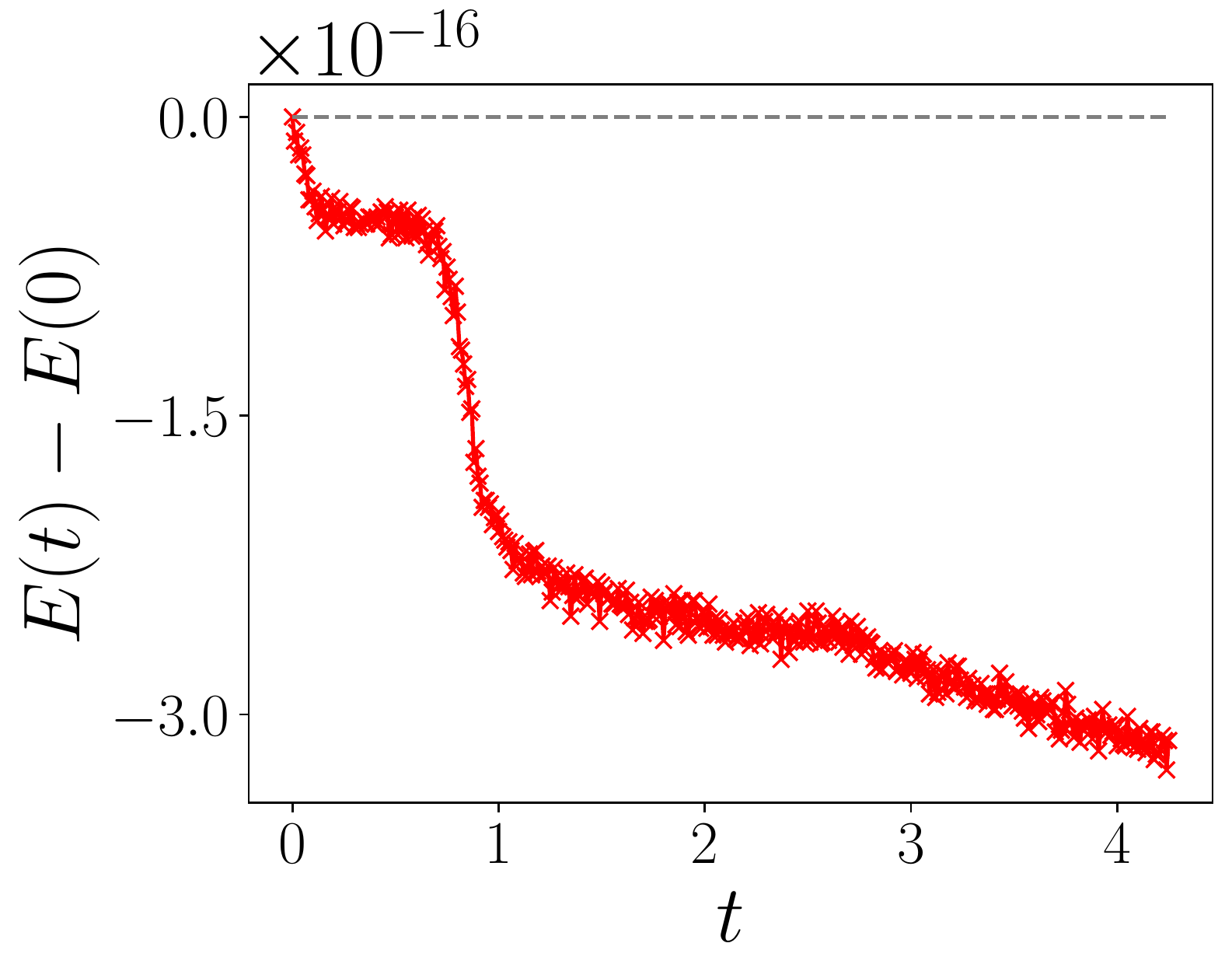}
\subcaption{total energy}
\end{subfigure}
    \caption{Global conservation properties for the test case in Section \ref{test:verification}.  The mass densities of each species, 
    the total momentum  ($M$) and total energy ($E$) have small oscillations on the order of $10^{-15}$ or less.  }
    \label{fig:conservation}
\end{figure}

\subsubsection{Hydrogen-Carbon test case } \label{test:HC} 
In this test case, we explore the effects of the velocity-dependent frequencies on the relaxation behavior of a multi-species problem in a more physically relevant setting, with dimensional formulas given in the cgs unit system.  

To define the collision frequency $\nu_{ij}$, we build on the formulas given in \cite{LeeMore}.  A simple model for the collision frequency is given by 
\begin{equation}\label{eq:nuMT}
    \nu_{ij}(\mbv) = n_j |\mbv - \mbu_\mix| \, \sigma_{\rm{mt}}(|\mbv - \mbu_\mix|),
\end{equation}
where $\sigma_{\rm{mt}}$ is the momentum transfer cross section for Coulomb collisions:
\begin{equation}
    \sigma_{\rm{mt}}(|\mbv - \mbu_\mix|) = 4\pi \left(\frac{Z_i Z_j e^2}{2\mu_{ij} |\mbv - \mbu_\mix|^2}\right)^2 L_{ij}(Z_i,Z_j,n_1, n_2, T_\mix).
\end{equation}
Here $\mu_{ij} = m_i m_j/(m_i + m_j)$ is the reduced mass; $Z_i e$ and  $Z_j e$ are the charges of the species $i$ and $j$ particles, respectively; and $L_{ij}$ is the Coulomb logarithm:
\begin{align}
\label{eq:Coulomb_log}
    L_{ij} = \frac{1}{2}\log\left(1+\frac{\lambda_D^2}{b_{90,ij}^2}\right),
\end{align}
where $b_{90,ij}$ is the distance of closest approach:
\begin{align}
    b_{90,ij} = \frac{Z_i Z_j e^2}{T_\mix},
\end{align}
where $e^2 = 1.44 \times 10^{-7}$ eV$\cdot$cm, in cgs units. For the Debye length $\lambda_D$ in \eqref{eq:Coulomb_log}, we use the following formulae
\begin{align}
    \lambda_D &= \left(\frac{1}{\lambda_e^2} + \frac{1}{\lambda_I^2}\right)^{-1/2} \quad \text{with} \quad
    \lambda_e = \left(\frac{T_\mix}{4 \pi n_e e^2}\right)^{1/2} \quad \text{and} \quad
    \lambda_I = \left( \frac{1}{\lambda_1^2}+\frac{1}{\lambda_2^2}\right)^{-1/2}, \quad \\\text{where}\quad
    \lambda_i &= \left(\frac{T_\mix}{4 \pi n_i Z_i^2 e^2}\right)^{1/2} \quad \text{and} \quad
    n_e =  Z_1 n_1 + Z_2 n_2.
\end{align}

For the purposes of evaluating the BGK model in this paper,  \eqref{eq:nuMT}:
\begin{align} \label{eq:colfreq_dep}
    \nu_{ij}(\mbv) =  4\pi n_j \left(\frac{Z_i Z_j e^2}{2\mu_{ij}} \right)^2 \left(\frac{1}{\delta_{ij} + |\mbv - \mbu_\mix|^3}\right) L_{ij}(Z_i,Z_j,n_1, n_2, T_\mix),
\end{align}
i.e, we add a small regularization parameter $\delta_{ij}>0$ in the denominator of \eqref{eq:colfreq_dep} to avoid a singularity at zero relative velocity. 
For the numerical experiments, one needs to ensure that $\delta_{ij}$ is much smaller than $|\mbv-\mbu_{\mix}|^3$, and thus we set $\delta_{ij} = 0.1 \cdot (\Delta v_{ij})^3$, where $\Delta v_{ij} = \frac{1}{4} \sqrt{k_BT_\mix/(2\mu_{ij})}$ and $k_B = 1.602\cdot 10^{-12}\,\text{erg}/\text{eV}$ is Boltzmann's constant in cgs units. This choice ensures the symmetry $n_1 \nu_{12} = n_2 \nu_{21}$.
The mixture quantities $\mbu_\mix$ and $T_\mix$ defined in \eqref{eq:u_mixture} and \eqref{eq:T_mixture} are inserted into these formulas to determine the collision frequencies used in the model. 

\begin{figure}[htb]
    \begin{subfigure}[c]{0.45\textwidth}
    \includegraphics[width=\textwidth]{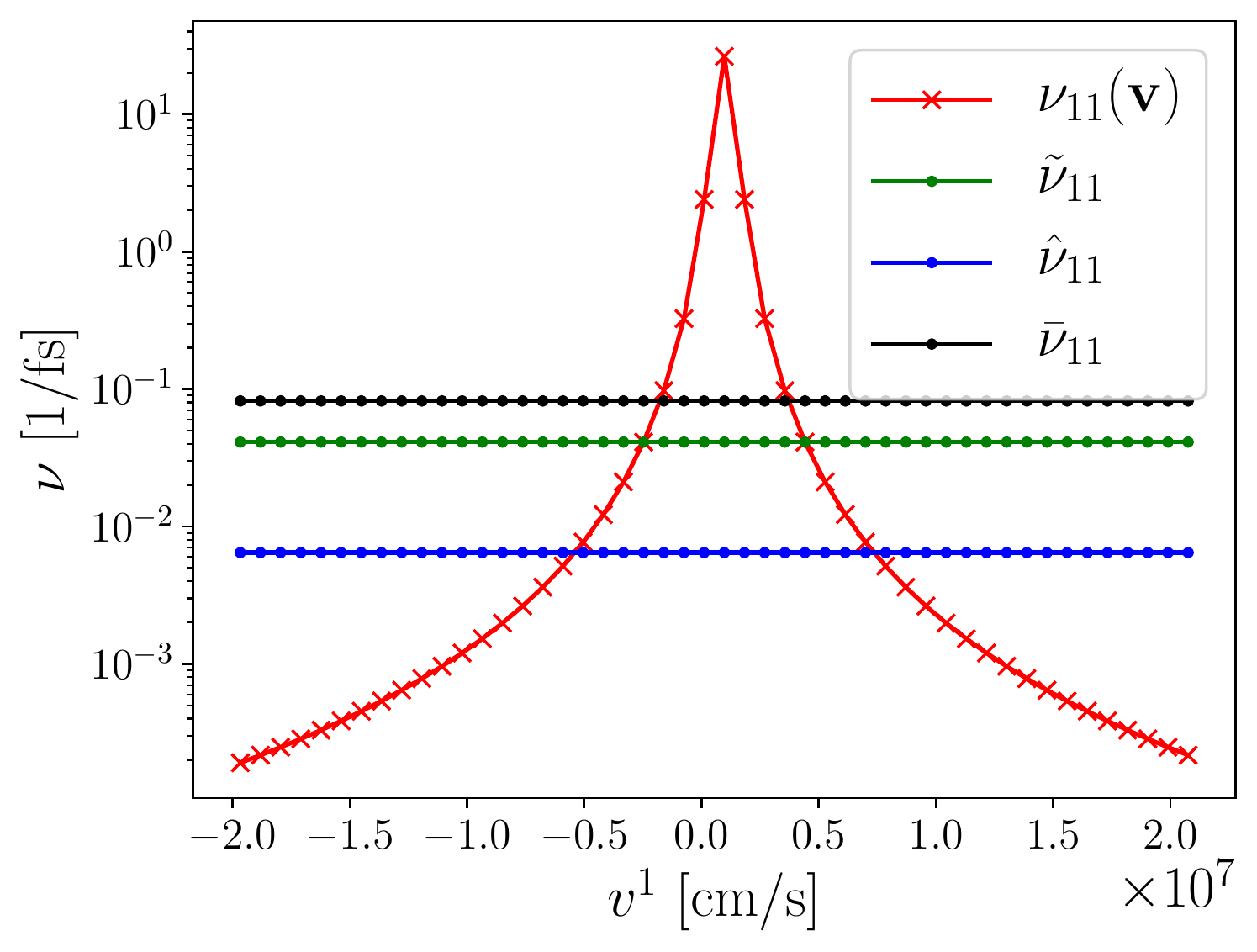}
    \end{subfigure}
    \hfill
    \begin{subfigure}[c]{0.45\textwidth}
    \includegraphics[width=\textwidth]{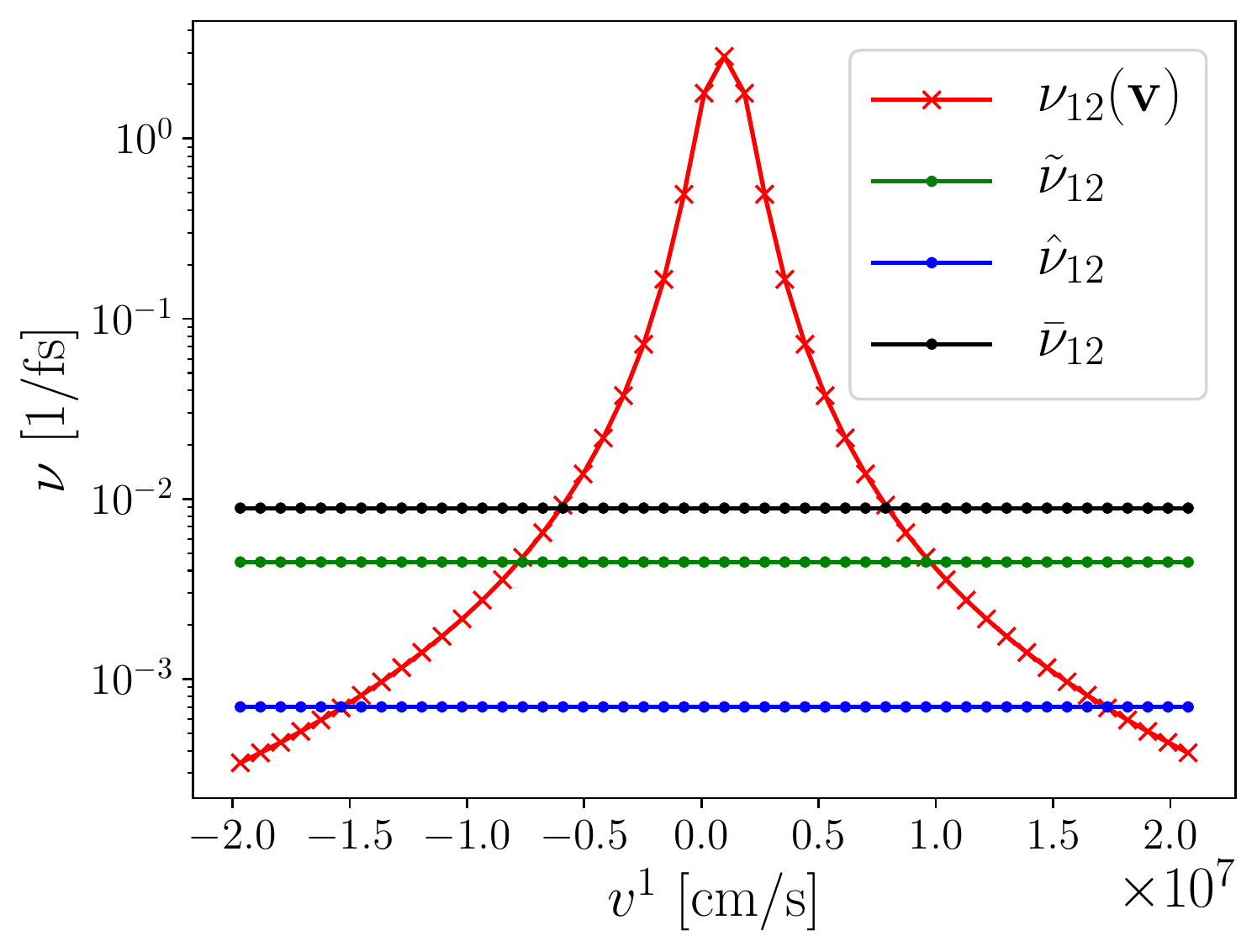}
    \end{subfigure} \\
    \begin{subfigure}[c]{0.45\textwidth}
    \includegraphics[width=\textwidth]{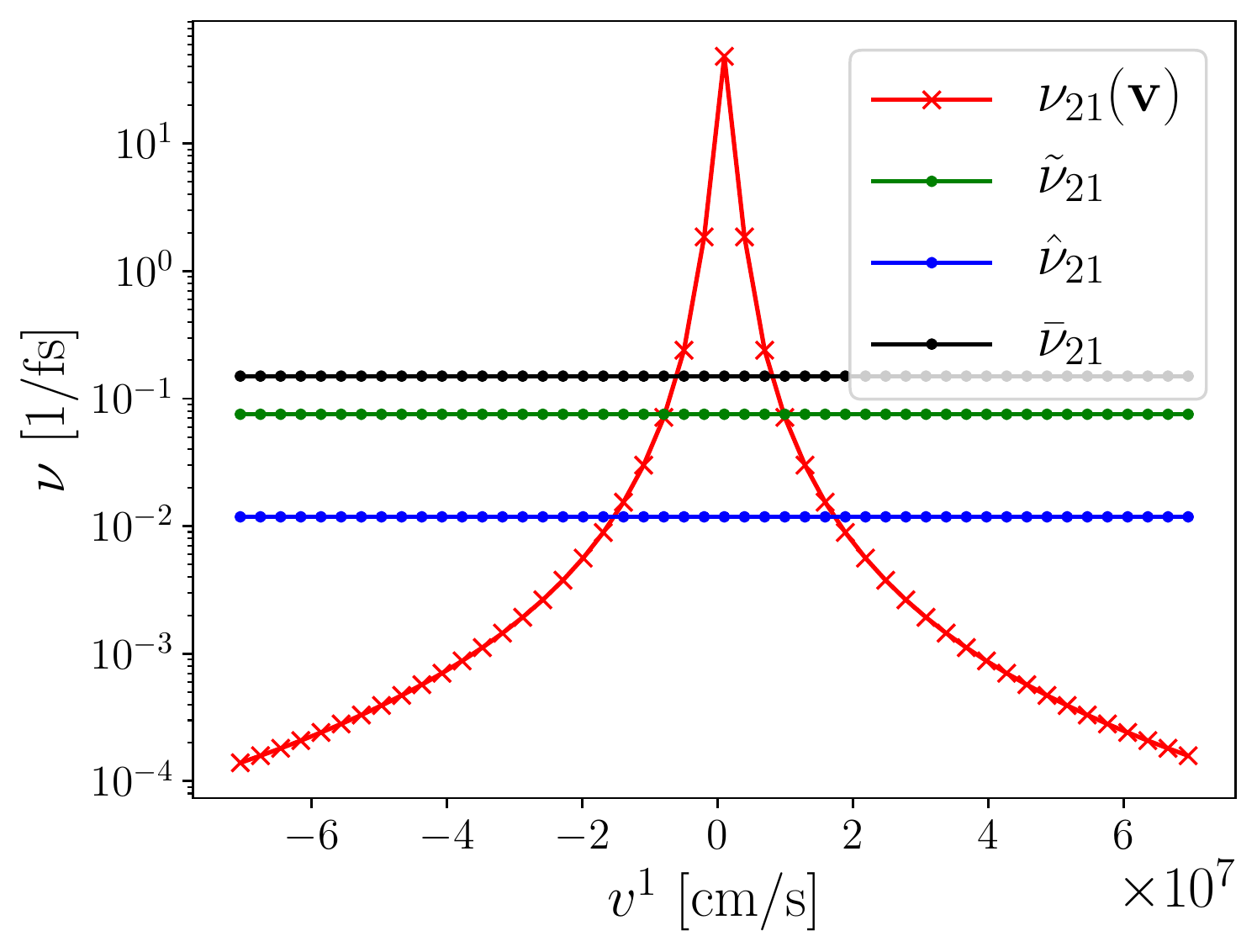}
    \end{subfigure}
    \hfill
    \begin{subfigure}[c]{0.45\textwidth}
    \includegraphics[width=\textwidth]{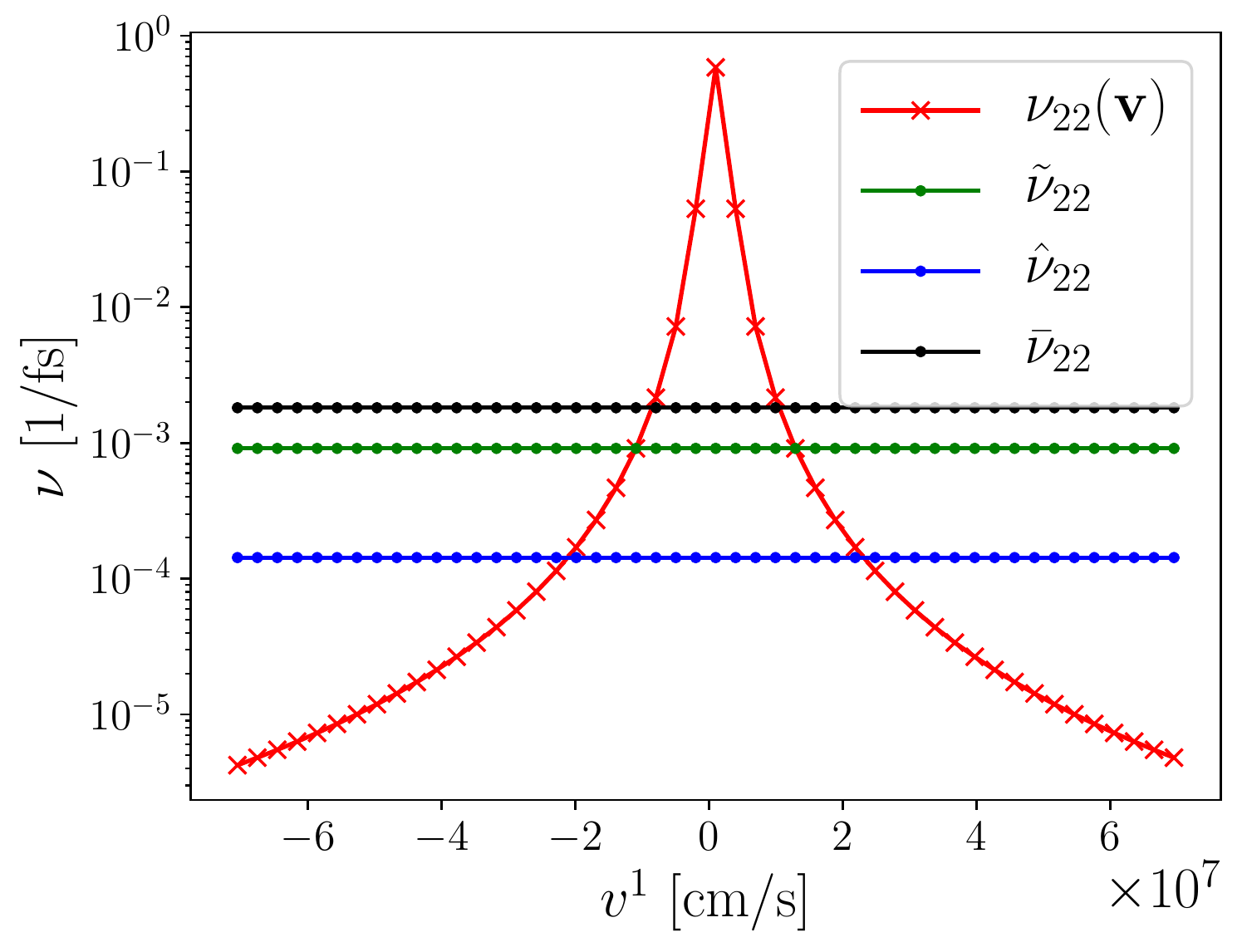}
    \end{subfigure}
    \caption{Collision frequencies given in \eqref{eq:colfreq_dep}, \eqref{eq:colfreq_indep_1}, \eqref{eq:colfreq_indep_2} and \eqref{eq:colfreq_indep_3} along the line $v^2 = v^3 = 0$. The large constant values for $\bar{\nu}$ correspond to the fastest relaxation process, see Figure \ref{fig:HC}.}
    \label{fig:HC_collision_frequencies}
\end{figure}

For comparison, we consider three velocity-independent collision frequencies that are often used as simpler alternatives to \eqref{eq:colfreq_dep}:
\begin{enumerate}
    \item Replacing $\vert\mbv-\mbu_\mix\vert$ by the thermal velocity $v_T =\sqrt{k_B T_\mix/(2\mu_{ij})}$ gives
\begin{align} \label{eq:colfreq_indep_1}
    \tilde{\nu}_{ij} = 4\pi n_j \left(\frac{Z_i Z_j e^2}{ 2\mu_{ij}} \right)^2 \left(\frac{1}{\delta_{ij} + v_T^3}\right) L_{ij}.
\end{align}
    \item Replacing $\vert\mbv-\mbu_\mix\vert^3$ by the weighted average
\begin{align} \label{eq:v_hat}
    \hat{v}^3=\frac{\int |\mbv -\mbu_\mix |^3 \mathcal{M}(\mbv)\dv}{\int \mathcal{M}(\mbv)\dv},
\end{align}
where 
\begin{align}
    \mathcal{M}(\mbv) = n_i \left( \frac{\mu_{ij}}{\pi T}\right)^{3/2} \exp\left(-\frac{\mu_{ij}|\mbv-\mbu_\mix|^2}{T_\mix}\right),
\end{align}
gives 
\begin{align}   \label{eq:colfreq_indep_2}
    \hat{\nu}_{ij} = 4\pi n_j \left(\frac{Z_i Z_j e^2}{2 \mu_{ij}} \right)^2 \left(\frac{1}{\delta_{ij} + \hat{v}^3}\right) L_{ij}.
\end{align}
\item Computing a weighted average of $\nu_{ij}$ directly gives
\begin{align} \label{eq:colfreq_indep_3}
    \bar{\nu}_{ij} = \frac{\int \nu_{ij}(\mbv) \mathcal M(\mbv)\dv}{\int \mathcal M(\mbv)\dv}.
\end{align}
\end{enumerate}
While the first option above is convenient and more common in applications \cite{SM}, it is somewhat arbitrary. The second and third options, on the other hand, provide a more consistent normalization. According to  Proposition \ref{prop:mixture_quantities}, the collision frequencies stay constant in time because the problem is spatially homogeneous. For purposes of illustration, we plot them in Figure \ref{fig:HC_collision_frequencies}.  

We consider relaxation between carbon (species 1) and hydrogen (species 2), with masses and charge numbers
\begin{alignat}{3}
m_1 &= 1.993\cdot 10^{-23} \,\text{g} , 
\qquad &&m_2 = 1.661\cdot 10^{-24} \,\text{g} , \\
Z_1&=6, \; \quad &&Z_2=1. \; \notag
\end{alignat}
Initially, the distribution functions are Maxwellians: $f_i=M_i[n_i,\mbu_i,T_i]$ with
\begin{alignat}{3}
n_1 &= 6.1\cdot 10^{22} \,\text{cm}^{-3} , \qquad  &&n_2 = 3.6133\cdot 10^{21} \,\text{cm}^{-3}  , \\
\mbu_1 &= (9.818\cdot 10^5, \, 0, \, 0)^\top\,\frac{\text{cm}}{\text{s}} , \qquad 
&&\mbu_2 = (0 , \, 0, \, 0)^\top \,\frac{\text{cm}}{\text{s}}, \notag\\
T_1 &= 150 \,\text{eV} , \quad &&T_2 = 100 \,\text{eV} .\notag
\end{alignat}

We simulate this test case using a velocity grid with $48^3$ nodes and the second-order IMEX Runge-Kutta scheme from Section \ref{subsec:secondorderIMEX} with time step $\Delta t = 0.8\, $fs.

\begin{figure}[htb]
    \begin{subfigure}[c]{0.45\textwidth}
    \includegraphics[width=\textwidth]{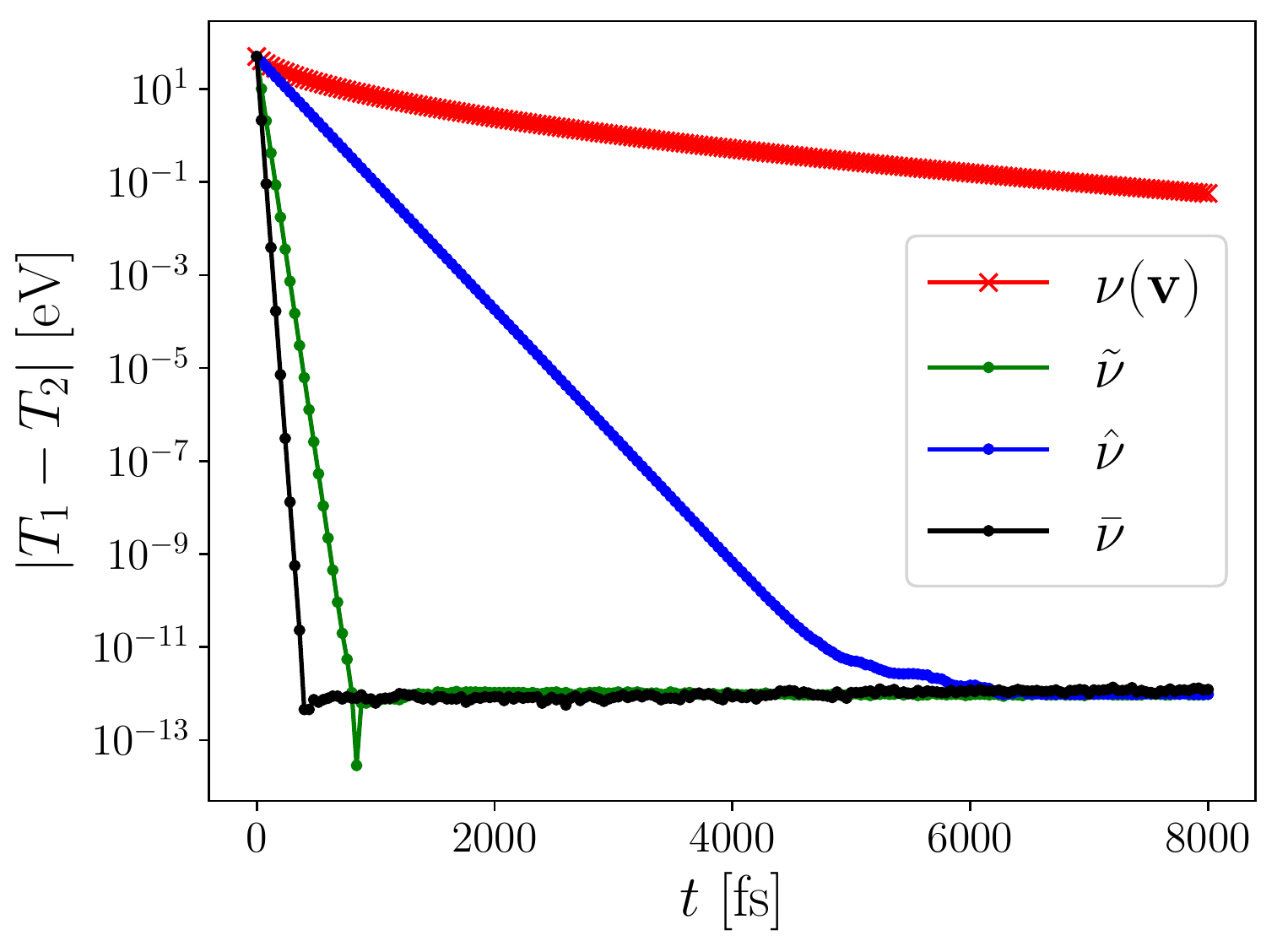}
    \subcaption{temperature}
    \end{subfigure}
    \hfill
    \begin{subfigure}[c]{0.45\textwidth}
    \includegraphics[width=\textwidth]{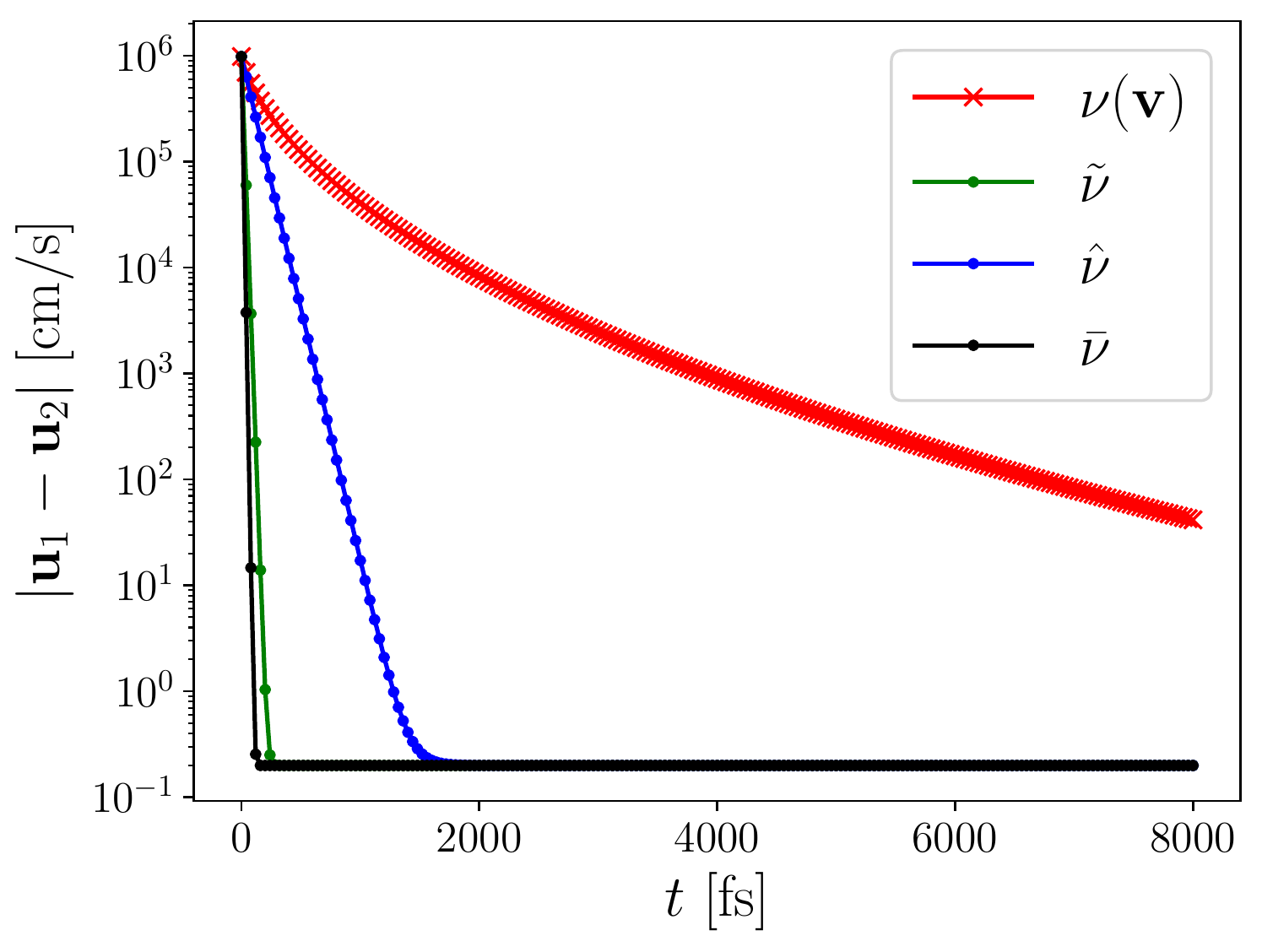}
    \subcaption{ mean velocity}
    \end{subfigure}
    \caption{Evolution of the difference in species temperatures and mean velocities for the Hydrogen-Carbon test case in Section \ref{test:HC}. The convergence for all velocity-independent collision frequencies---$\tilde{\nu}$ in \eqref{eq:colfreq_indep_1}, $\hat{\nu}$ in \eqref{eq:colfreq_indep_2}, and $\bar{\nu}$ in \eqref{eq:colfreq_indep_3}---appears exponential.  However, the convergence for velocity-dependent collision frequency  $\nu$ given in \eqref{eq:colfreq_dep} is significantly longer and notably different.}
    \label{fig:HC}
\end{figure}

\begin{figure}[htb]
    \begin{subfigure}[c]{0.45\textwidth}
    \includegraphics[width=\textwidth]{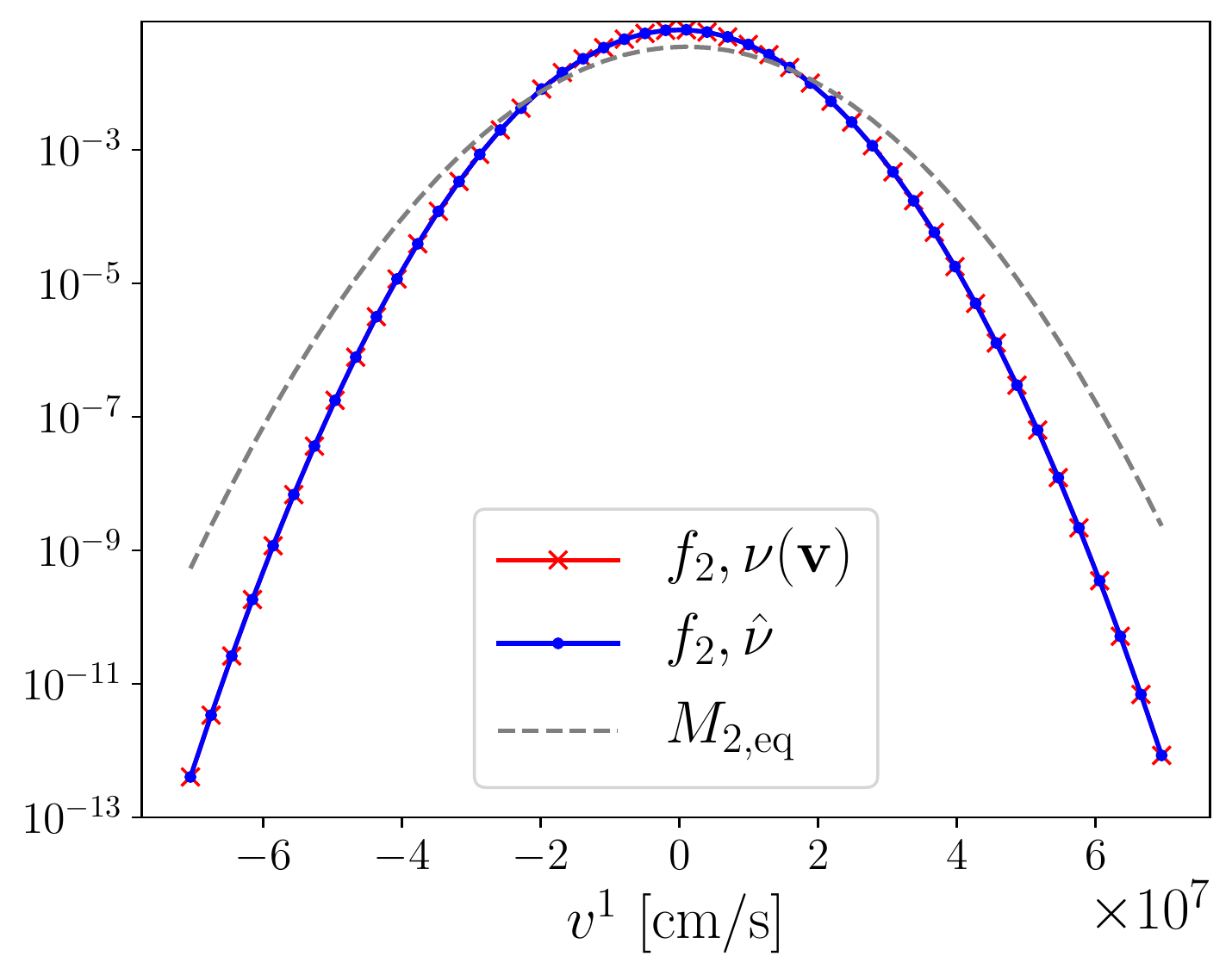}
    \subcaption{$t=0$ }
    \end{subfigure}
    \hfill
    \begin{subfigure}[c]{0.45\textwidth}
    \includegraphics[width=\textwidth]{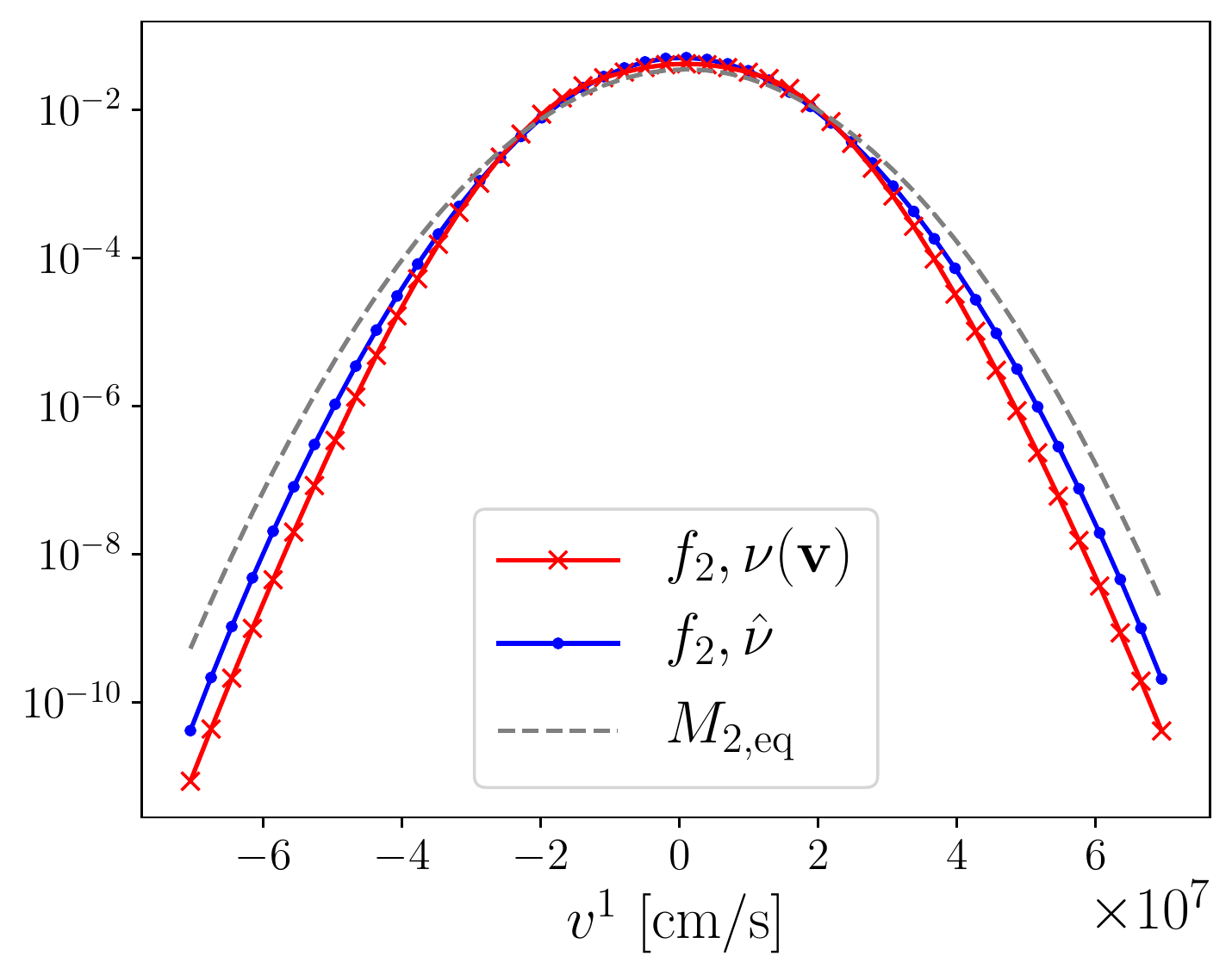}
    \subcaption{$t=100\Delta t$ }
    \end{subfigure} \\
    \begin{subfigure}[c]{0.45\textwidth}
    \includegraphics[width=\textwidth]{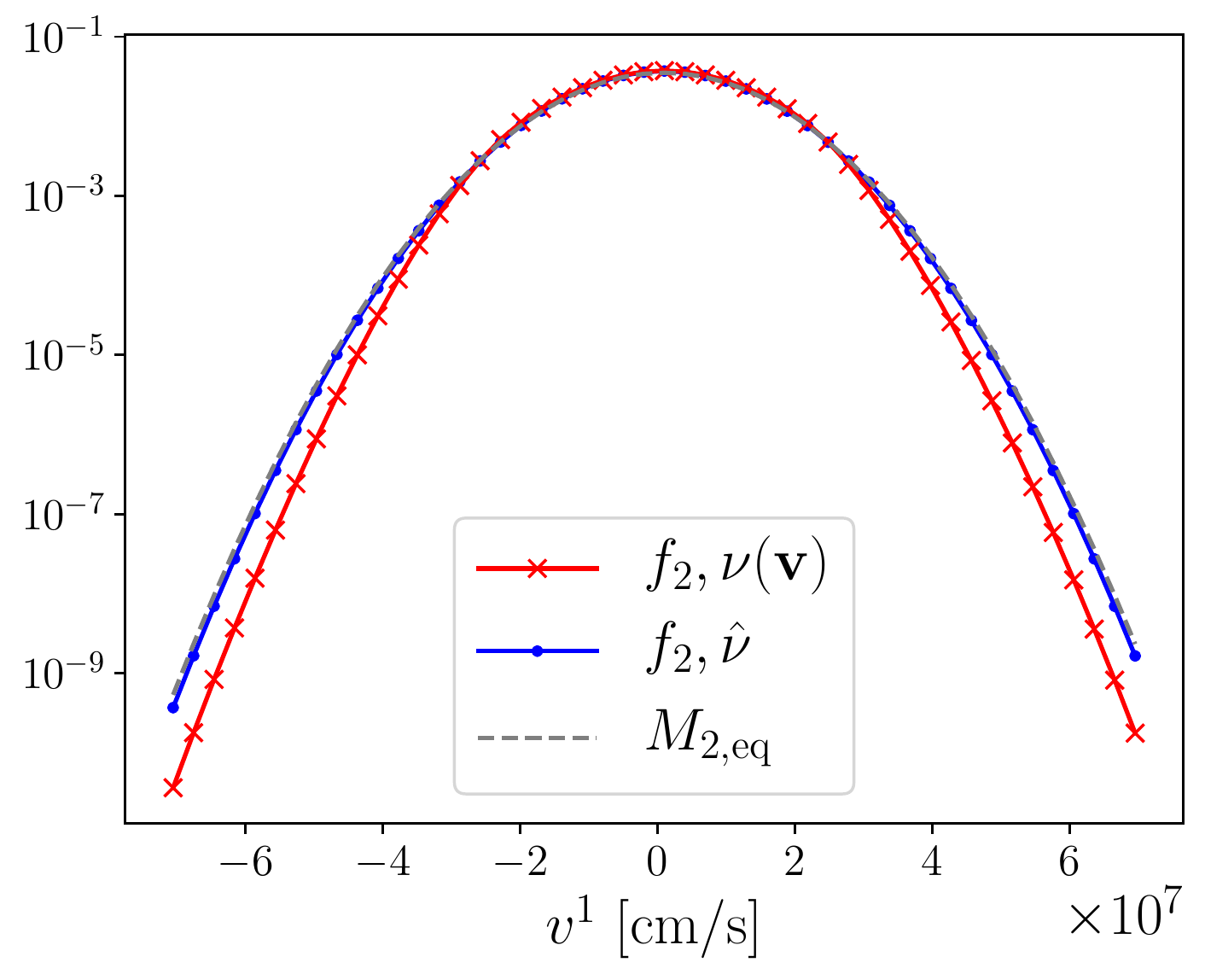}
    \subcaption{$t=500\Delta t$ }
    \end{subfigure}
    \hfill
    \begin{subfigure}[c]{0.45\textwidth}
    \includegraphics[width=\textwidth]{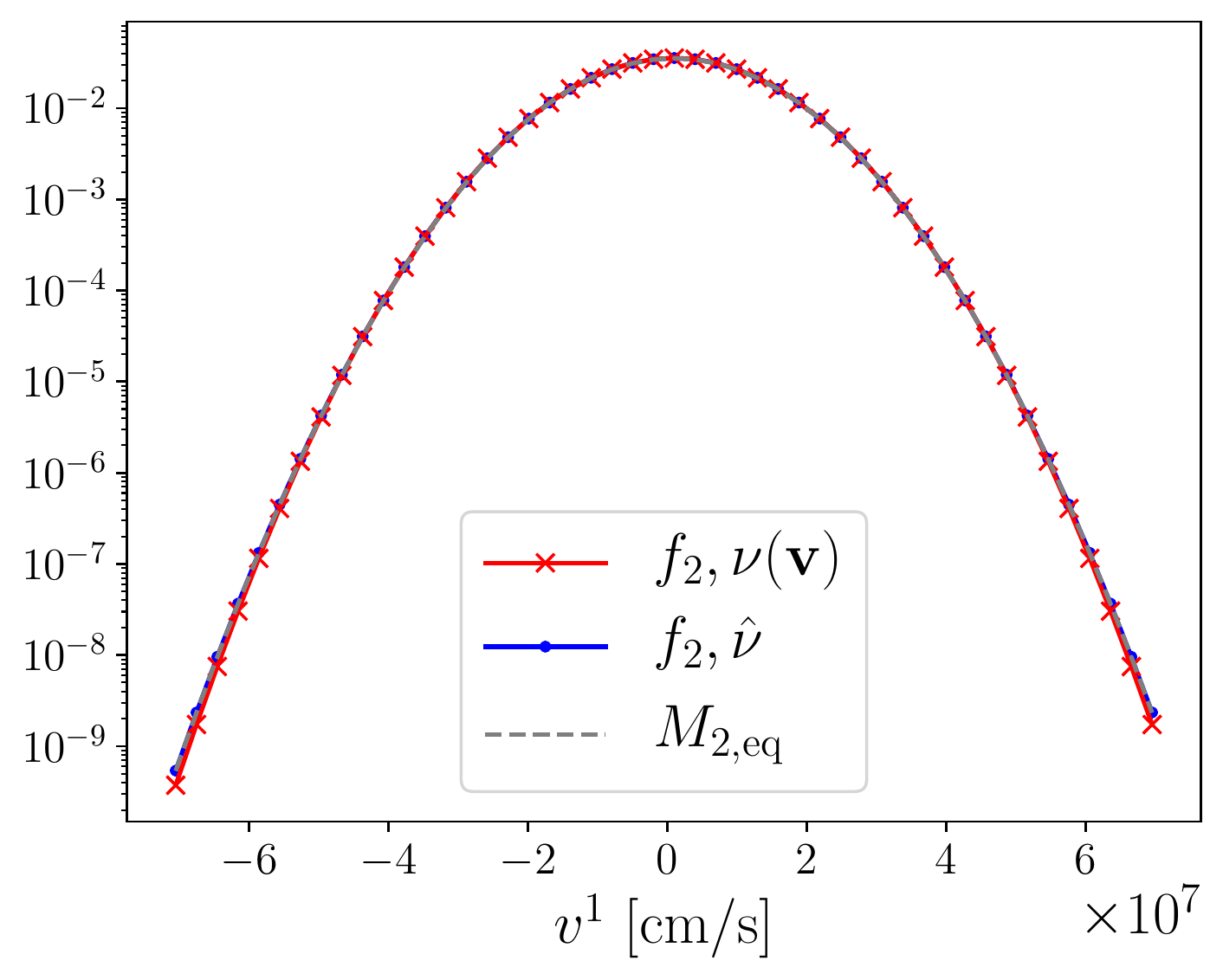}
    \subcaption{$t=10000\Delta t$}
    \end{subfigure}
    \caption{Relaxation of the kinetic distribution function of hydrogen in the Hydrogen-Carbon test case of Section \ref{test:HC} with time step $\Delta t = 0.8\, $fs.  The velocity components $v^2$ and $v^3$ are fixed at $u^2=0$ and $u^3=0$, respectively. The dashed line corresponds to the {Maxwellian $M_{2,\rm eq}=M_2[n_2,\mbu_\mix,T_\mix]$} for the hydrogen species.  The blue and red lines are results computed with the velocity-independent collision frequencies $\hat{\nu}$ in \eqref{eq:colfreq_indep_2} and the velocity-dependent collision frequencies $\nu(\mbv)$ in \eqref{eq:colfreq_dep}, respectively. The tails of the distribution converge more slowly for the velocity-dependent collision frequencies. }
    \label{fig:HC_distributions_center}
\end{figure}

In Figure~\ref{fig:HC}, we plot the evolution of the differences between species temperatures and mean velocities.  For constant collision frequencies, the convergence is known to be exponential \cite{CrestettoKlingenbergPirner2020}; this behavior can be clearly observed numerically. However, the convergence of these quantities for the velocity-dependent cross-section appears much slower and distinctly different in form.

In Figure \ref{fig:HC_distributions_center}, we plot the kinetic distribution of the hydrogen species for $\nu(\mbv)$ and $\hat{\nu}$, the latter giving the slowest relaxation of the velocity-independent collision frequencies described above.   Since  the  macroscopic  quantities  of  the  heavy  species (carbon) hardly change, we only show the results for the lighter species (hydrogen). 
The relaxation process is weighted by the collision frequencies. 
Because the velocity-dependent cross-section is maximal at $\mbv=\mbu_\mix$ and decays at larger relative velocity, relaxation to equilibrium in the tails of the distribution is slower when using a velocity-dependent cross-section.  

\subsection{Riemann problems}

\subsubsection{Sod problem} \label{test:Sod}
We run a kinetic version of the well-known Sod problem \cite{Sod} in the fluid regime (i.e., with large collision frequencies). 
In the limit of large collision frequencies, the distribution functions can be approximated by Maxwellians:
\begin{equation}
    f_i \simeq { M_i[n_i,\mbu_i,T_i]},
\end{equation}
where $M_i$ is defined in \eqref{eq:Maxwellian}.
With this approximation, the conservation laws \eqref{eq:cons_macros} reduce to the Euler equations. 
We further reduce the problem to the single species case by assuming $m_1=m_2=m$, $\rho_1=\rho_2=\rho$, $\mbu_1=\mbu_2=\mbu$ and $T_1=T_2=T$. 
In one space dimension, with 3 translational degrees of freedom, the single species Euler equations are
\begin{subequations}\label{eq:Euler}
\begin{gather} 
    \partial_t \rho + \nabla_x \cdot  (\rho \mbu) = 0,  \\
    \partial_t  (\rho \mbu) + \nabla_x \cdot (\rho \mbu \otimes \mbu) + \nabla_x p = 0, \\
     \partial_t  \left(\frac{\rho |\mbu|^2}{2}  + \frac{3\rho T}{2m} \right) 
    + \nabla_x \cdot \left(\left(\frac{\rho |\mbu|^2}{2}  + \frac{3\rho T}{2m} +p \right) \mbu  \right)=0,
\end{gather} 
\end{subequations}
where $p = \frac{\rho T}{m}$ denotes the pressure.

This single-species problem can be implemented with the multi-species model by simply treating each species as the same type of particle.  
We set $m_1=m_2=1$ and consider two collision frequencies: one that depends on $\mbv$ 
\begin{equation}
\label{eq:sod_nu_v}
    \nu_{ij}(x,\mbv,t) = 2\cdot 10^4 \frac{n_j}{\delta_{ij}+|\mbv-\mbu_\mix|^3}
\end{equation}
and one that does not: 
\begin{equation}
\label{eq:sod_nu_const}
   \hat{\nu}_{ij}(x,t) = 2\cdot 10^4 \,
   \frac{n_j}{\delta_{ij}+\hat{v}^3} ,
\end{equation}
where the formula for the averaged relative velocity $\hat{v}$ can be found in \eqref{eq:v_hat}. Again we use the regularization parameter $\delta_{ij} = 0.1 \cdot (\Delta v_{ij})^3$ where $\Delta v_{ij} = \frac{1}{4} \sqrt{T_\mix/(2\mu_{ij})}$ and $\mu_{ij} = m_im_j/(m_i+m_j)$.

The initial data is given by $f_i = M_i[n_i,\mbu_i,T_i]$, where
\begin{align}
n_1 = n_2 = 1, 
\qquad \mbu_1 = \mbu_2 = 0,
\qquad T_1 = T_2 = 1,
\end{align}
for $x \leq 0$ and
\begin{align}
n_1 =  n_2 = 0.1, 
\qquad \mbu_1 = \mbu_2 = 0,
\qquad T_1 = T_2 = 0.8.
\end{align}
for $x>0$.

The simulations are run using a velocity grid with $48^3$ points and 400 equally spaced cells in $x$.  We use the second-order IMEX Runge-Kutta scheme from Section \ref{subsec:secondorderIMEX} combined with the second-order finite volume scheme from Section \ref{sec:space}.

\begin{figure}[htb]
\centering
\begin{subfigure}[c]{0.45\textwidth}
\includegraphics[width=\textwidth]{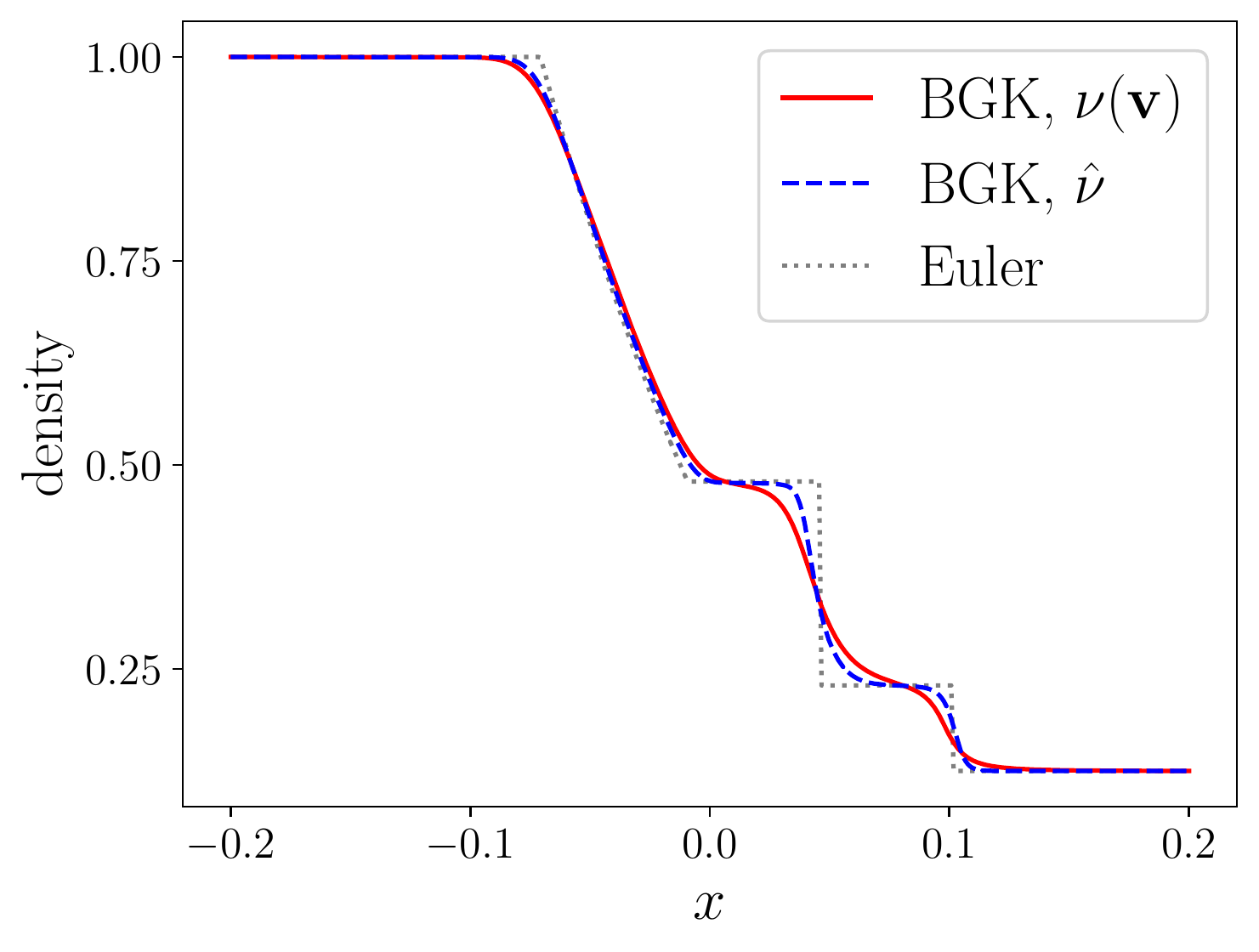}
\caption{density}
\end{subfigure}
~
\begin{subfigure}[c]{0.45\textwidth}
\includegraphics[width=\textwidth]{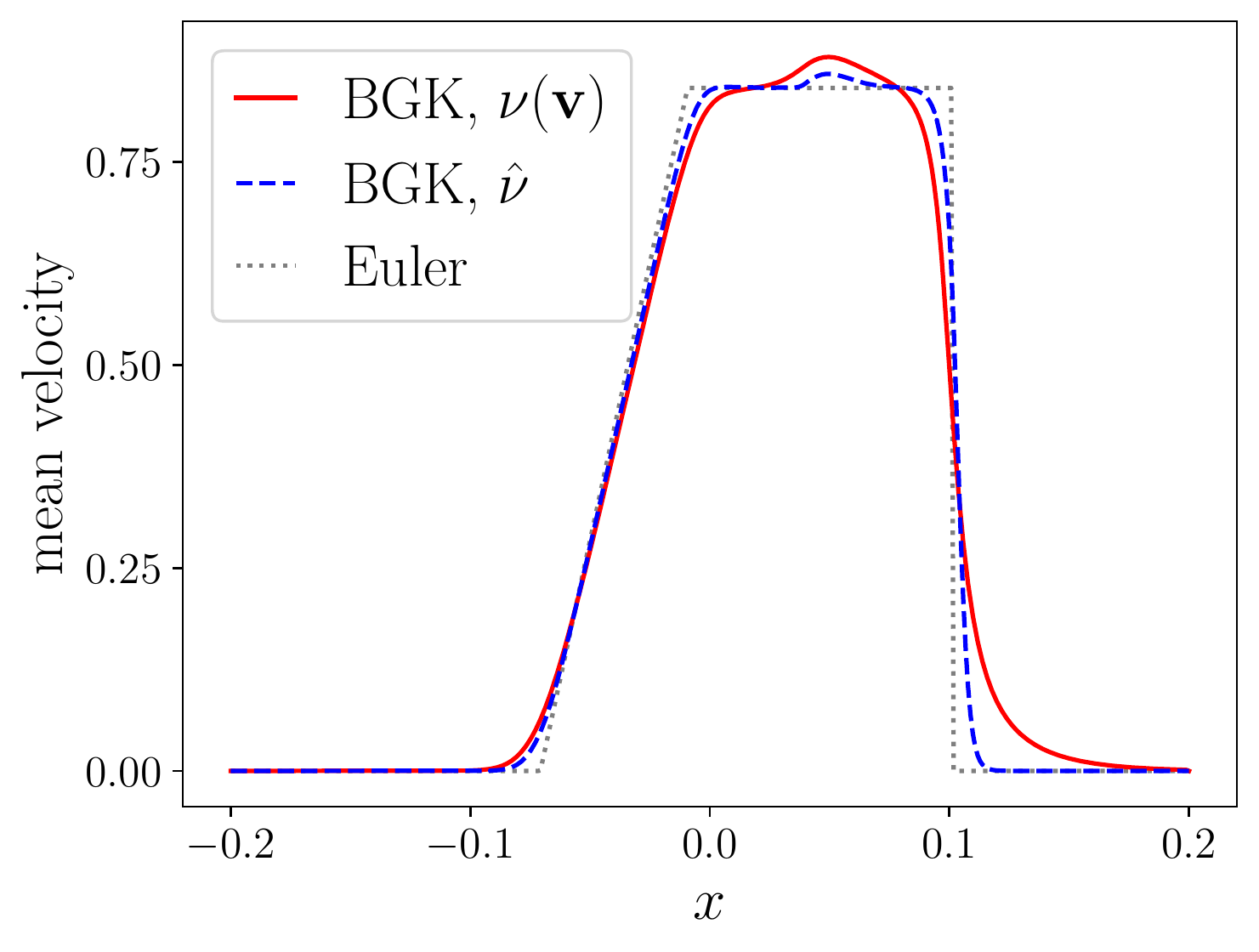}
\caption{mean velocity}
\end{subfigure}
\vskip10pt
\begin{subfigure}[c]{0.45\textwidth}
\includegraphics[width=\textwidth]{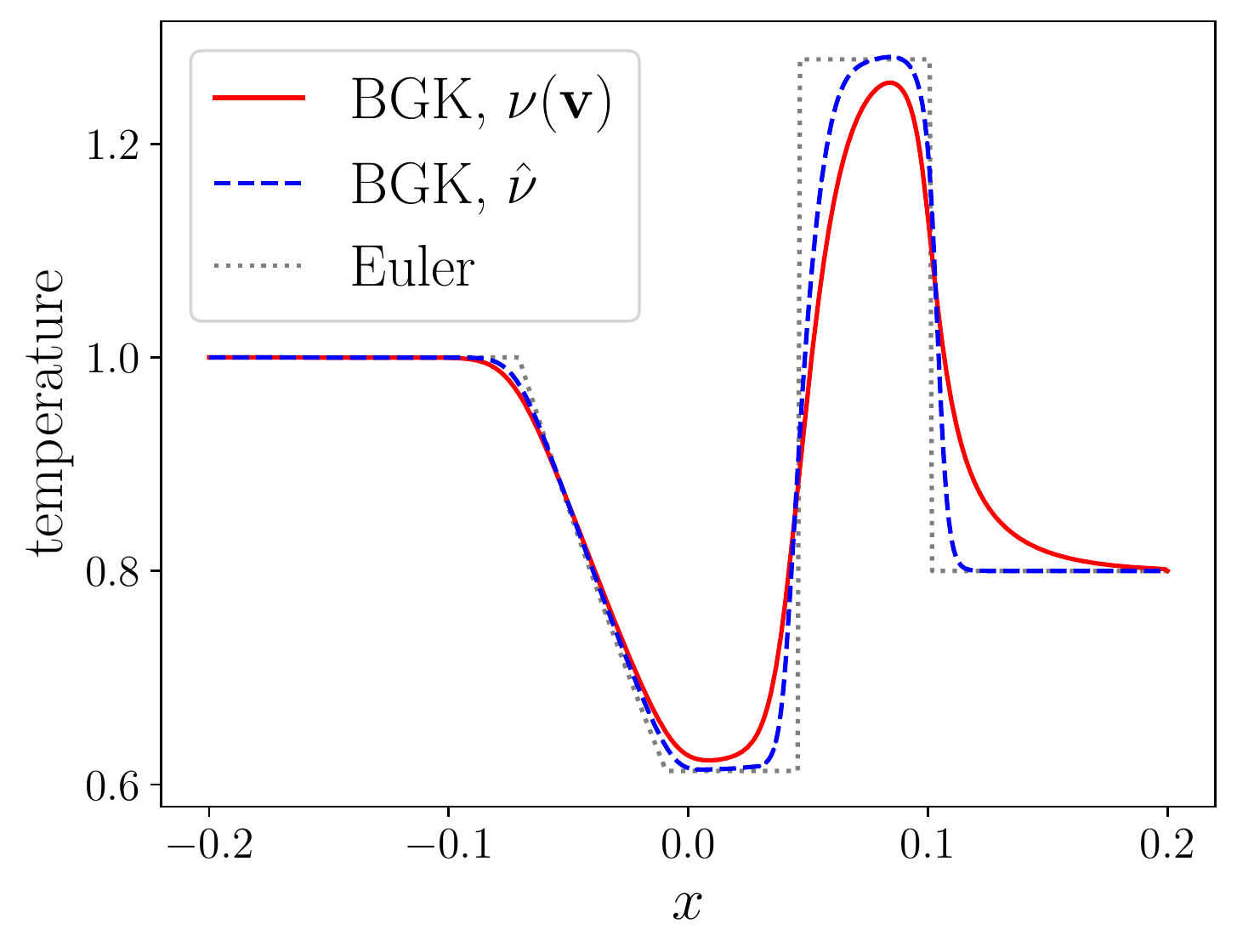}
\caption{temperature}
\end{subfigure}
~
\begin{subfigure}[c]{0.45\textwidth}
\includegraphics[width=\textwidth]{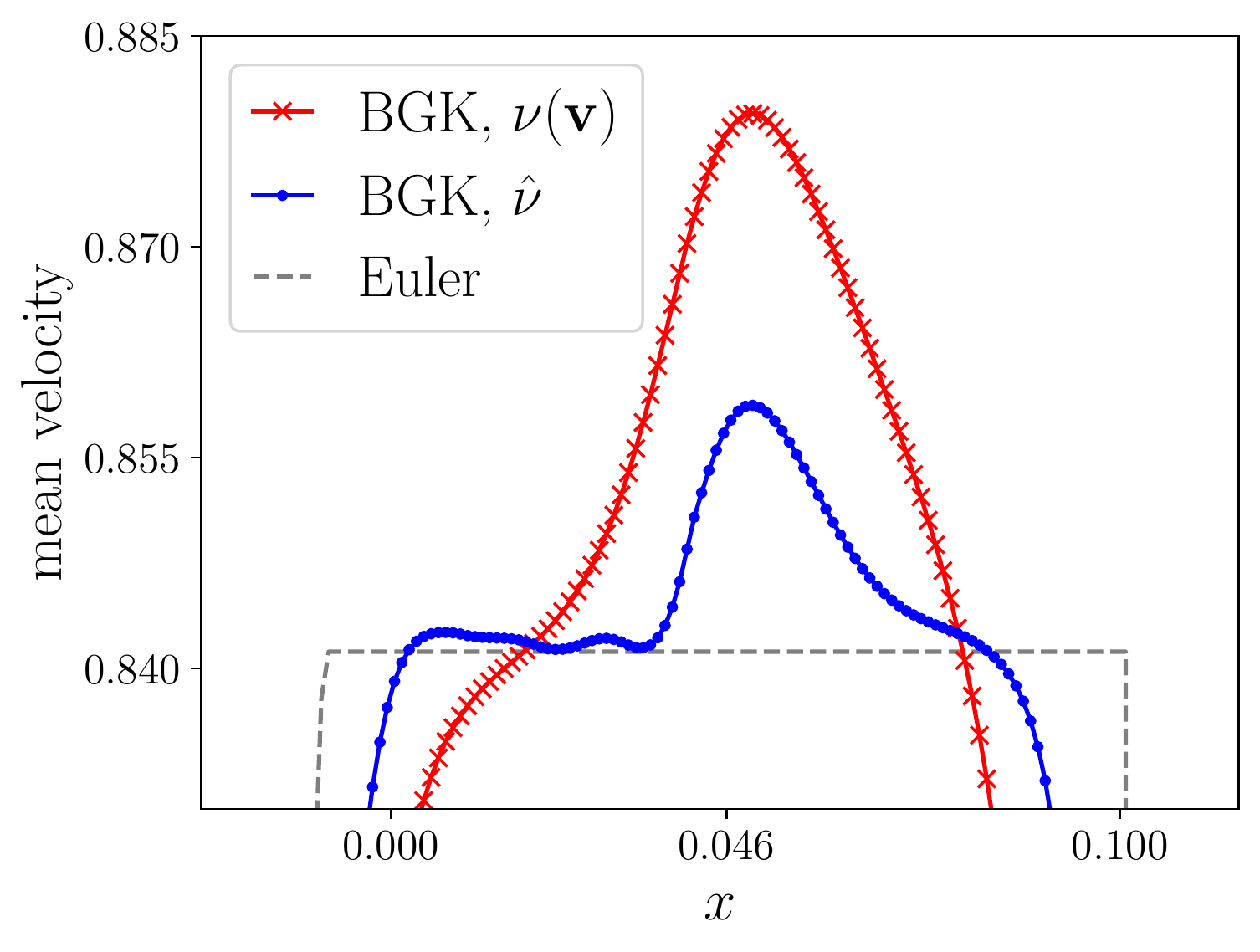}
\caption{mean velocity, closeup}
\end{subfigure}
    \caption{Numerical solution  at $t=0.055$ of the Sod problem in Section \ref{test:Sod}.   We show results for a 2-species kinetic simulation using the velocity-dependent collision frequency $\nu(\mbv)$ in \eqref{eq:sod_nu_v} (red solid line) and the velocity-independent collision frequency $\hat{\nu}$ in \eqref{eq:sod_nu_const} (dashed blue line).  The solutions for both species are identical; we show only the species 1 results.  For reference, the exact solution for the Euler equations \eqref{eq:Euler} is also provided (dotted gray line).  Both kinetic solutions recover the fluid limit fairly well, but the velocity-dependent frequencies contribute to more kinetic behavior around transitions.}
    \label{fig:Sod}
\end{figure}

Numerical simulations of the density, mean velocity, and temperature are given in Figure \ref{fig:Sod}.  We include results using the BGK model with both $\nu(\mbv)$ and $\hat{\nu}$, as well as the analytic solution for the Euler equations in \eqref{eq:Euler}.  Both of the collision frequencies $\nu(\mbv)$ and $\hat{\nu}$ give similar results, but the deviations from the Euler solution near the discontinuities in the fluid model are more pronounced when using $\nu(\mbv)$.

\subsubsection{Mach 1.7 Shock wave problem} \label{test:shock_Mach1.7}

In this example, we compute the flow across a standing Mach 1.7 normal shock wave in a mixture of hydrogen (species 1) and helium (species 2). 
The shock wave structure is difficult to capture in standard hydrodynamic schemes with a single material/species; in mixtures we further expect species separation to occur due to the mass difference between the two species. 
The shock conditions are calculated via the Rankine-Hugoniot jump conditions for a monoatomic gas \cite{Hypersonics}. 
We take a domain size of 6 microns ($6\cdot 10^{-4}\,\text{cm}$) and compute the solution in the frame of the shock.
The masses and charges are (units in cgs) 
\begin{align}
m_1  &= 1.655\cdot 10^{-24} \,\text{g}, 
\quad m_2 = 3.308\cdot 10^{-24} \,\text{g} ,
\quad
Z_1 = 1 ,\quad Z_2 = 2 .
\end{align}
The initial conditions are {$f_i = M_i[n_i,\mbu_i,T_i]$ with}: 
\begin{align}
n_1 = n_2 = 6.666\cdot 10^{19}\,\text{cm}^{-3}  ,
\qquad 
u_1 = u_2 = 1.7634411\cdot 10^{7}\,\frac{\text{cm}}{\text{s}}, 
\qquad 
T_1 = T_2 = 100\,\text{eV},
\end{align} 
for $x \leq 0$ and
\begin{align}
n_1 = n_2 = 1.308\cdot 10^{20}\,\text{cm}^{-3}, 
\qquad
u_1 = u_2 = 8.985007\cdot 10^6 \,\frac{\text{cm}}{\text{s}},
\qquad 
T_1 = T_2 = 171.32\,\text{eV}
\end{align}
for  $x>0$.

The simulations are run using a velocity grid with $48^3$ nodes and spatial mesh with 200 cells.  We use the second-order IMEX Runge-Kutta scheme from Section \ref{subsec:secondorderIMEX} and the second-order spatial discretization in Section \ref{sec:space}, with the limiter given in \eqref{eq:second-order-limiter}. The time step $\Delta t = 22$ fs is set according to the  CFL condition in \eqref{eq:CFL}.

\begin{figure}[htb]
\begin{subfigure}[c]{0.33\textwidth}
\includegraphics[width=\textwidth]{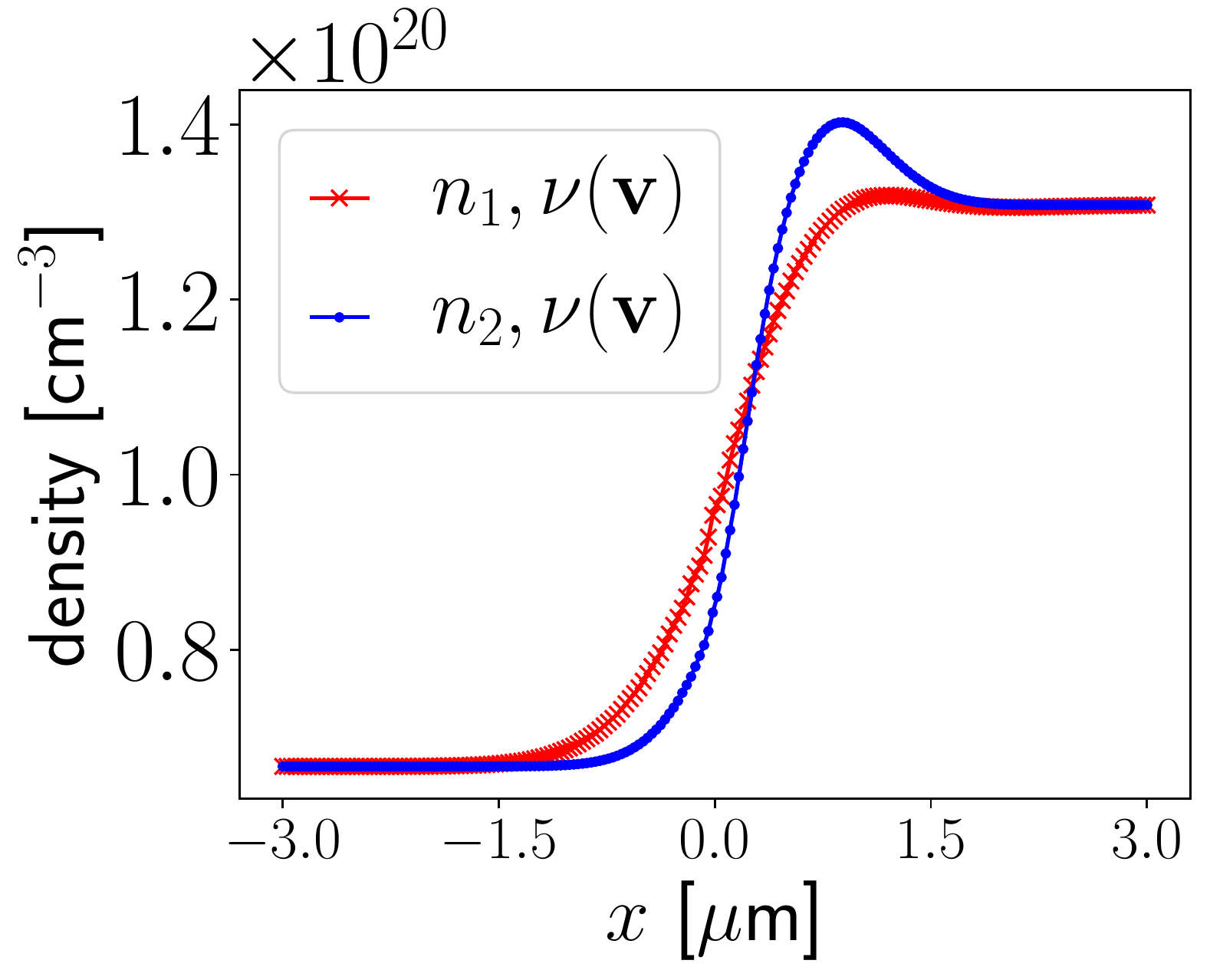}
\end{subfigure}
\begin{subfigure}[c]{0.33\textwidth}
\includegraphics[width=\textwidth]{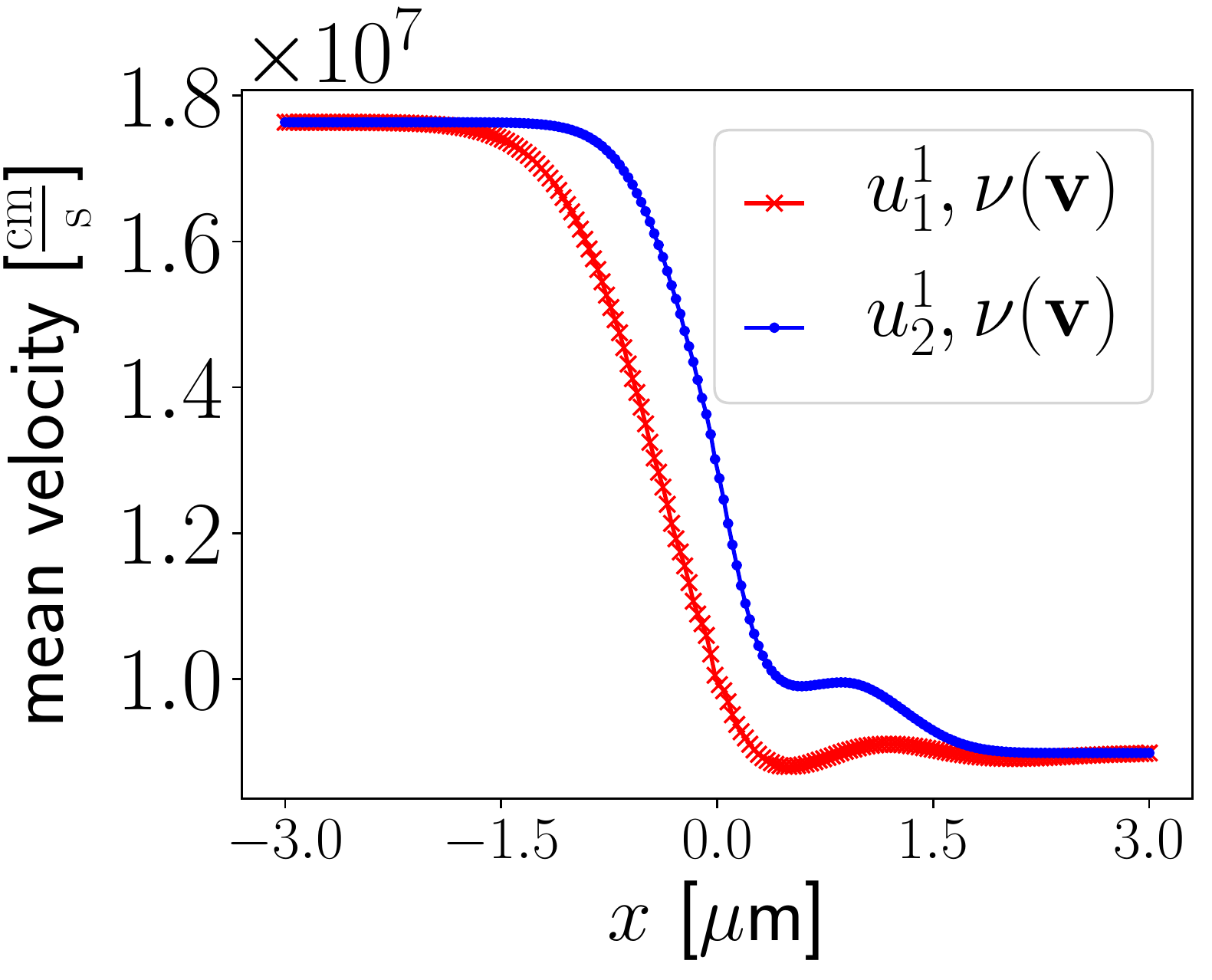}
\end{subfigure}
\begin{subfigure}[c]{0.33\textwidth}
\includegraphics[width=\textwidth]{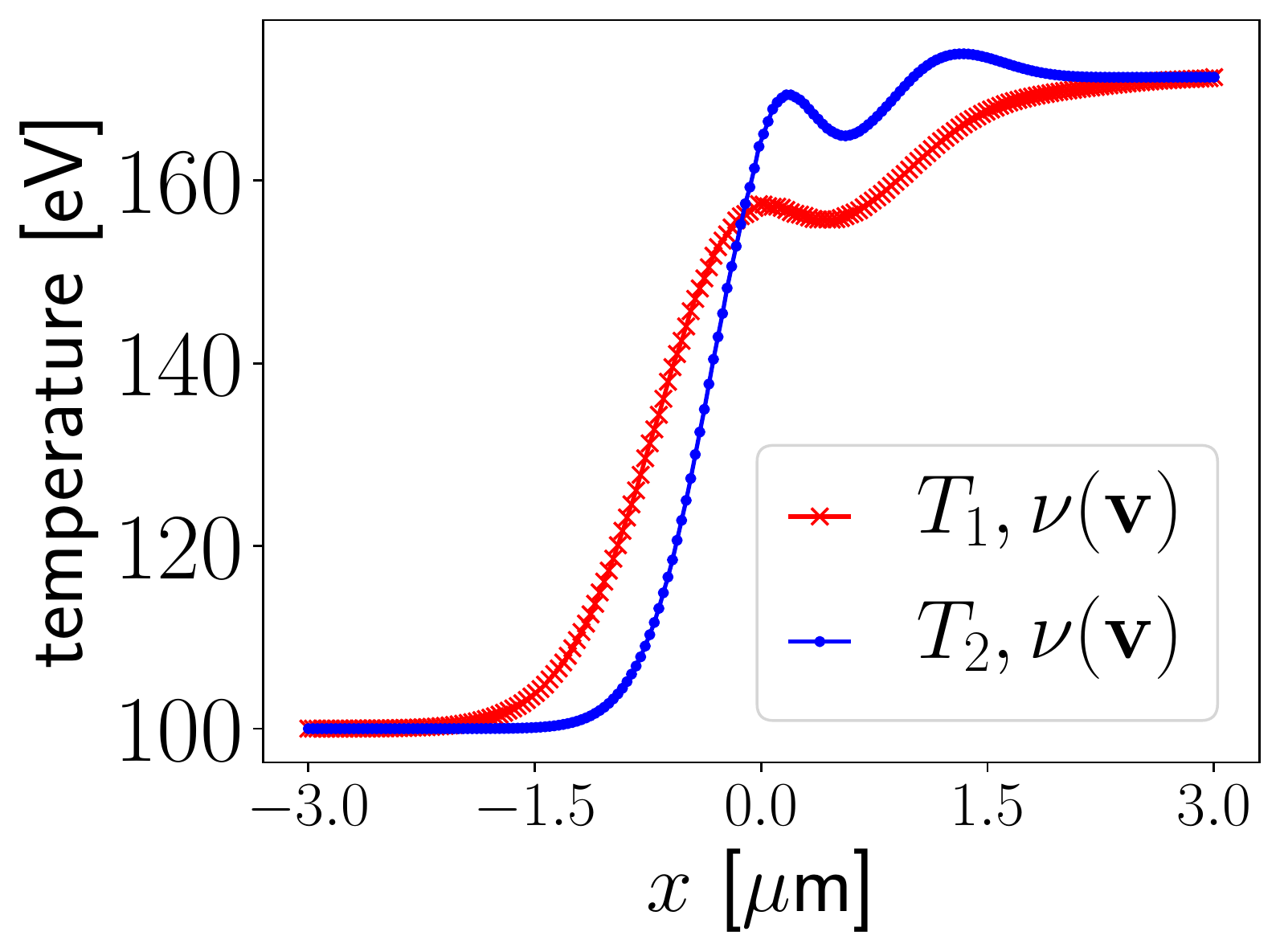}
\end{subfigure}

\begin{subfigure}[c]{0.33\textwidth}
\includegraphics[width=\textwidth]{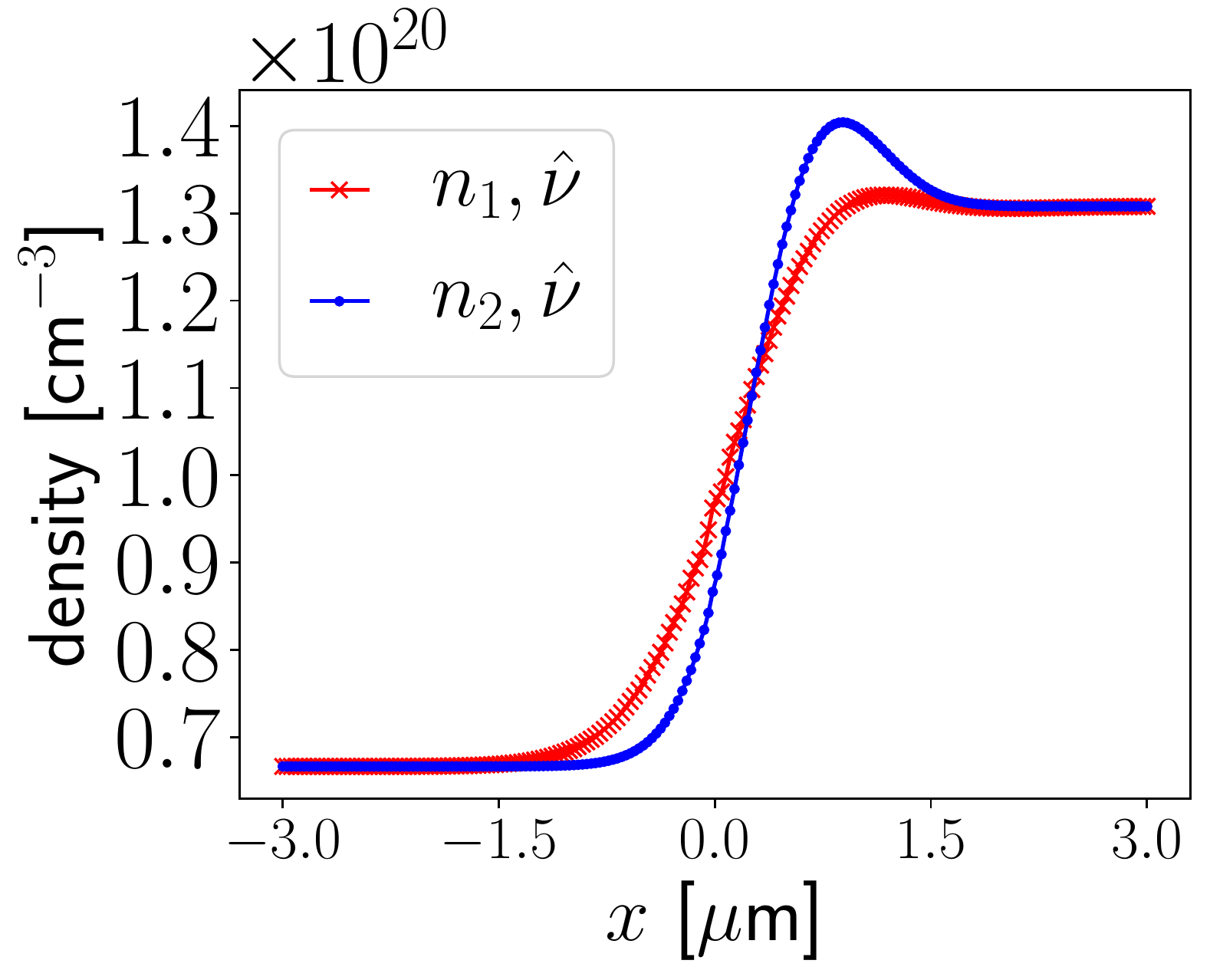}
\end{subfigure}
\begin{subfigure}[c]{0.33\textwidth}
\includegraphics[width=\textwidth]{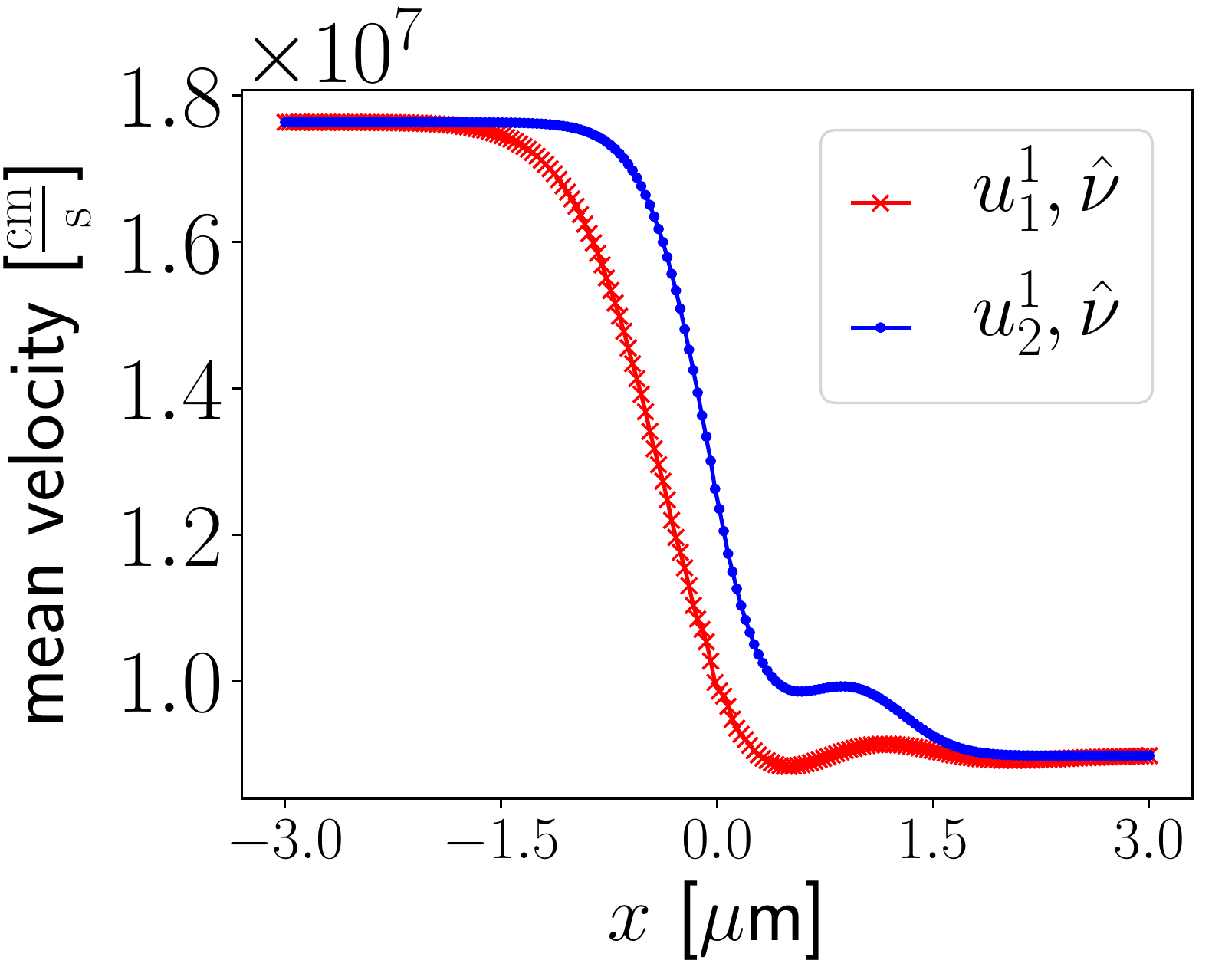}
\end{subfigure}
\begin{subfigure}[c]{0.33\textwidth}
\includegraphics[width=\textwidth]{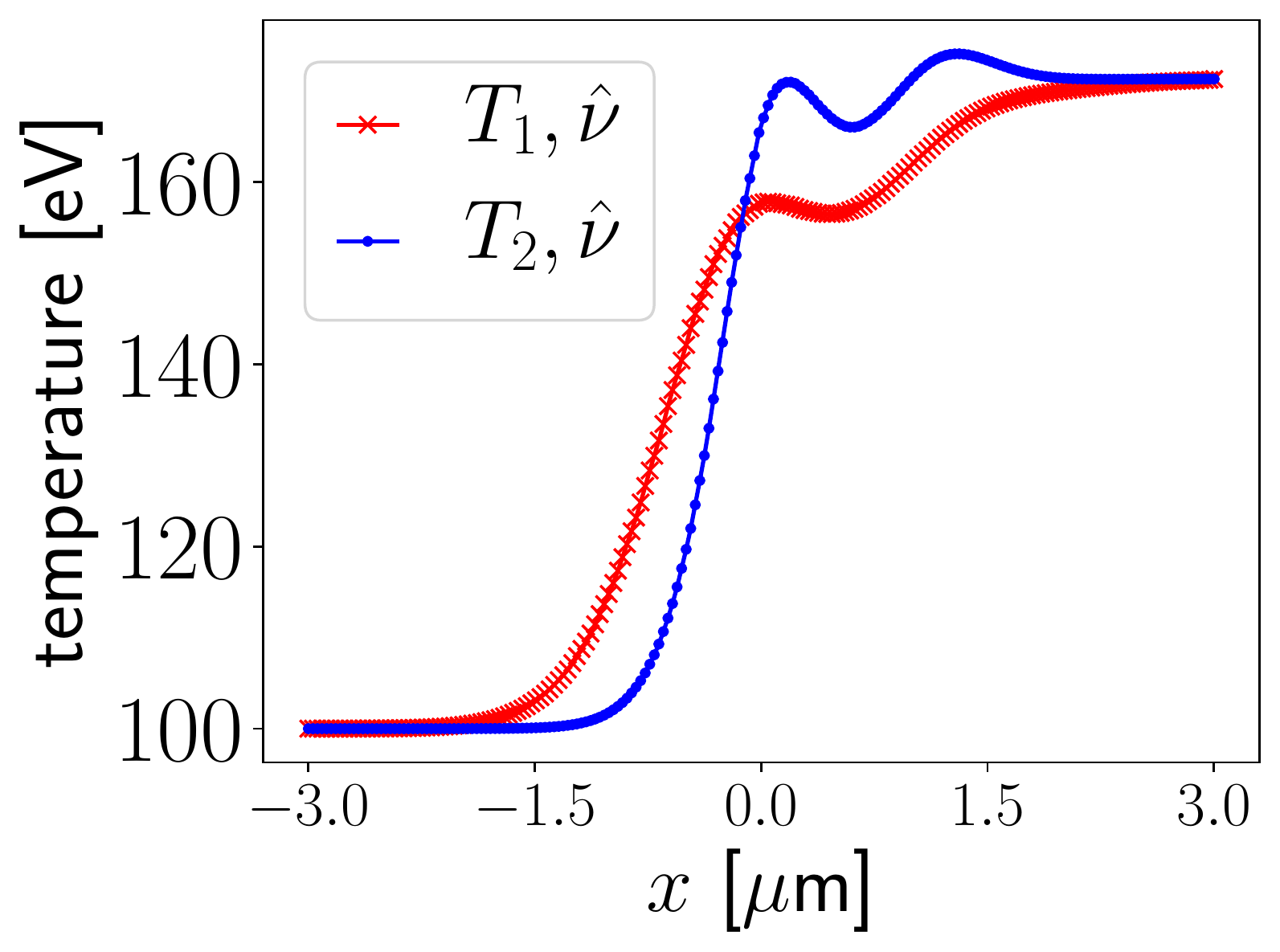}
\end{subfigure}


\begin{subfigure}[c]{0.33\textwidth}
\includegraphics[width=\textwidth]{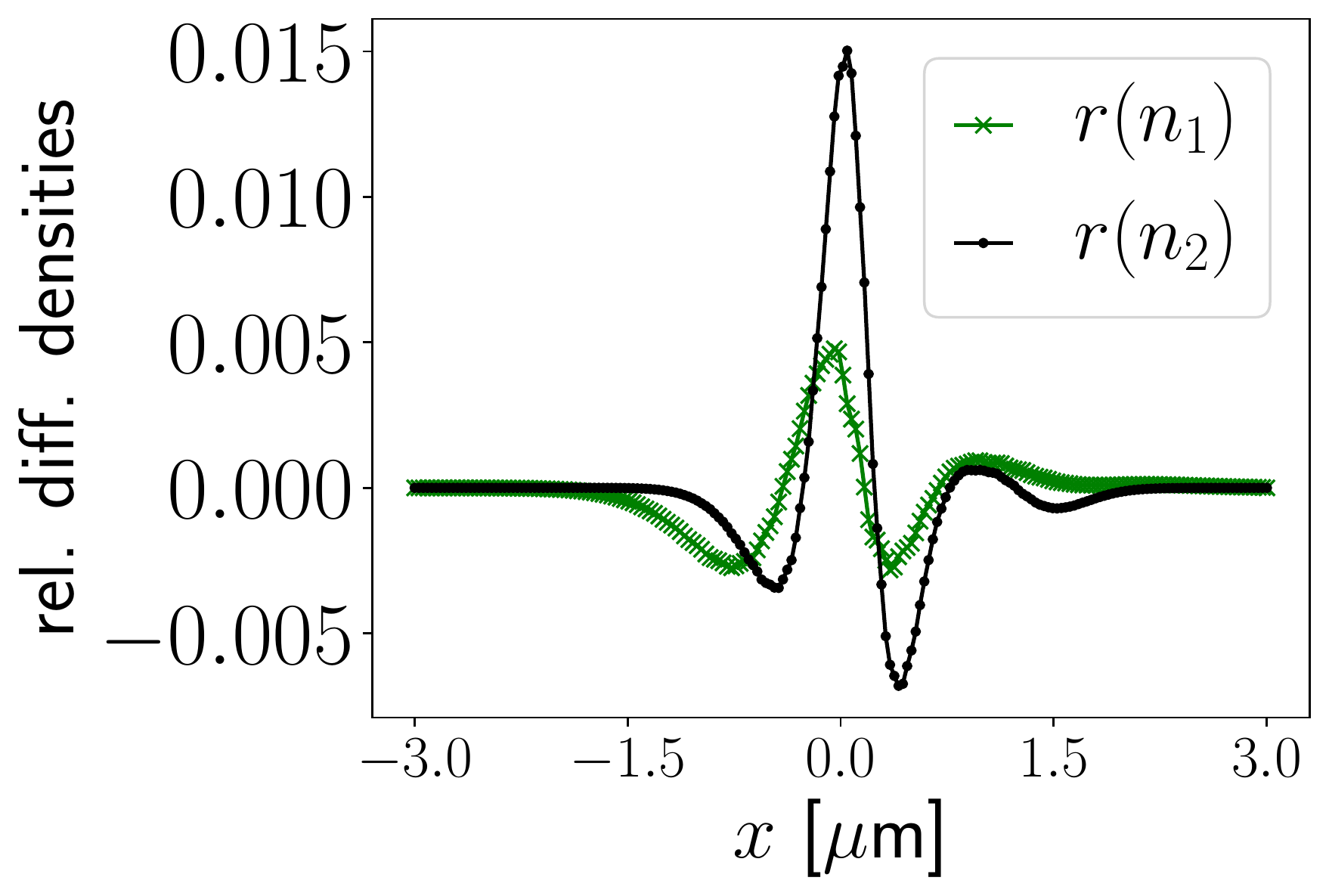}
\end{subfigure}
\begin{subfigure}[c]{0.33\textwidth}
\includegraphics[width=\textwidth]{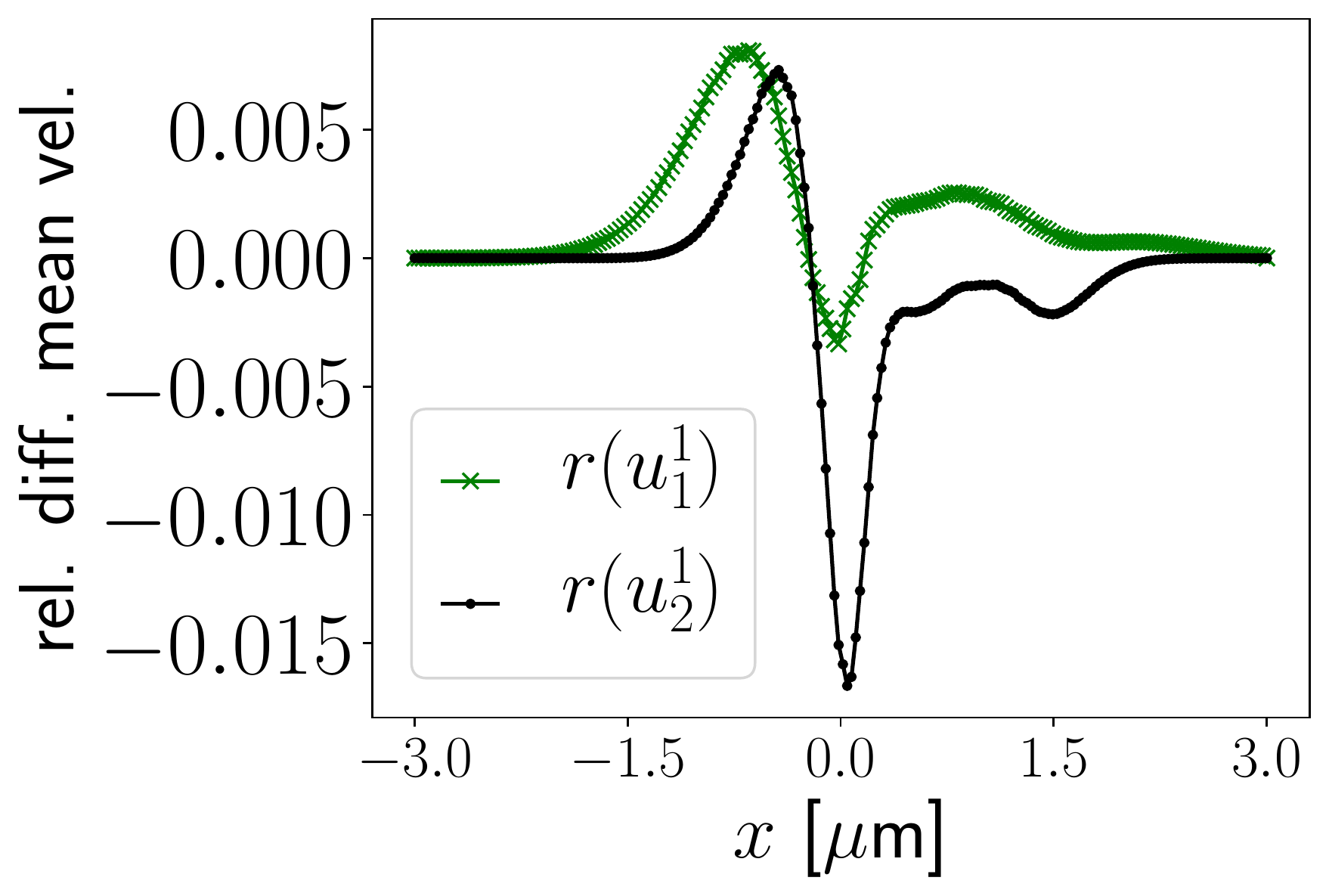}
\end{subfigure}
\begin{subfigure}[c]{0.33\textwidth}
\includegraphics[width=\textwidth]{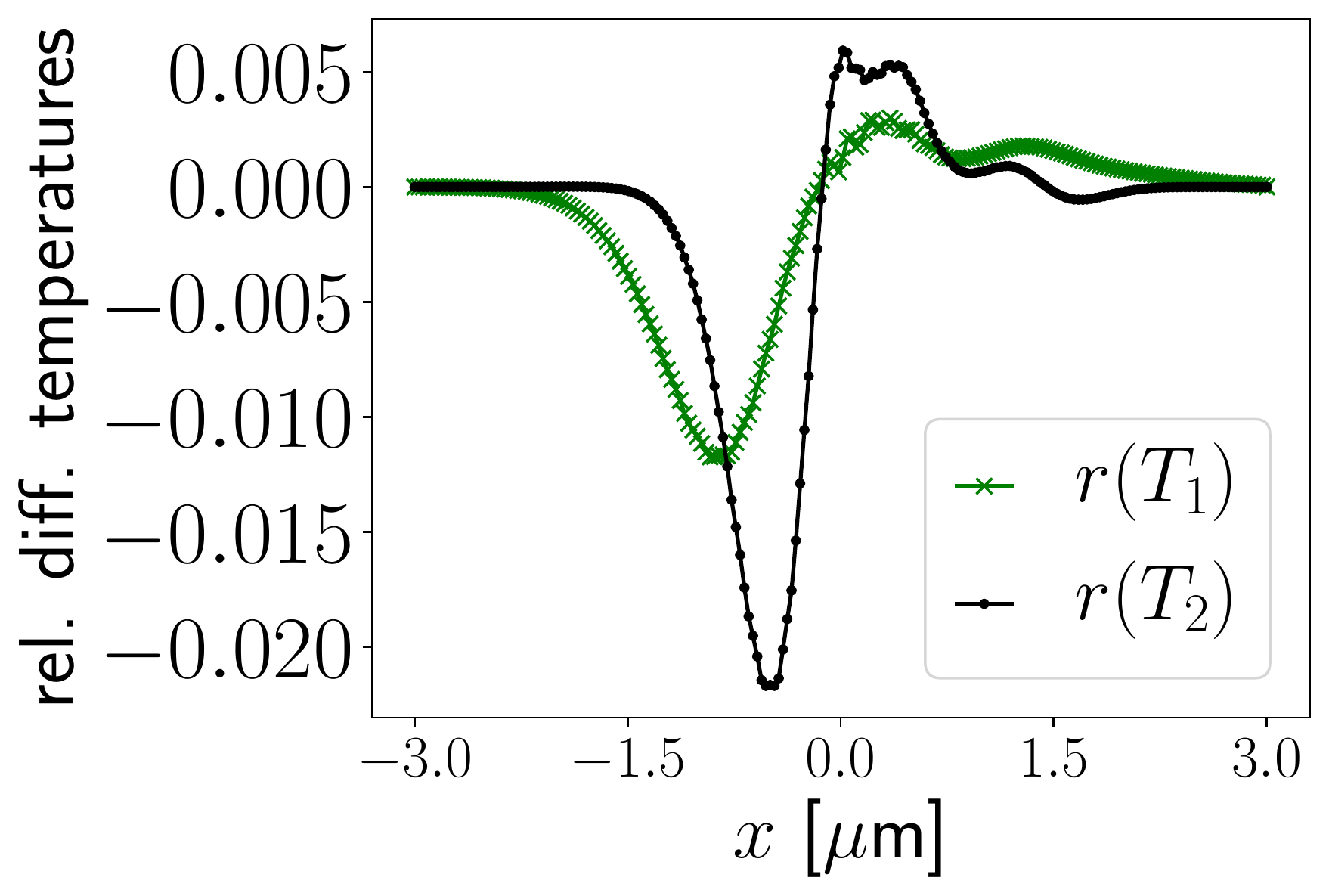}
\end{subfigure}

    \caption{ Fluid quantities at time $t=5.390$ ps for the  Mach 1.7 shock wave problem from Section \ref{test:shock_Mach1.7}. The initial Riemann is evolving to a multi-species normal shock wave which exhibits separation between hydrogen (species 1) and helium (species 2) particles. Top row: velocity-dependent collision frequencies $\nu_{ij}$, given in \eqref{eq:colfreq_dep}; middle row: the constant collision frequencies $\hat{\nu}_{ij}$, given in \eqref{eq:colfreq_indep_2}; bottom row: relative difference (see \eqref{eq:rel_diff_new}).  These differences are typically within  $2\%$.} 
    \label{fig:shockwave_HHe_Mach1.7}
\end{figure}
In Figure \ref{fig:shockwave_HHe_Mach1.7} we compare numerical results during the transient regime at time $t=5.390$ ps using the velocity-dependent collision frequency $\nu(\mbv)$, given in \eqref{eq:colfreq_dep}, with those using the constant collision frequencies $\hat{\nu}$, given in \eqref{eq:colfreq_indep_2}. 
In addition to these results, we plot the relative difference 
\begin{equation}
\label{eq:rel_diff_new}
   r(q)=\frac{q(\{\hat{\nu}_{ij}\})-q(\{\nu_{ij}(\mbv)\})}{|q(\{\hat{\nu}_{ij}\})|+|q(\{\nu_{ij}(\mbv)\})|} 
\end{equation} 
for the densities ($q=n_i$), mean velocities ($q=u_i^1$), and temperatures ($q=T_i$). 
As expected, both the velocity-dependent and constant collision frequency models show a species separation.
For all hydrodynamic quantities, the differences are within a few percent.
While we expect a difference in output profiles between the two models due to the tail particles relaxing more slowly than the bulk, it is likely that the collision frequencies outside of the `kinetic' region of the shock interface are high enough to suppress large deviations from equilibrium.

\subsubsection{ Mach 4 Shock wave problem} \label{test:shock_Mach4}

In this test case, we repeat the normal shock wave in a hydrogen-helium mixture from the previous test case, but increase the shock strength in the mixture to Mach 4 with the expectation that the distributions will be further out of equilibrium than the previous case.   
The species masses and charges are the same as in the Mach 1.7 case, but we widen the domain size to 12 microns and modify the initial conditions to construct a Mach 4 shock, again using the Rankine-Hugoniot relations. 
Specifically we set $f_i = M_i[n_i,\mbu_i,T_i]$ where 
\begin{equation}
n_1 = n_2 = 3.3488\cdot 10^{19}\,\text{cm}^{-3},
\qquad
u_1 = u_2=5.06 \cdot 10^7 \,\frac{\text{cm}}{\text{s}},
\qquad 
T_1 = T_2 = 100\,\text{eV} ,
\end{equation}
for $x \leq 0$,  and 
\begin{equation}
n_1 = n_2 = 1.128\cdot 10^{20}\,\text{cm}^{-3},
\qquad
u_1 = u_2=1.50\cdot 10^7 \,\frac{\text{cm}}{\text{s}},
\qquad 
T_1 = T_2 = 586.3\,\text{eV} 
\end{equation}
for $x>0$.

As in the previous case, simulations are run using a velocity grid with $48^3$ nodes and a spatial mesh with 200 cells.  
We use the second-order IMEX Runge-Kutta scheme from Section \ref{subsec:secondorderIMEX}  and the second-order spatial discretization in Section \ref{sec:space}, with the limiter given in \eqref{eq:second-order-limiter}.
The time step $\Delta t = 25$ fs is set according to the  CFL condition in \eqref{eq:CFL}. 

\begin{figure}[htb]
\begin{subfigure}[c]{0.33\textwidth}
\includegraphics[width=\textwidth]{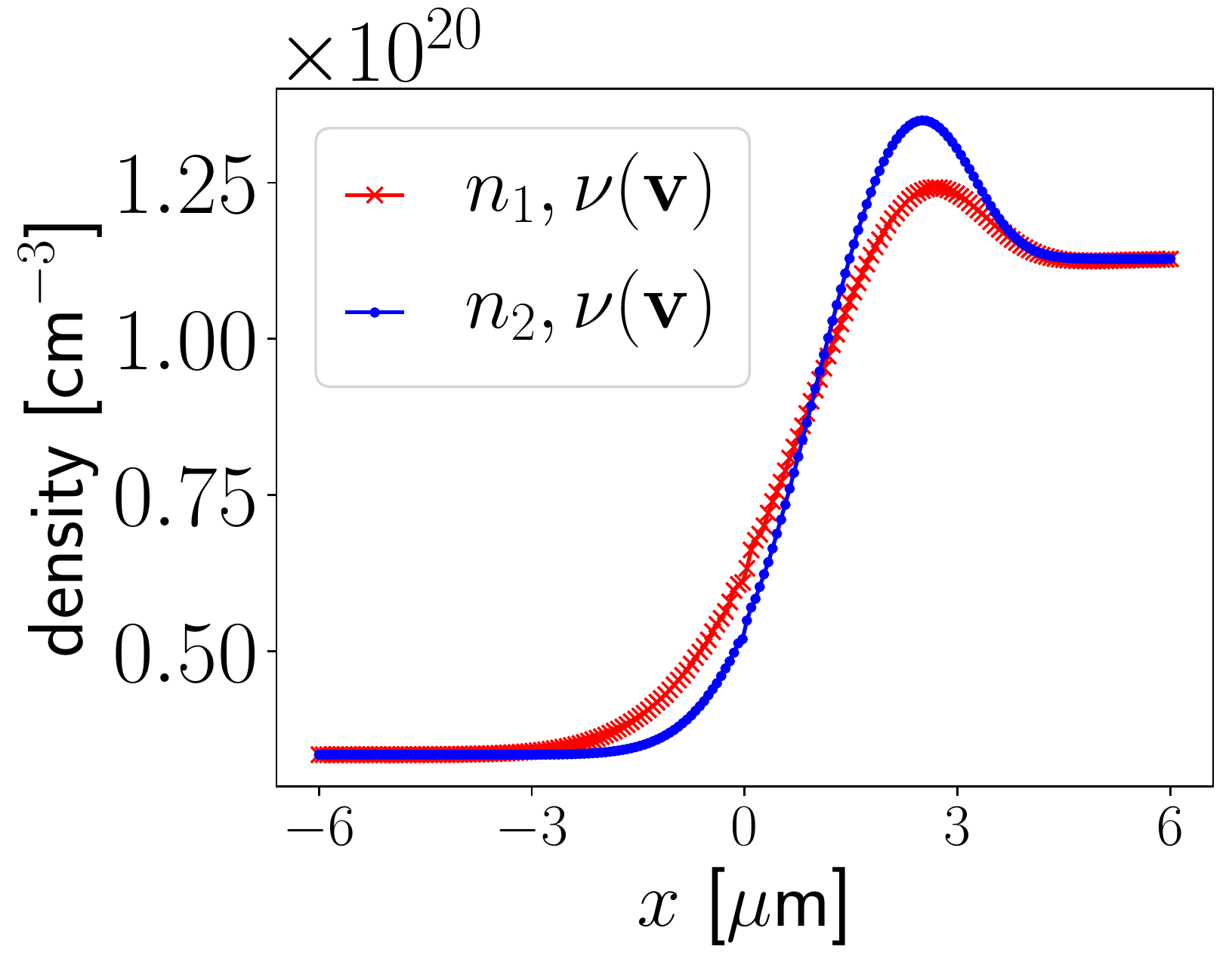}
\end{subfigure}
\begin{subfigure}[c]{0.33\textwidth}
\includegraphics[width=\textwidth]{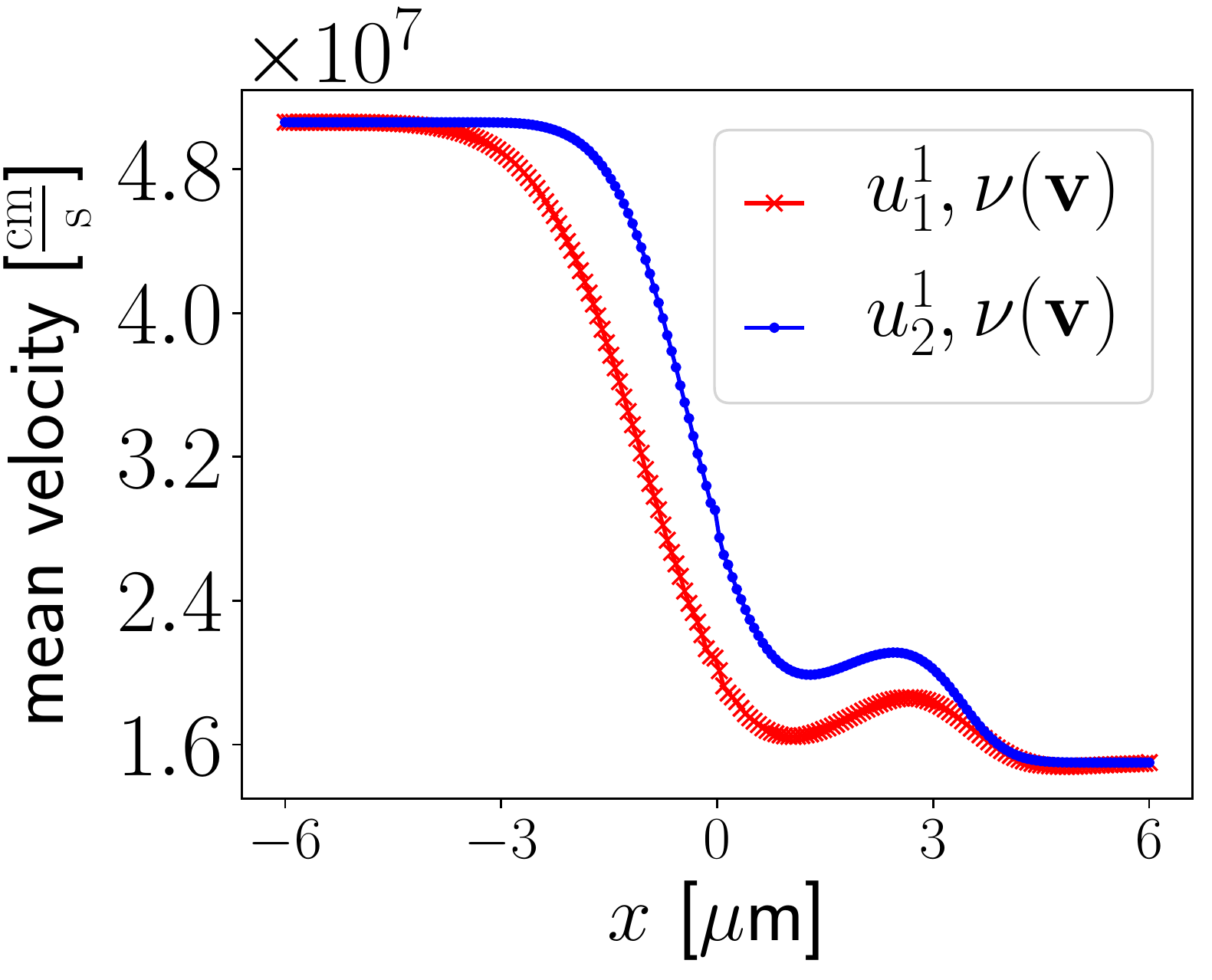}
\end{subfigure}
\begin{subfigure}[c]{0.33\textwidth}
\includegraphics[width=\textwidth]{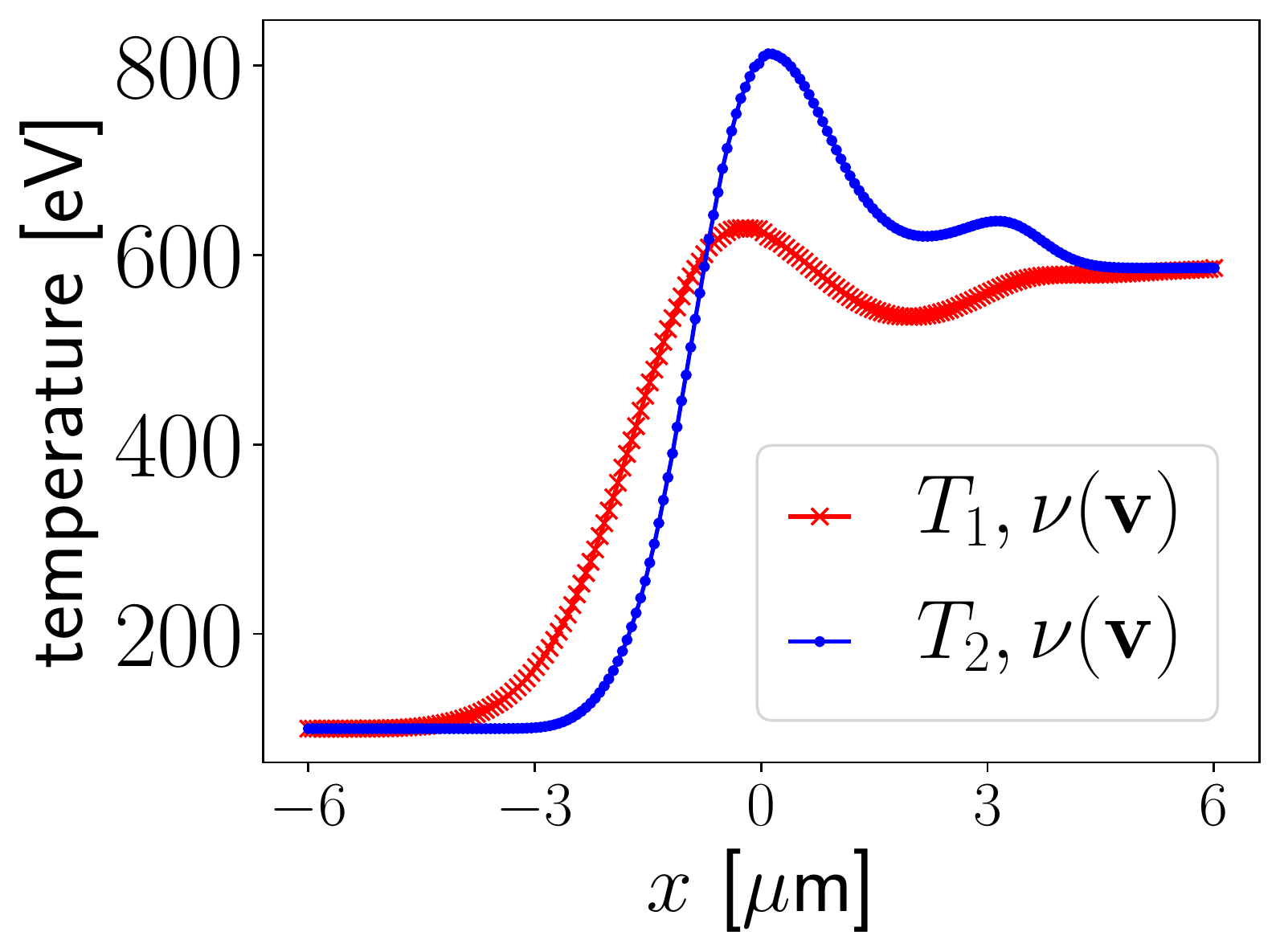}
\end{subfigure}

\begin{subfigure}[c]{0.33\textwidth}
\includegraphics[width=\textwidth]{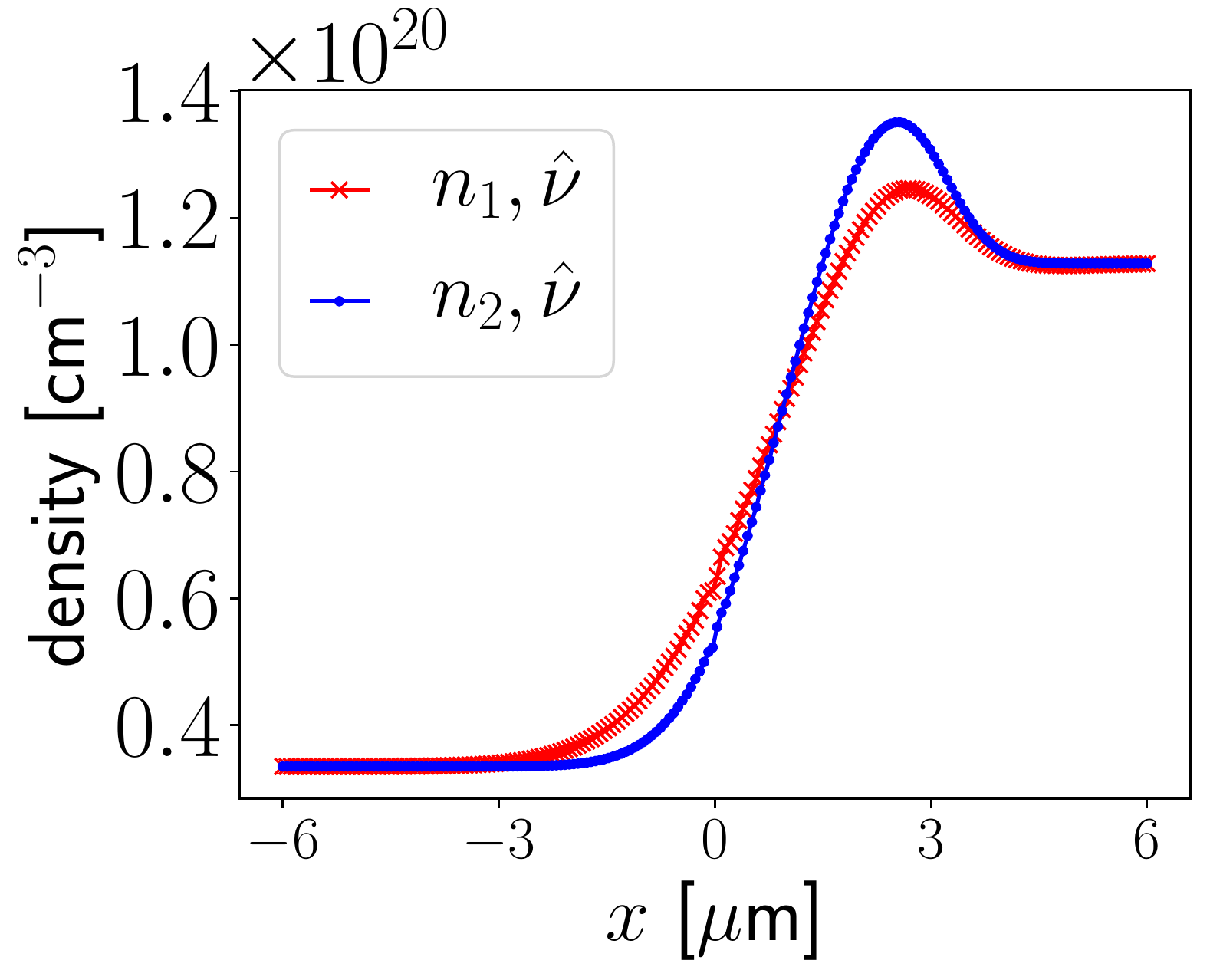}
\end{subfigure}
\begin{subfigure}[c]{0.33\textwidth}
\includegraphics[width=\textwidth]{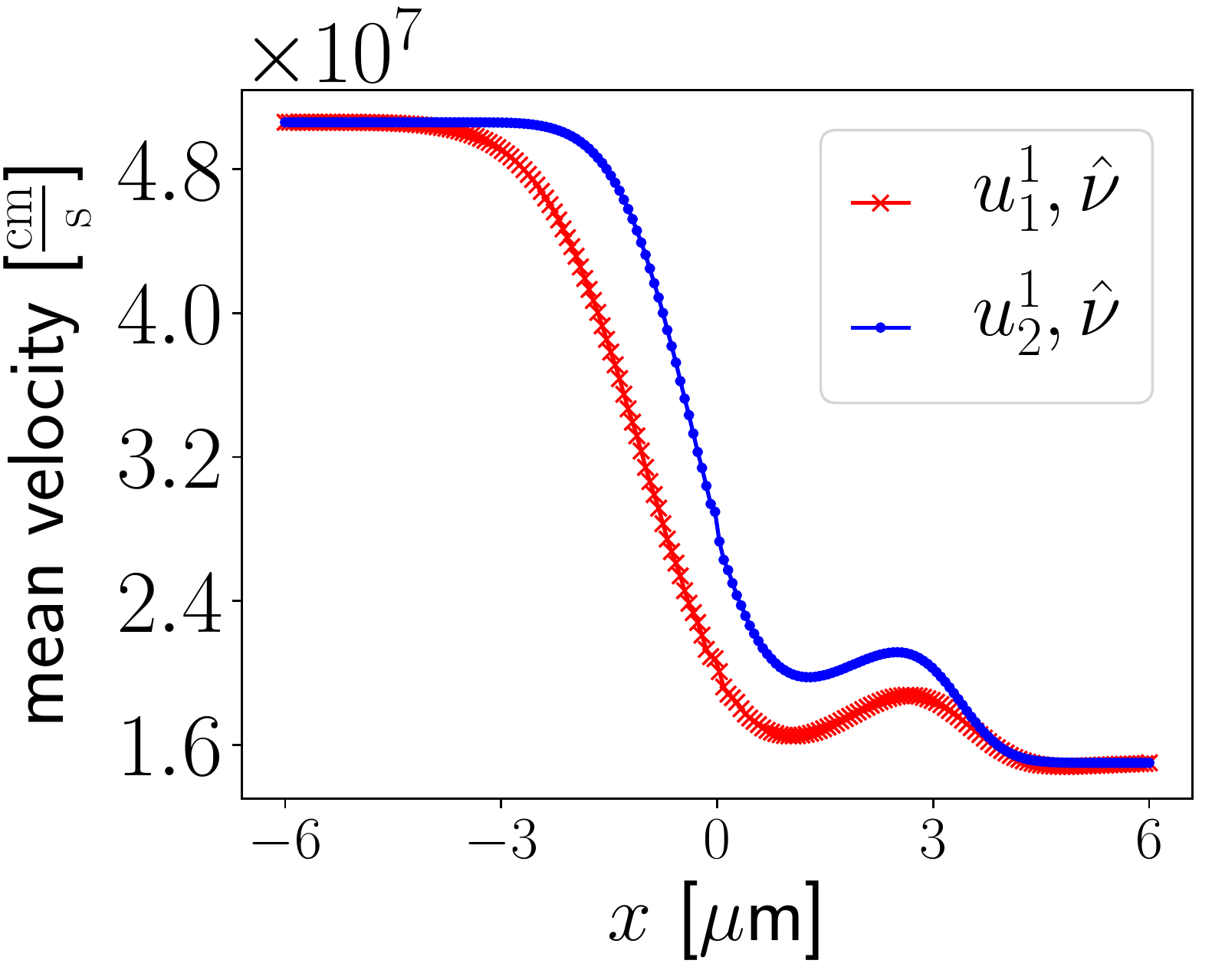}
\end{subfigure}
\begin{subfigure}[c]{0.33\textwidth}
\includegraphics[width=\textwidth]{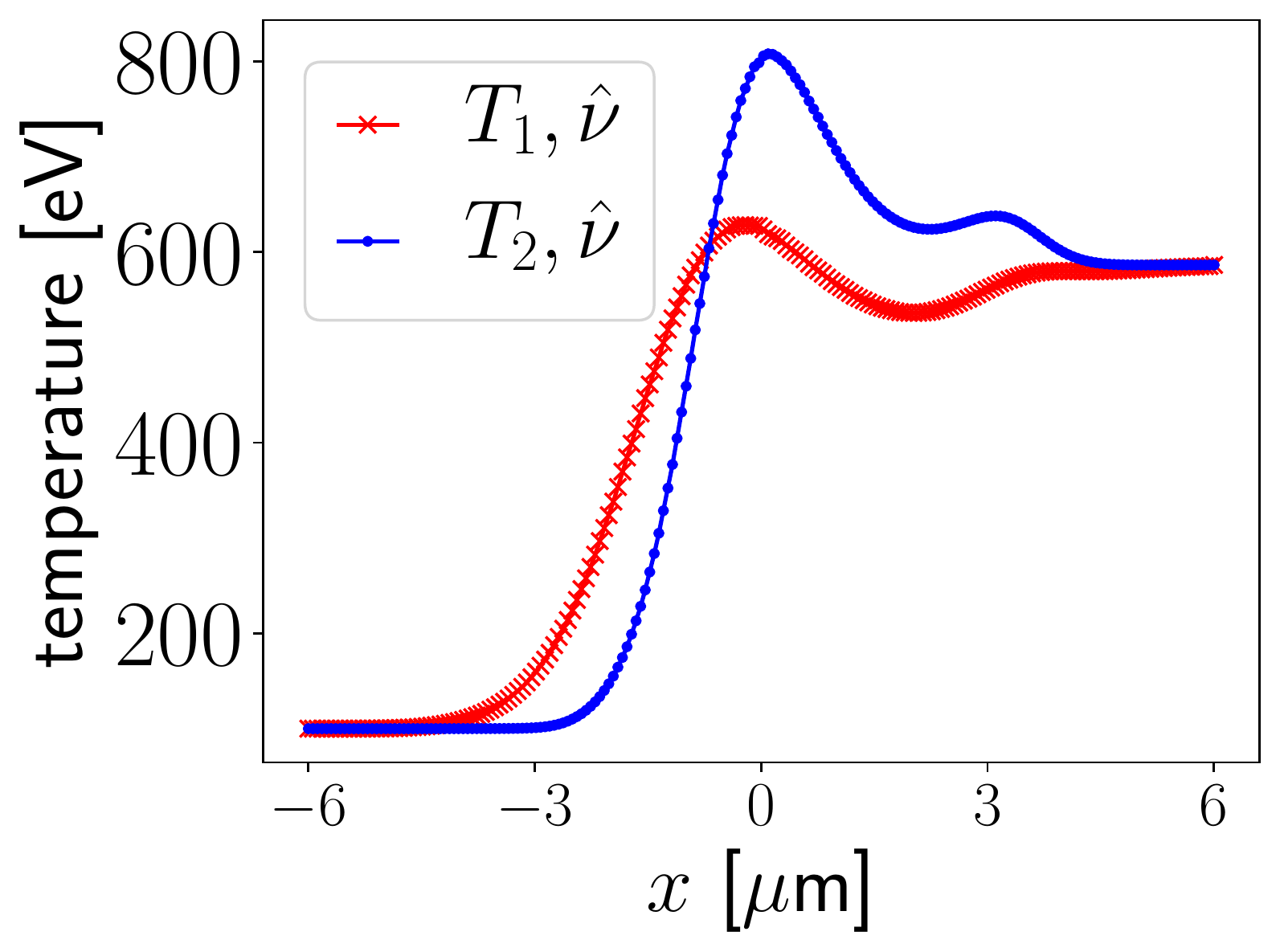}
\end{subfigure}


\begin{subfigure}[c]{0.33\textwidth}
\includegraphics[width=\textwidth]{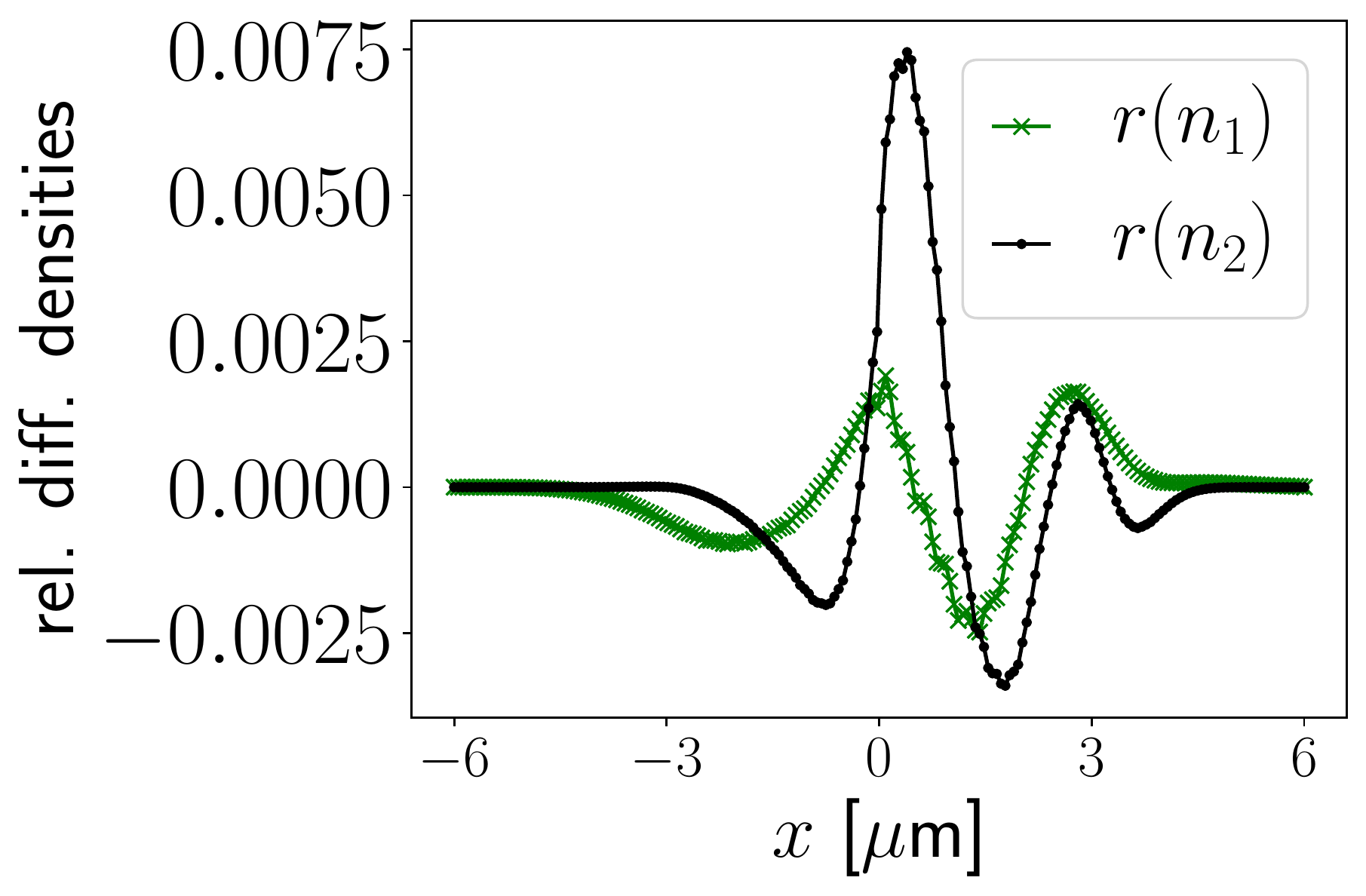}
\end{subfigure}
\begin{subfigure}[c]{0.33\textwidth}
\includegraphics[width=\textwidth]{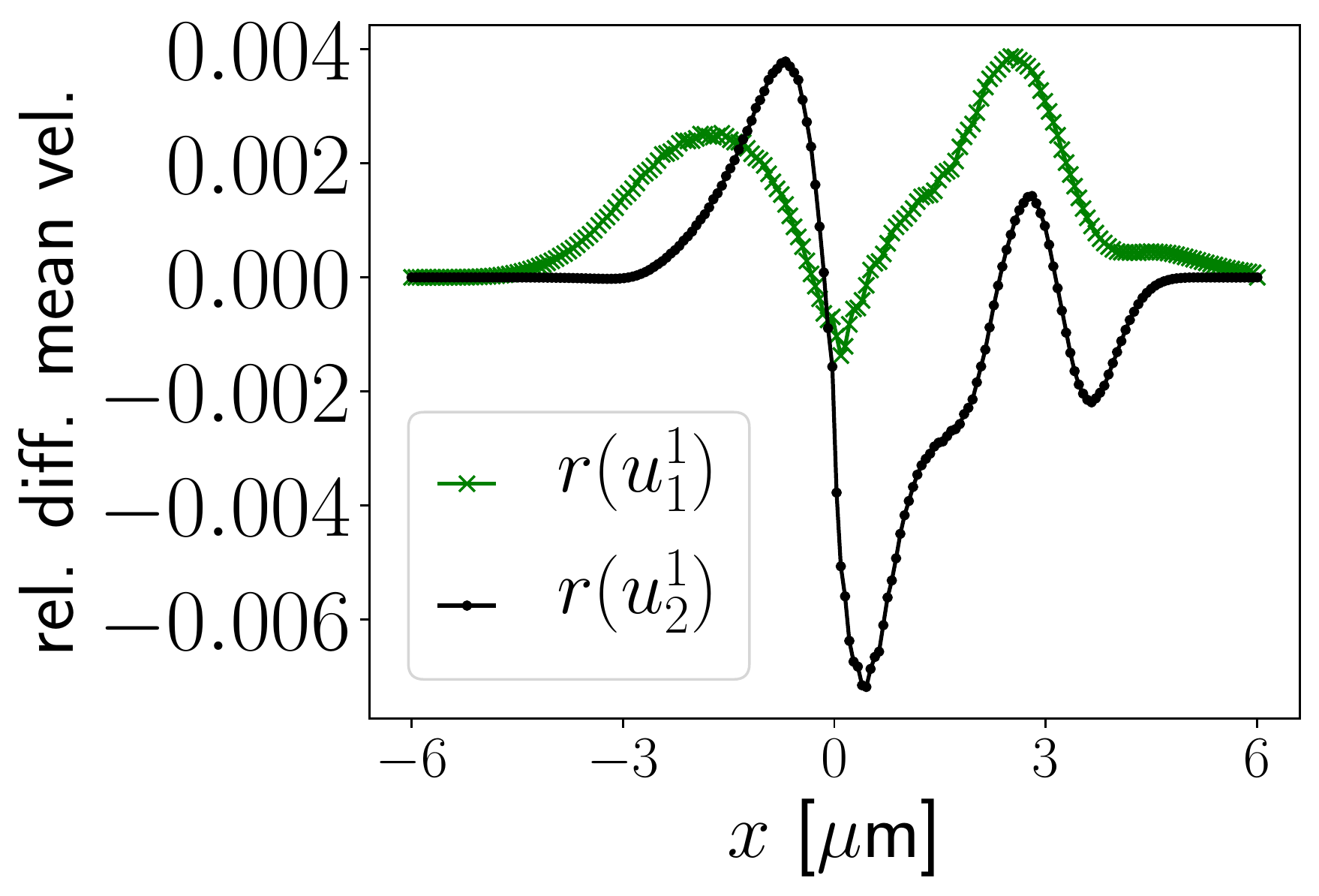}
\end{subfigure}
\begin{subfigure}[c]{0.33\textwidth}
\includegraphics[width=\textwidth]{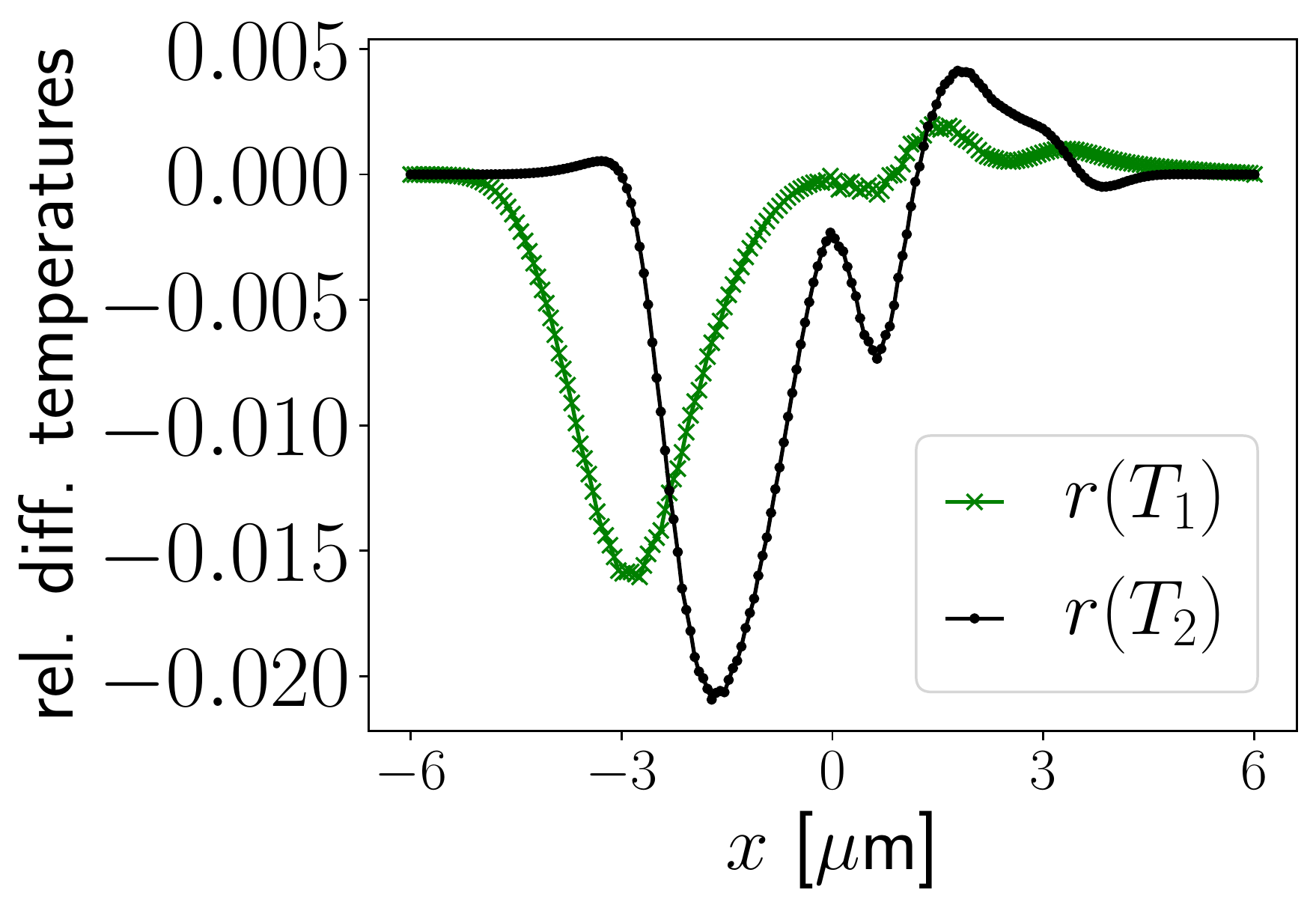}
\end{subfigure}
    \caption{The fluid quantities for the Mach 4 shock wave problem from Section \ref{test:shock_Mach4} (Mach 4) are presented at time $t=6.345$ ps. Top row: velocity-dependent collision frequencies $\nu_{ij}$, given in \eqref{eq:colfreq_dep}; middle row: the constant collision frequencies $\hat{\nu}_{ij}$, given in \eqref{eq:colfreq_indep_2}; bottom row: relative difference (see \eqref{eq:rel_diff_new}). The results for the velocity-dependent collision frequency and the constant collision frequency $\hat{\nu}_{ij}$ look very similar at first glance. However, the relative differences clarify the disparities. They differ up to 2 \%, which is similar to the relative difference seen in the weaker shock wave problem in Section \ref{test:shock_Mach1.7} (Mach 1.7).} 
    \label{fig:shockwave_HHe_Mach4}
\end{figure}

In Figure \ref{fig:shockwave_HHe_Mach4} we compare numerical results during the transient regime at time $t=6.345$ ps using the velocity-dependent collision frequency $\nu(\mbv)$, given in \eqref{eq:colfreq_dep}, with those using the constant collision frequencies $\hat{\nu}$, given in \eqref{eq:colfreq_indep_2}. 
We again observe the evolution towards a standing shock wave for both the velocity-dependent collision frequency $\nu(\mbv)$  and the constant collision frequency $\hat{\nu}$.  
As in the Mach 1.7 case above, while we expect a difference in output profiles between the two models due to the tail particles relaxing more slowly than the bulk, it is likely that the collision frequencies outside of the `kinetic' region of the shock interface are high enough to suppress large deviations from equilibrium for this test problem.

\subsubsection{Interpenetration problem:  high density} \label{test:interpenetration1}

Standard hydrodynamic models have great difficulty in capturing interpenetrating flows of rarefied gases. 
For example in inertial confinement fusion (ICF) simulations, colliding streams of blown-off hohlraum wall particles and capsule ablator particles result in an unphysical density spike due to the lack of interpenetration in hydrodynamic models, which interferes with laser energy propogation in the integrated simulation. 
This discrepancy has been proposed as a cause of symmetry discrepancies in capsule drive between experiments and simulations in ICF \cite{HopkinsPop}.

For this numerical example, we simulate the dynamics of two counter-streaming beams of different species. 
We take a domain size of 50 microns ($50\cdot 10^{-4}\,\text{cm}$) and compute the solution when hydrogen (species 1) interpenetrate with helium (species 2) particles.
We include a trace amount of each species in the whole domain as a background for ease of computation.
The masses and charges are (units in cgs) 
\begin{align}
m_1  &= 1.655\cdot 10^{-24} \,\text{g}, 
\quad m_2 = 3.308\cdot 10^{-24} \,\text{g} ,
\quad
Z_1 = 1 ,\quad Z_2 = 2.
\end{align}
The initial conditions are {$f_i = M_i[n_i,\mbu_i,T_i]$ with}: 
\begin{align}
n_1 = 10^{20}\,\text{cm}^{-3},
\qquad 
n_2 = 10^{17}\,\text{cm}^{-3},
\qquad 
u_1 = u_2 = 2.2\cdot 10^{6}\,\frac{\text{cm}}{\text{s}}, 
\qquad 
T_1 = T_2 = 10\,\text{eV},
\end{align} 
for $x \leq 0$ and
\begin{align}
n_1 = 10^{17}\,\text{cm}^{-3},
\qquad 
n_2 = 10^{20}\,\text{cm}^{-3}, 
\qquad
u_1 = u_2 = -2.2\cdot 10^6 \,\frac{\text{cm}}{\text{s}},
\qquad 
T_1 = T_2 = 10\,\text{eV}
\end{align}
for  $x>0$.

The simulations are run using a velocity grid with $48^3$ nodes and a spatial mesh with 200 cells.  We use the second-order IMEX Runge-Kutta scheme from Section \ref{subsec:secondorderIMEX}  and the second-order spatial discretization in Section \ref{sec:space}, with the limiter given in \eqref{eq:second-order-limiter}. The time step $\Delta t = 806$ fs is set according to the  CFL condition in \eqref{eq:CFL}. 

\begin{figure}


\begin{subfigure}[c]{0.33\textwidth}
\includegraphics[width=\textwidth]{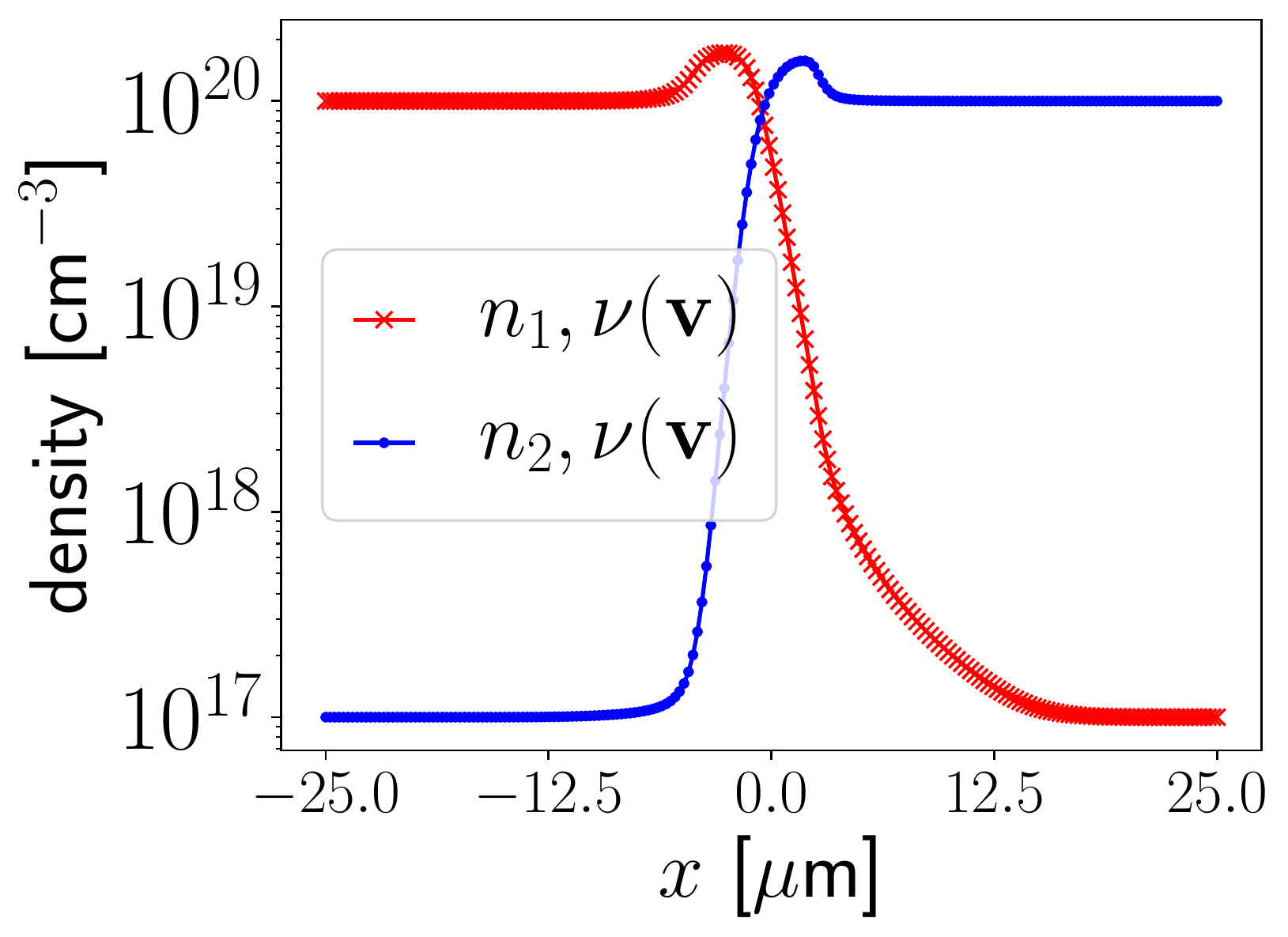}
\end{subfigure}
\begin{subfigure}[c]{0.33\textwidth}
\includegraphics[width=\textwidth]{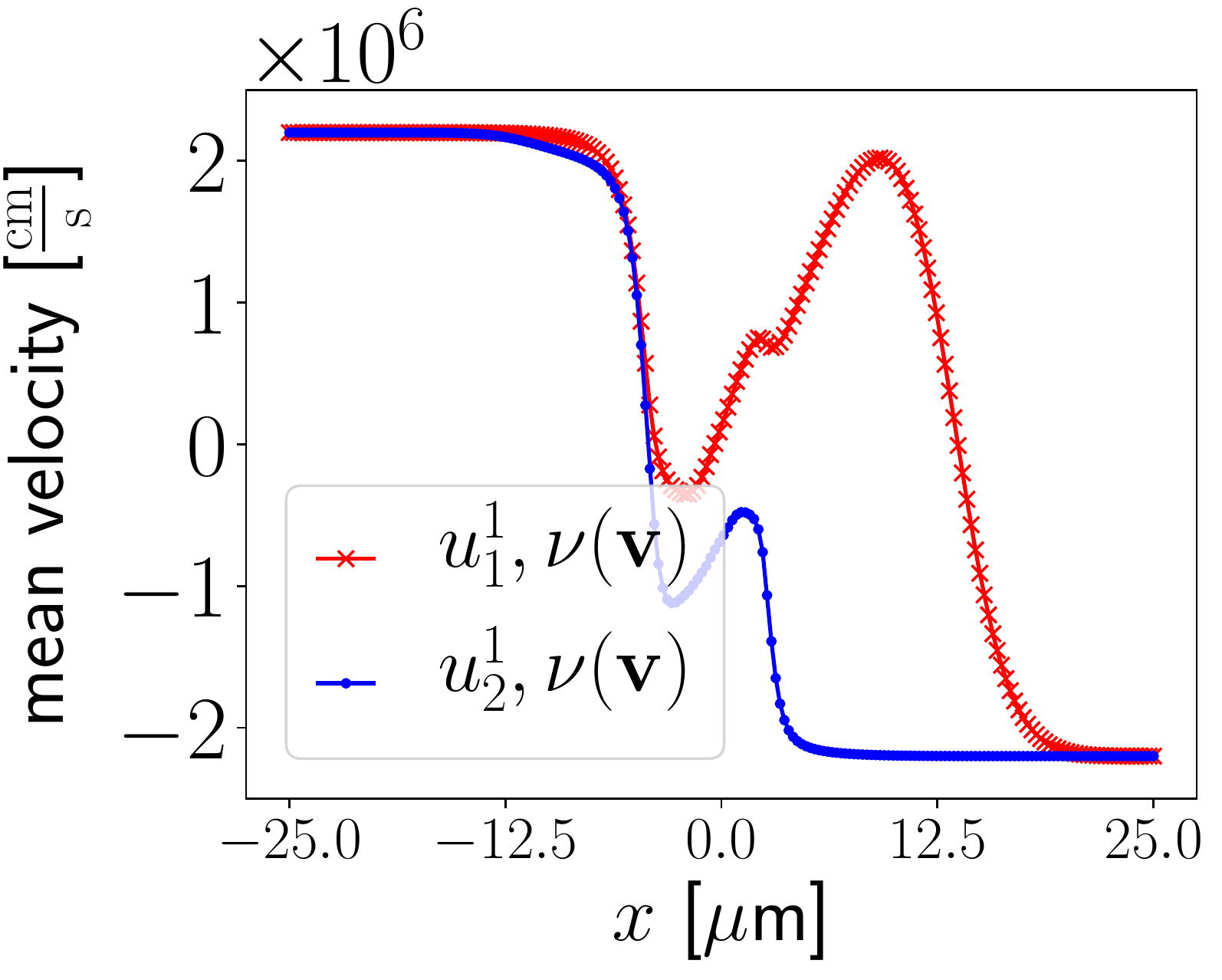}
\end{subfigure}
\begin{subfigure}[c]{0.33\textwidth}
\includegraphics[width=\textwidth]{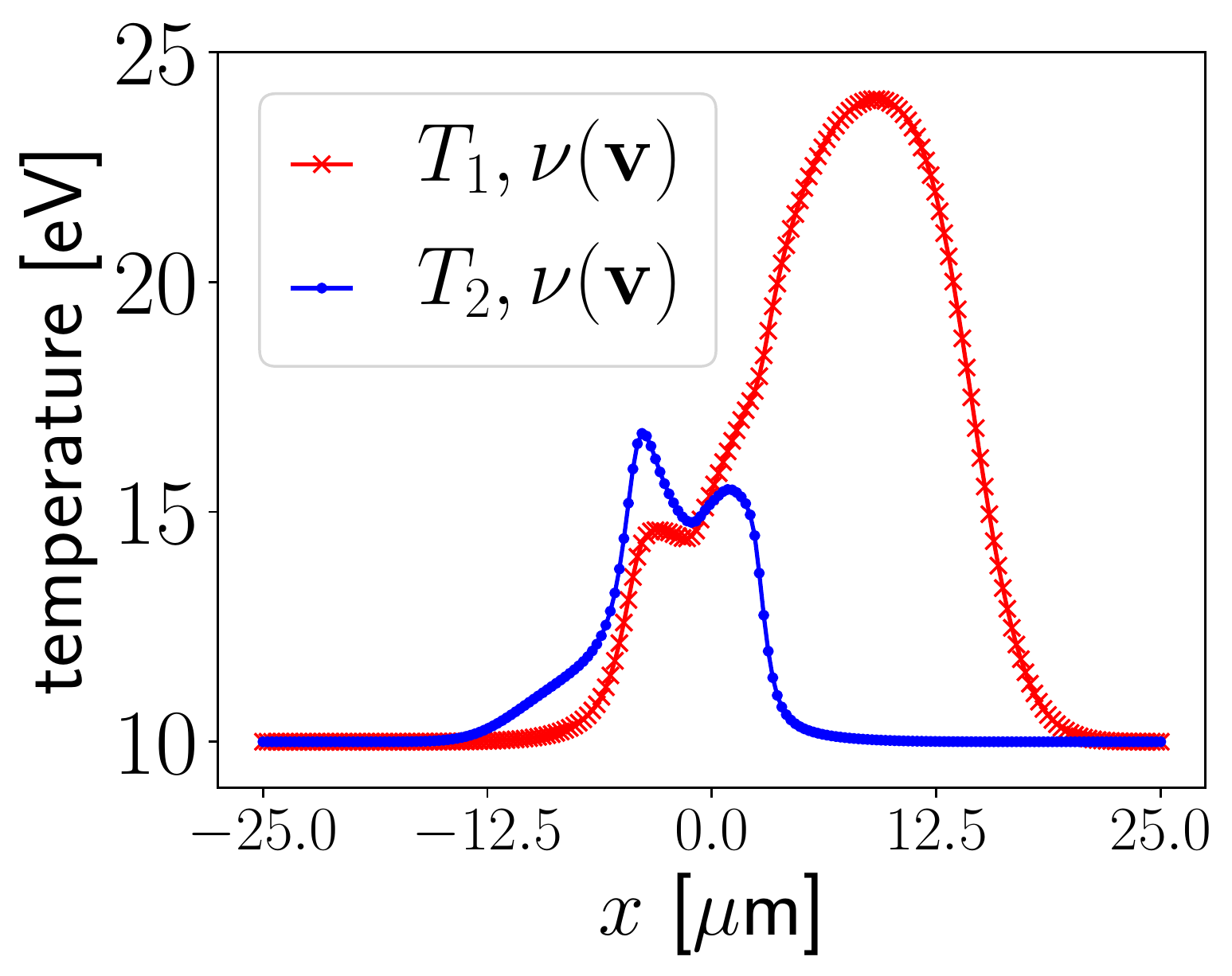}
\end{subfigure}

\begin{subfigure}[c]{0.33\textwidth}
\includegraphics[width=\textwidth]{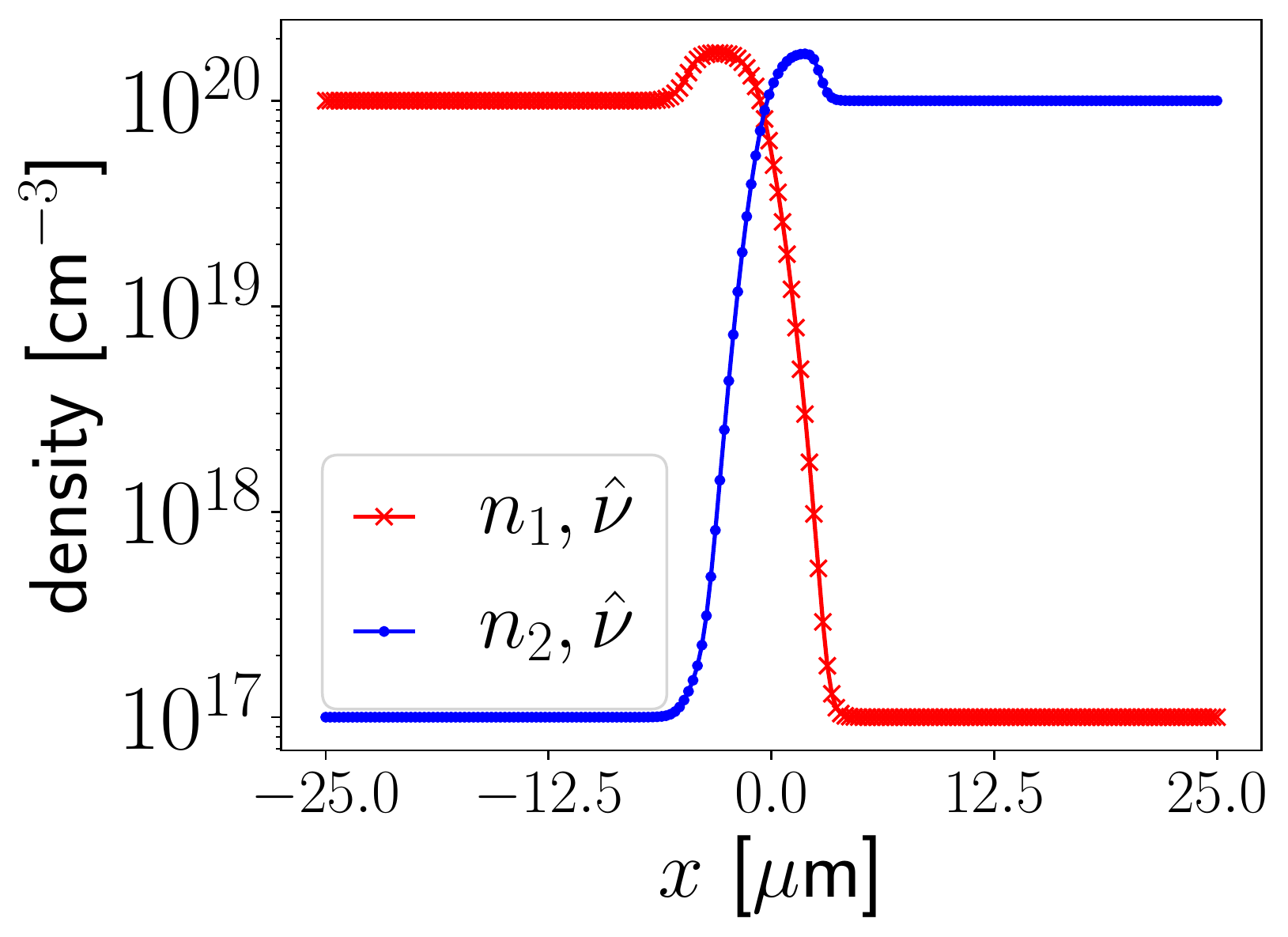}
\end{subfigure}
\begin{subfigure}[c]{0.33\textwidth}
\includegraphics[width=\textwidth]{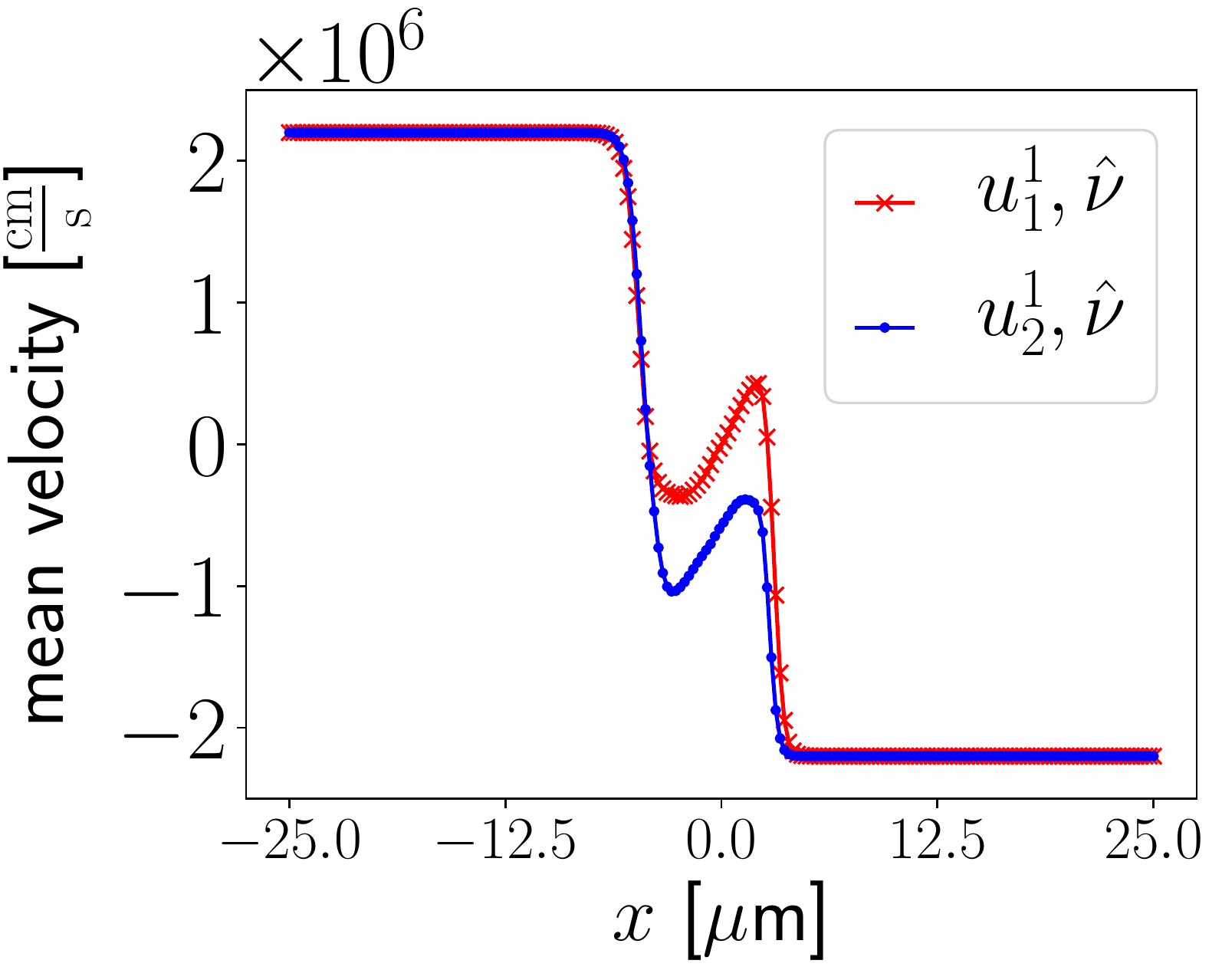}
\end{subfigure}
\begin{subfigure}[c]{0.33\textwidth}
\includegraphics[width=\textwidth]{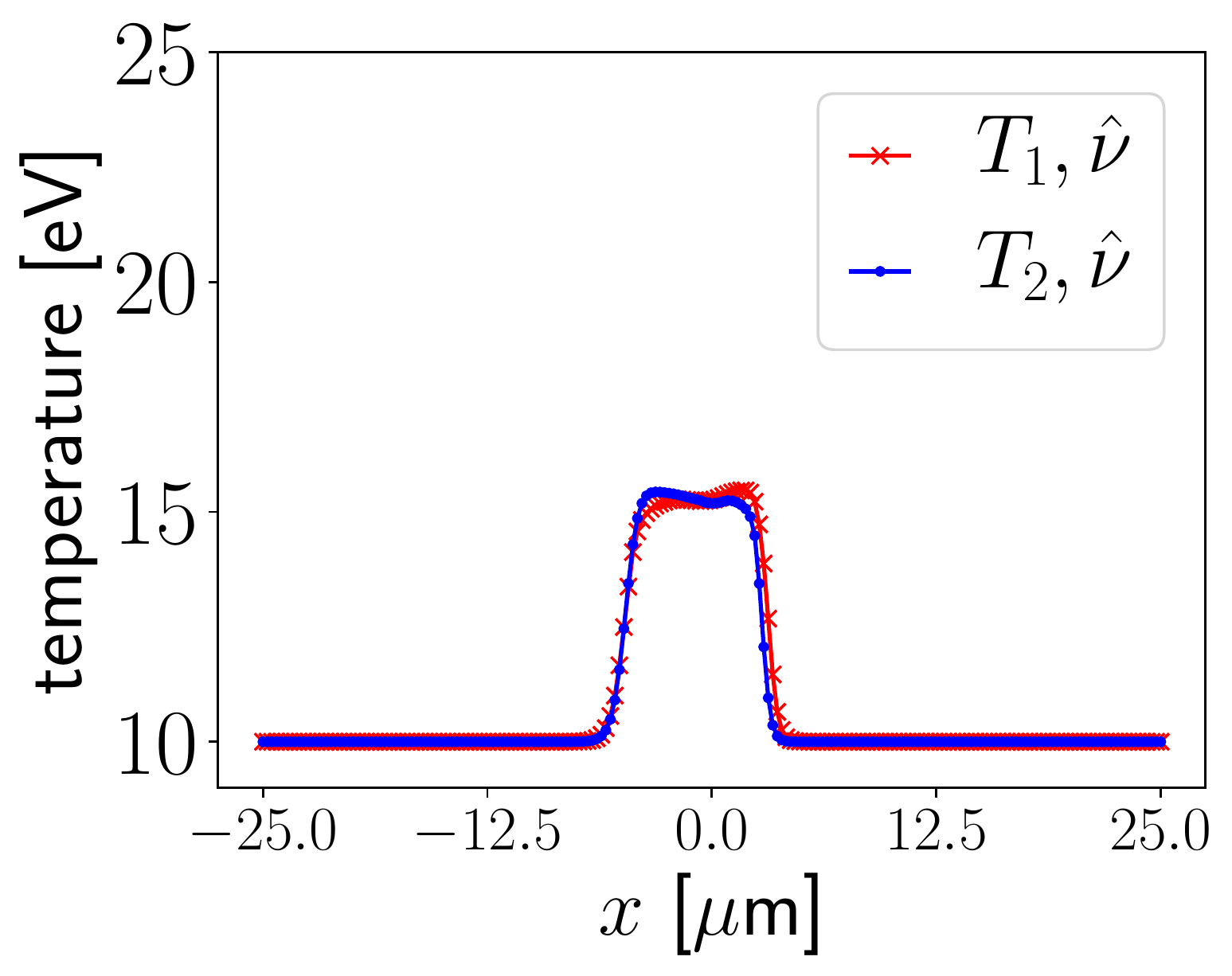}
\end{subfigure}


\begin{subfigure}[c]{0.33\textwidth}
\includegraphics[width=\textwidth]{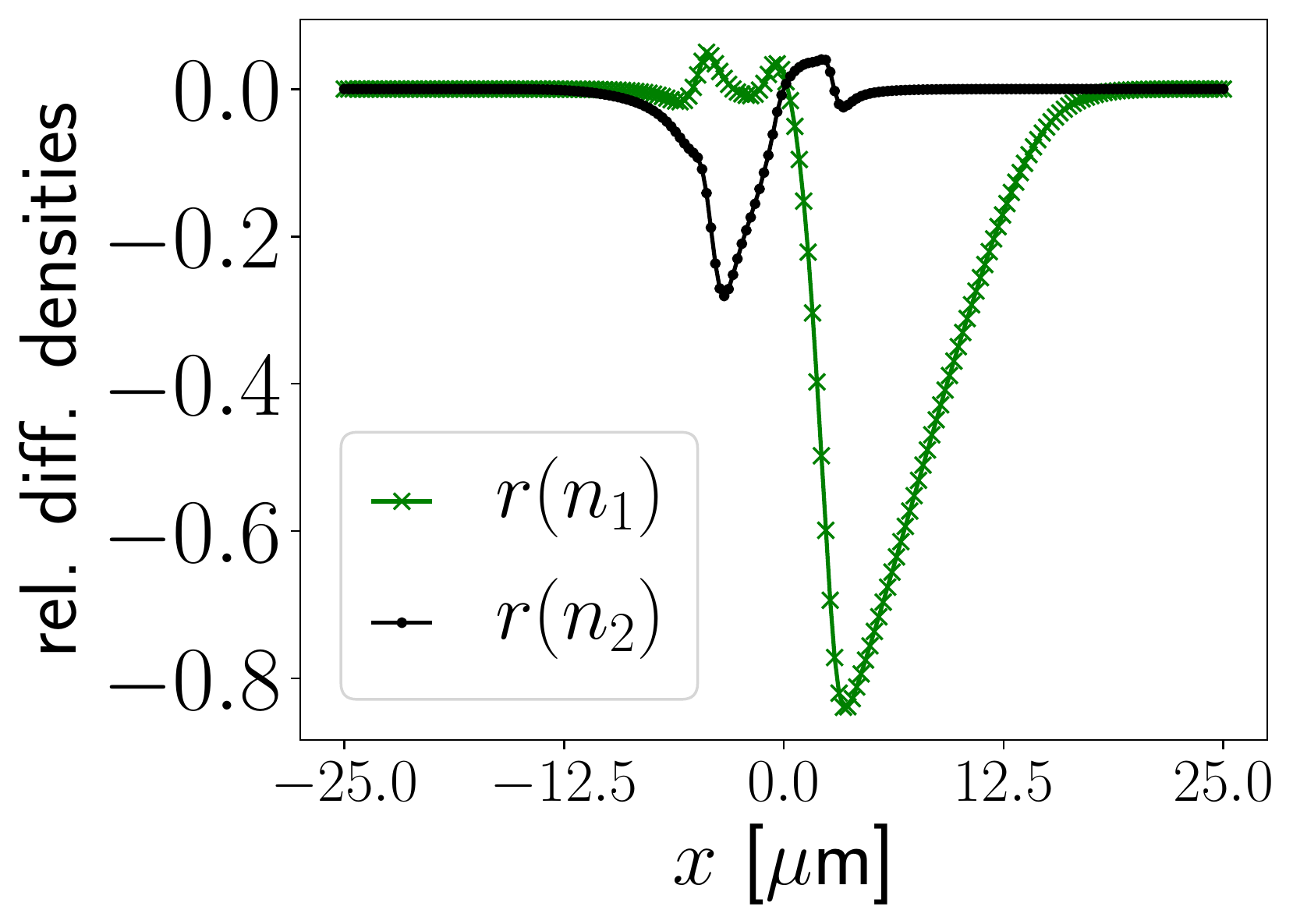}
\end{subfigure}
\begin{subfigure}[c]{0.33\textwidth}
\includegraphics[width=\textwidth]{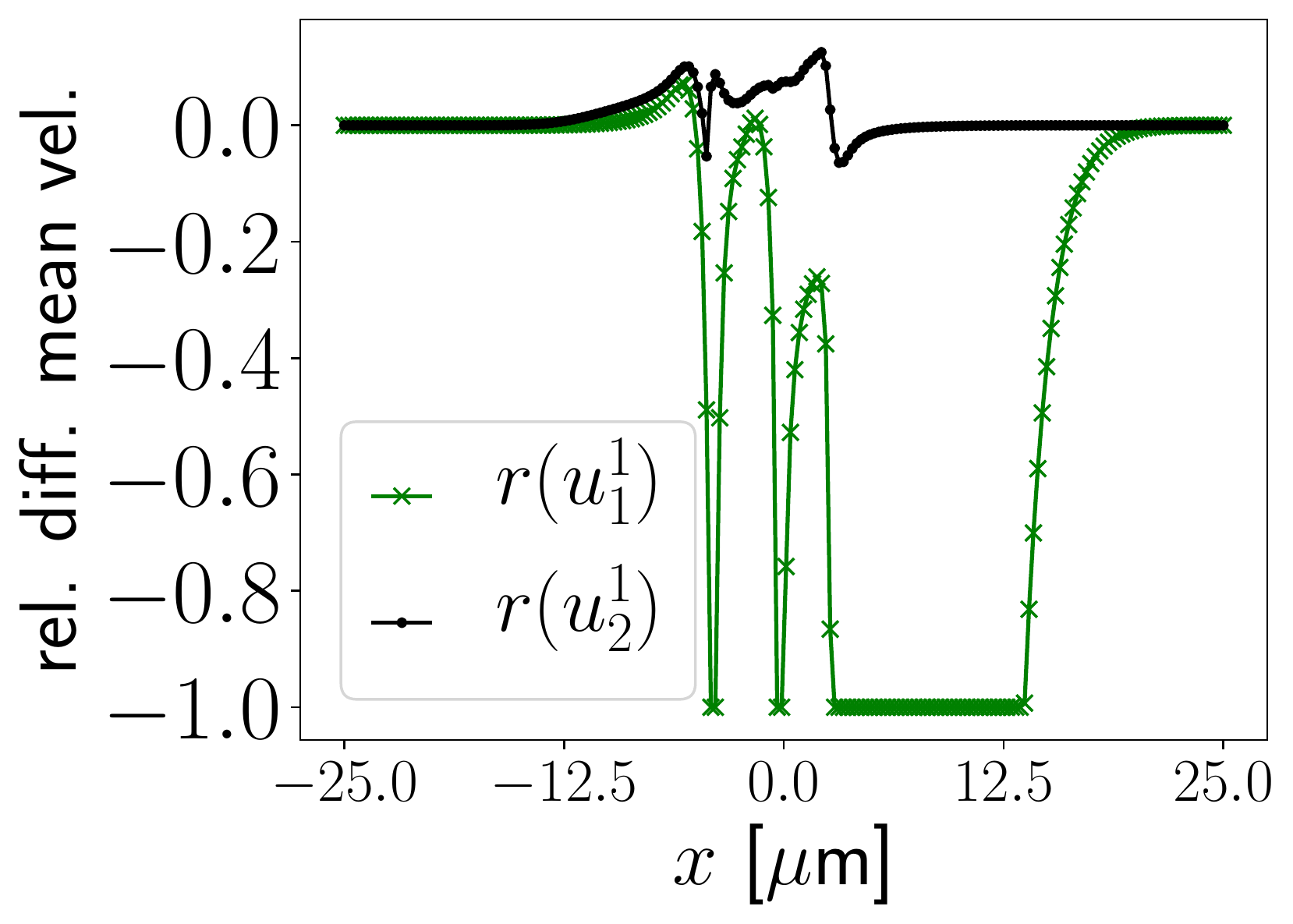}
\end{subfigure}
\begin{subfigure}[c]{0.33\textwidth}
\includegraphics[width=\textwidth]{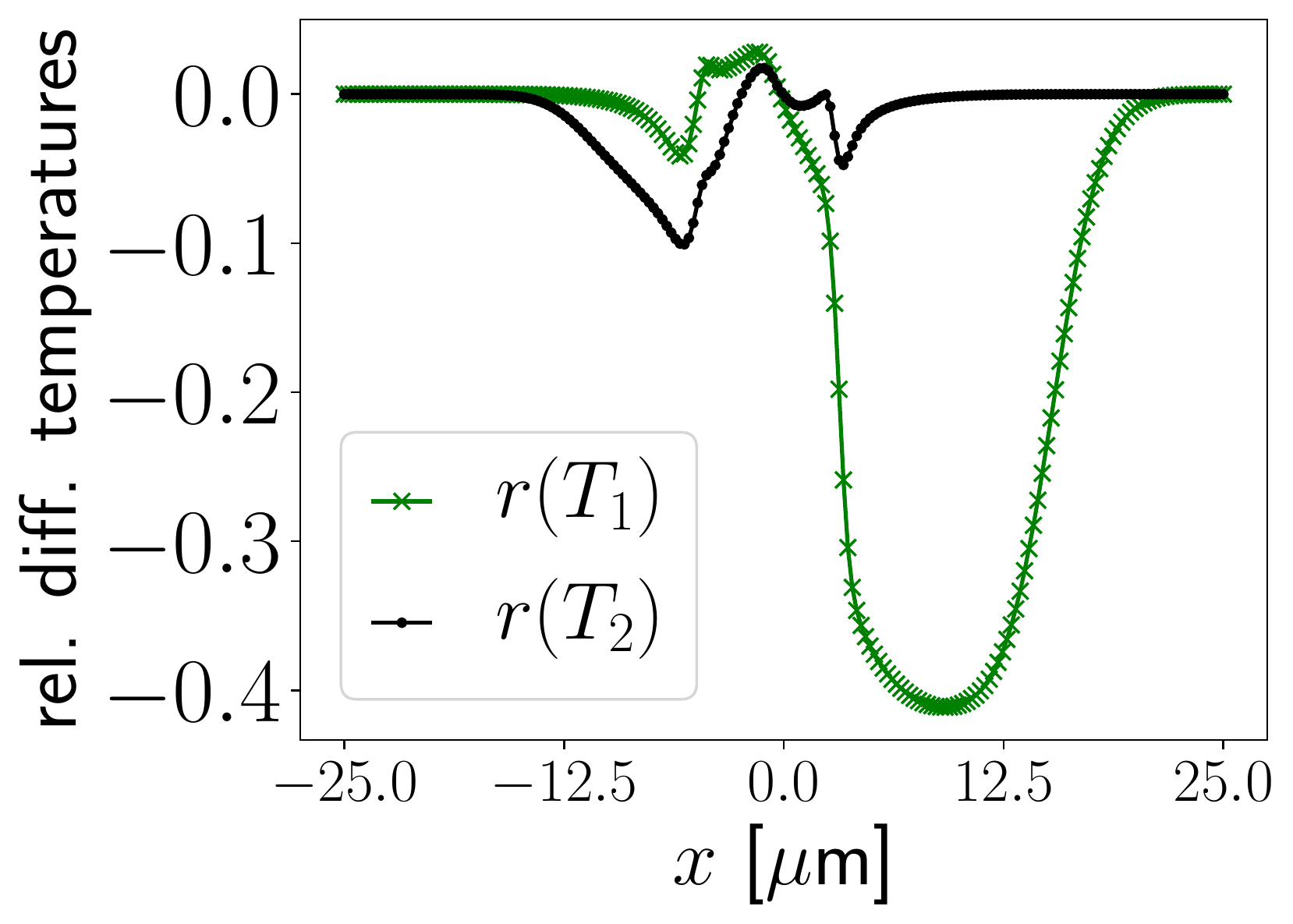}
\end{subfigure}

    \caption{The fluid quantities for the interpenetration problem from Section \ref{test:interpenetration1} are presented at time $t=120.870$ ps. First row: velocity-dependent collision frequencies $\nu_{ij}$, given in \eqref{eq:colfreq_dep}. Second row: the constant collision frequencies $\hat{\nu}_{ij}$, given in \eqref{eq:colfreq_indep_2}. Third row: relative difference between rows 1 and 2 according to \eqref{eq:rel_diff_new}. Red line: hydrogen. Blue line: helium.  Variations in the collision frequency induce significant differences in the profile of the hydrogen, which penetrates much further into the right side of the domain when the collision frequency is velocity-dependent.  Due to relatively higher mass and charge state, the helium species undergoes more collisions and is less sensitive to variations in the collision frequency. }
    \label{fig:interpenetration1}
\end{figure}

We  compare the numerical results  at time $t=120.870$ ps using the velocity-dependent collision frequency $\nu(\mbv)$, given in \eqref{eq:colfreq_dep}, with those using the constant collision frequencies $\hat{\nu}$, given in \eqref{eq:colfreq_indep_2}, in Figure \ref{fig:interpenetration1}. 
The lighter hydrogen species shows a significant difference in profiles between the two species, and displays much more penetration into the helium beam. 
Due to its relatively higher mass and charge state, the helium species is much more collisional than the hydrogen species, and presents a more hydrodynamic-like profile.

\subsubsection{Interpenetration problem: low density} \label{test:interpenetration2}
We repeat the interpenetration problem from above, but reduce the initial densities by two orders of magnitude, which leads to fewer collisions. 
We expect to see a greater interpenetration of the two beams, with less of a density spike at the interface point.
The domain size, masses and charges are the same as before.
The initial conditions are {$f_i = M_i[n_i,\mbu_i,T_i]$ with}: 
\begin{align}
n_1 = 10^{18}\,\text{cm}^{-3},
\qquad 
n_2 = 10^{15}\,\text{cm}^{-3},
\qquad 
u_1 = u_2 = 2.2\cdot 10^{6}\,\frac{\text{cm}}{\text{s}}, 
\qquad 
T_1 = T_2 = 10\,\text{eV},
\end{align} 
for $x \leq 0$ and
\begin{align}
n_1 = 10^{15}\,\text{cm}^{-3},
\qquad 
n_2 = 10^{18}\,\text{cm}^{-3}, 
\qquad
u_1 = u_2 = -2.2\cdot 10^6 \,\frac{\text{cm}}{\text{s}},
\qquad 
T_1 = T_2 = 10\,\text{eV}
\end{align}
for  $x>0$.

As before, the simulations are run using a velocity grid with $48^3$ nodes and a spatial mesh with 200 cells.  We use the second-order IMEX Runge-Kutta scheme from Section \ref{subsec:secondorderIMEX}  and the second-order spatial discretization in Section \ref{sec:space}, with the limiter given in \eqref{eq:second-order-limiter}. The time step $\Delta t = 806$ fs is set according to the  CFL condition in \eqref{eq:CFL}. 

\begin{figure}


\begin{subfigure}[c]{0.33\textwidth}
\includegraphics[width=\textwidth]{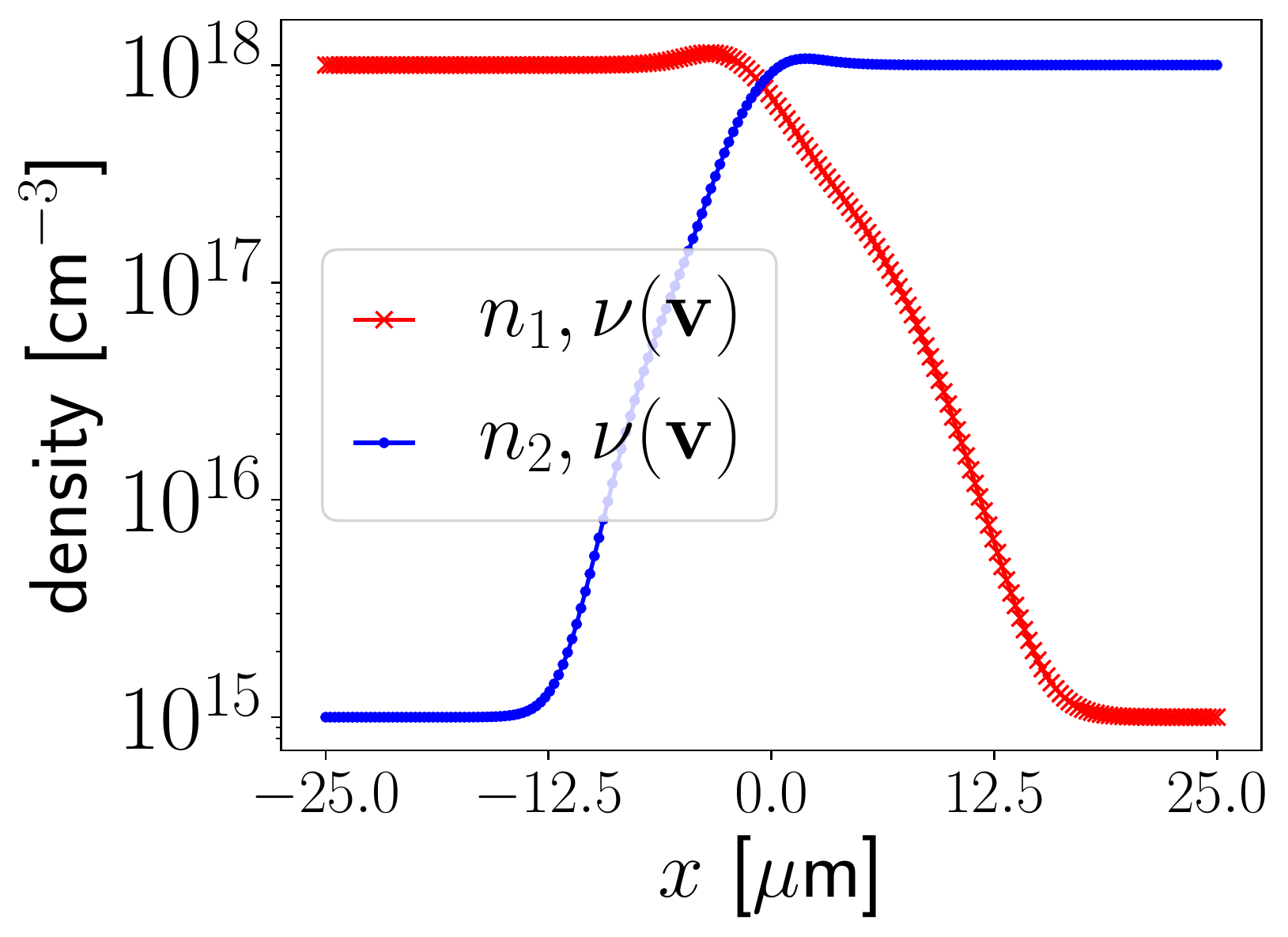}
\end{subfigure}
\begin{subfigure}[c]{0.33\textwidth}
\includegraphics[width=\textwidth]{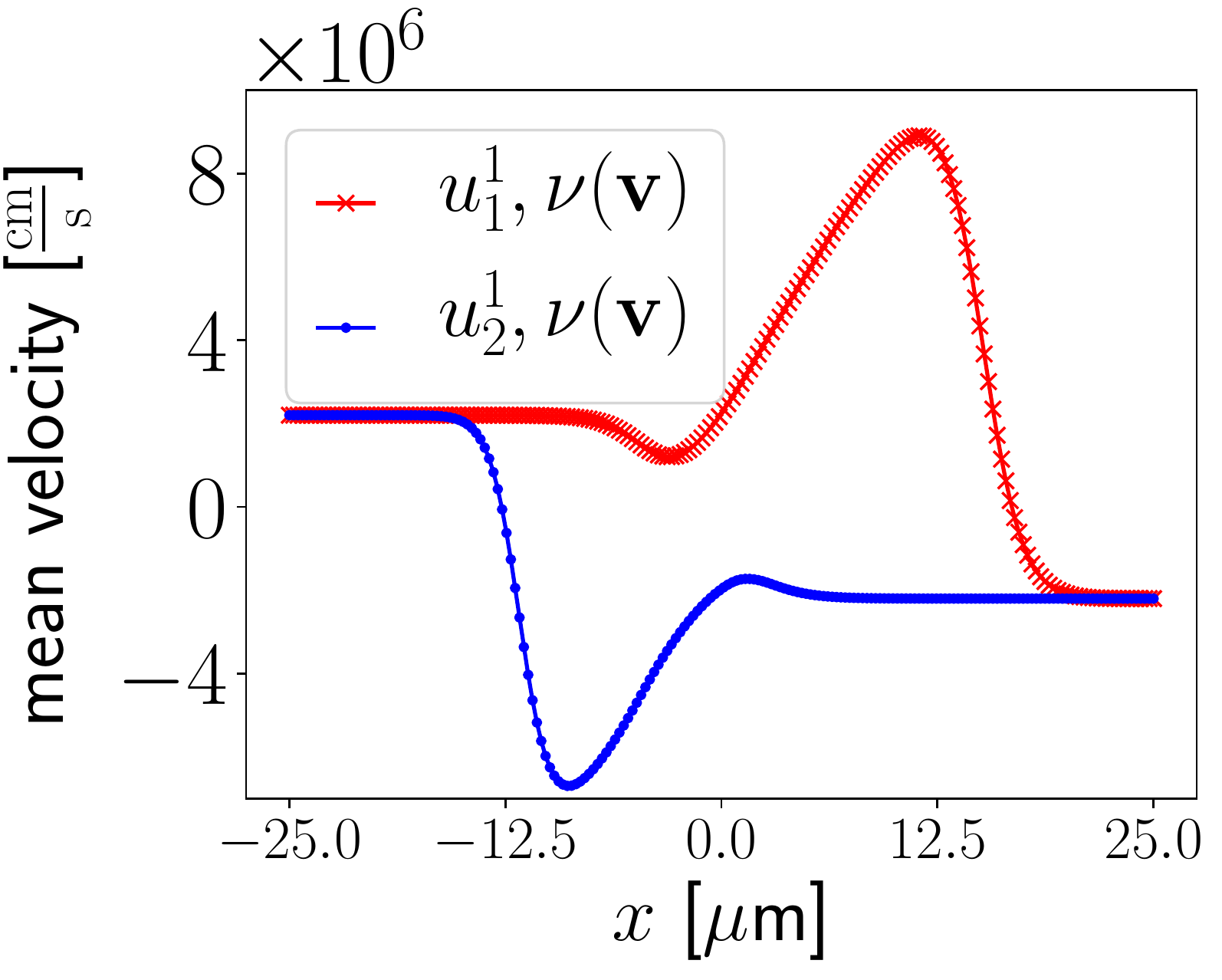}
\end{subfigure}
\begin{subfigure}[c]{0.33\textwidth}
\includegraphics[width=\textwidth]{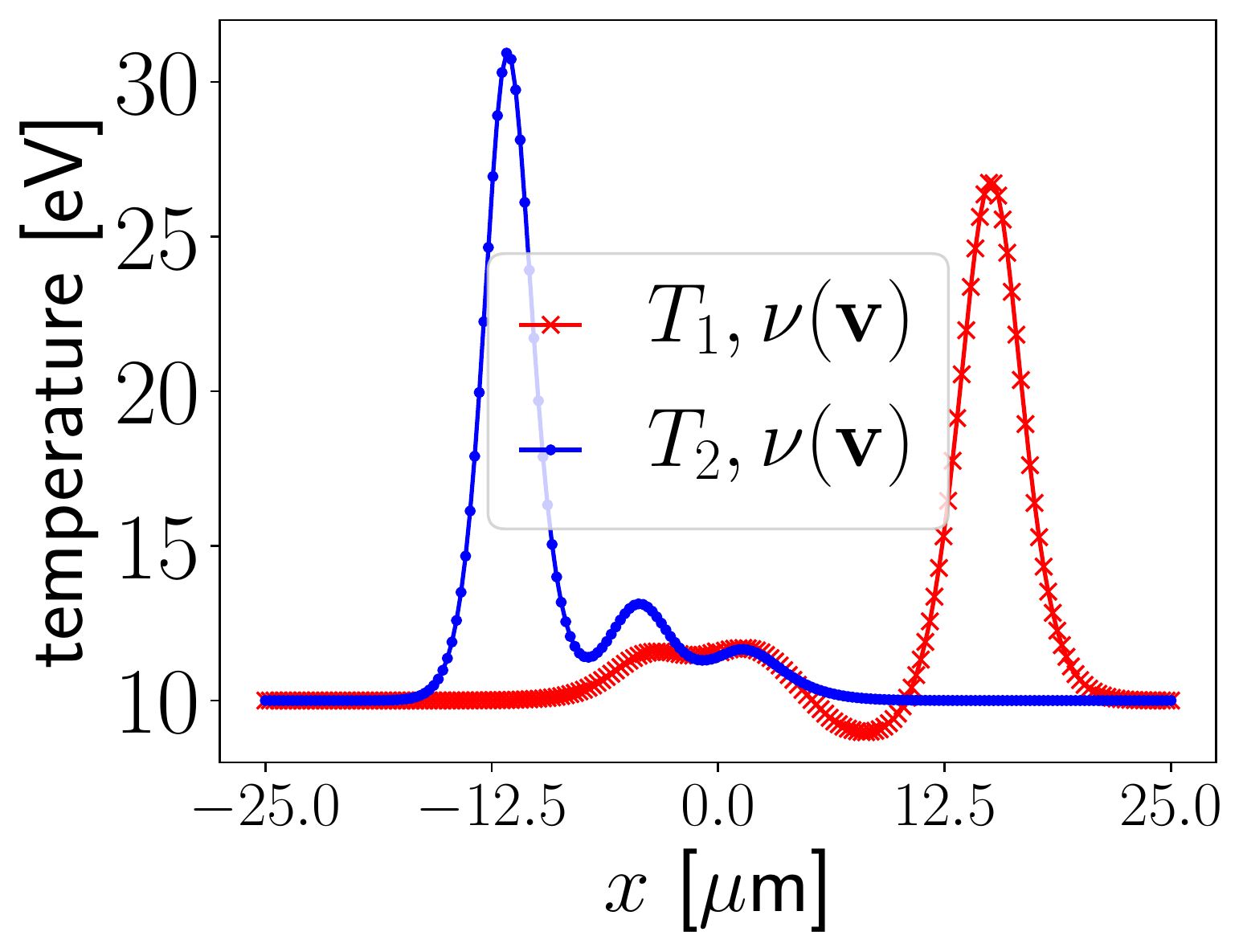}
\end{subfigure}

\begin{subfigure}[c]{0.33\textwidth}
\includegraphics[width=\textwidth]{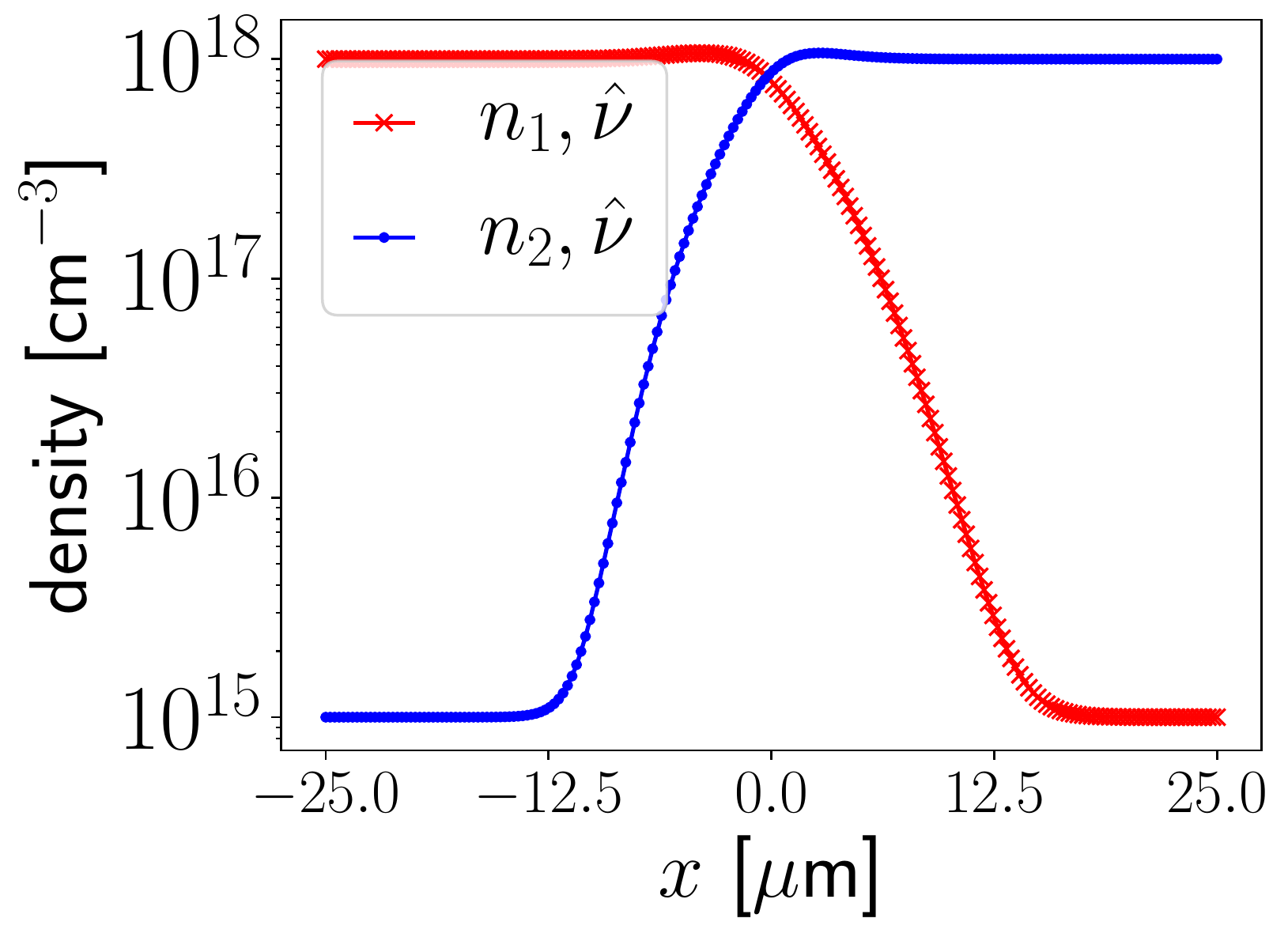}
\end{subfigure}
\begin{subfigure}[c]{0.33\textwidth}
\includegraphics[width=\textwidth]{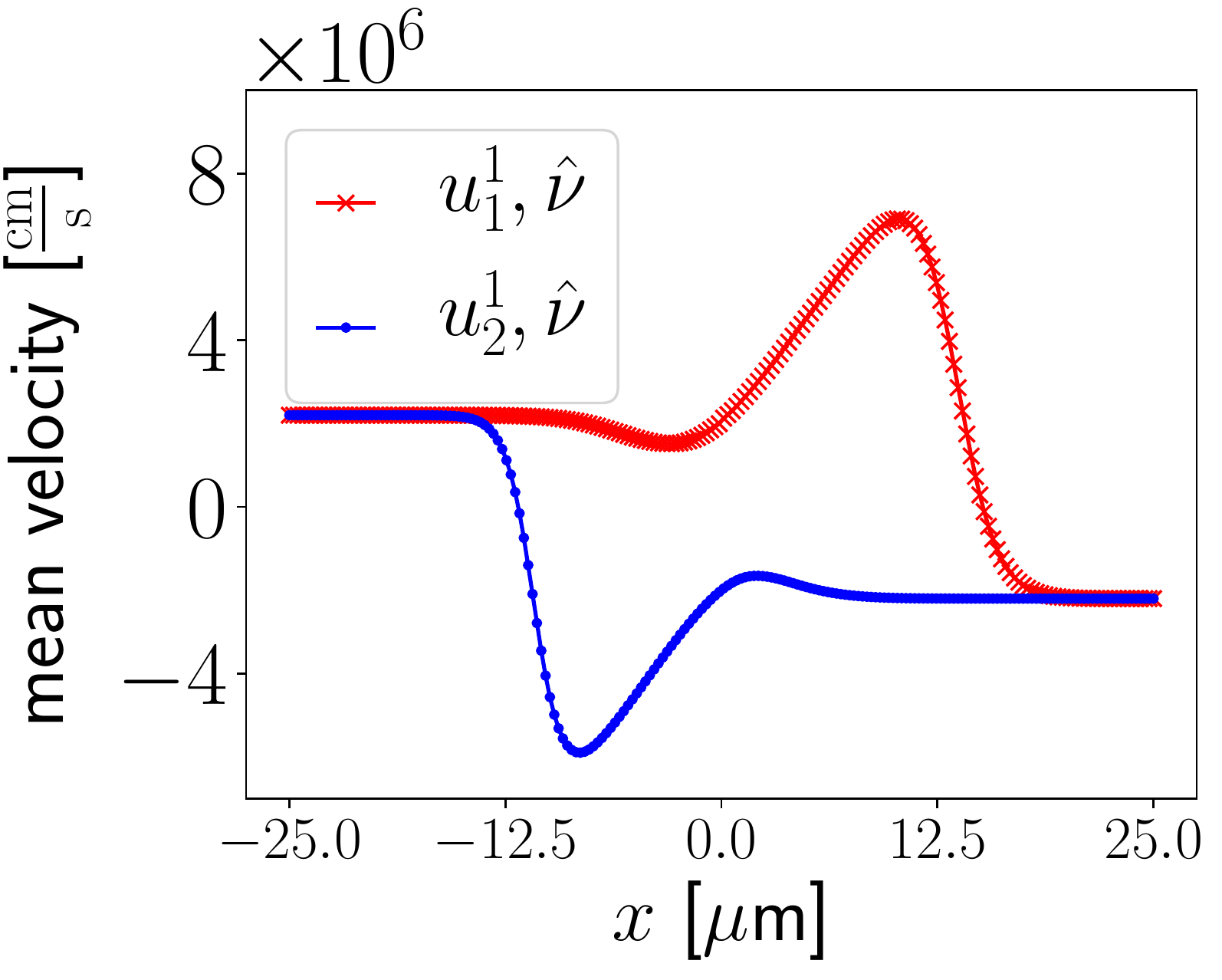}
\end{subfigure}
\begin{subfigure}[c]{0.33\textwidth}
\includegraphics[width=\textwidth]{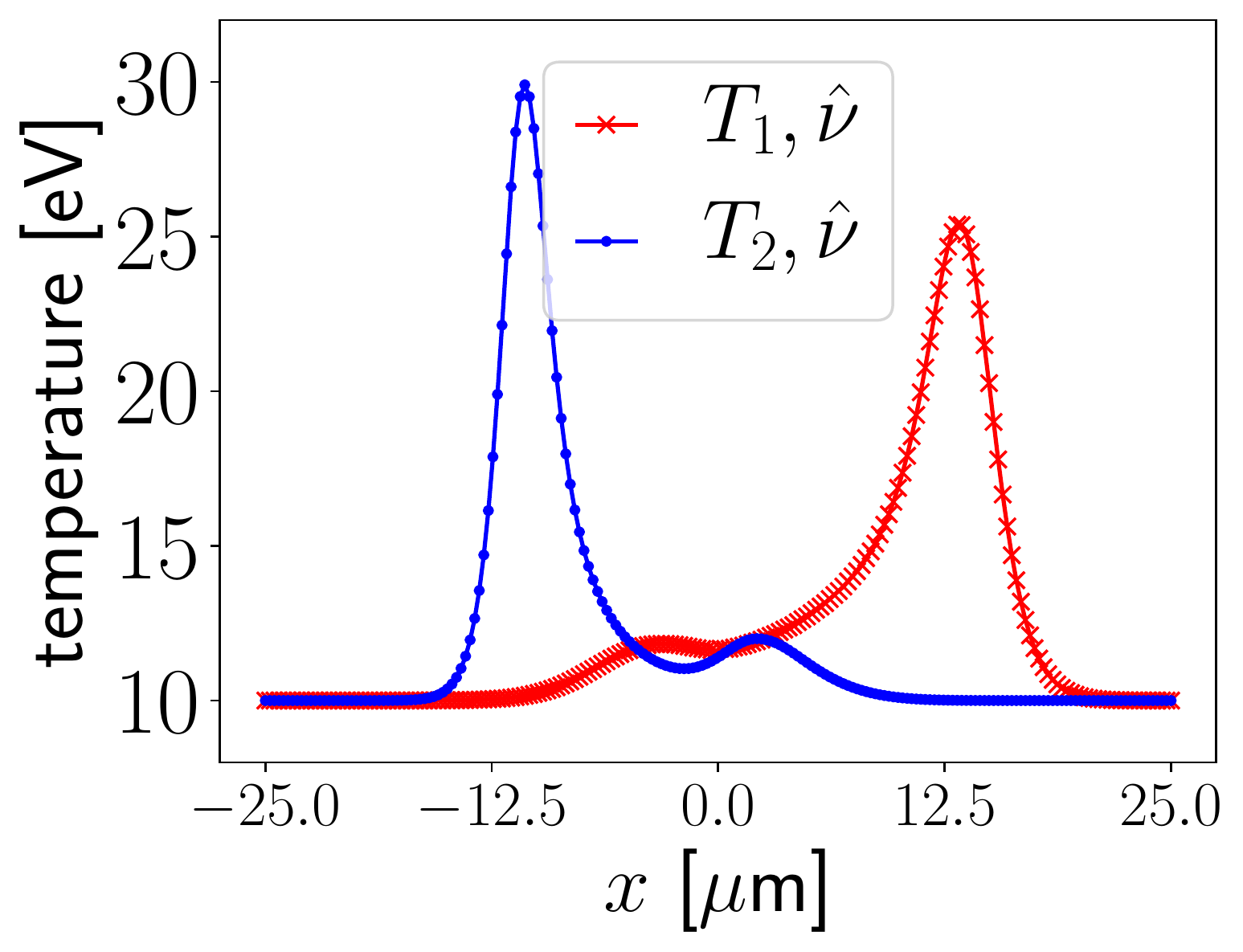}
\end{subfigure}


\begin{subfigure}[c]{0.33\textwidth}
\includegraphics[width=\textwidth]{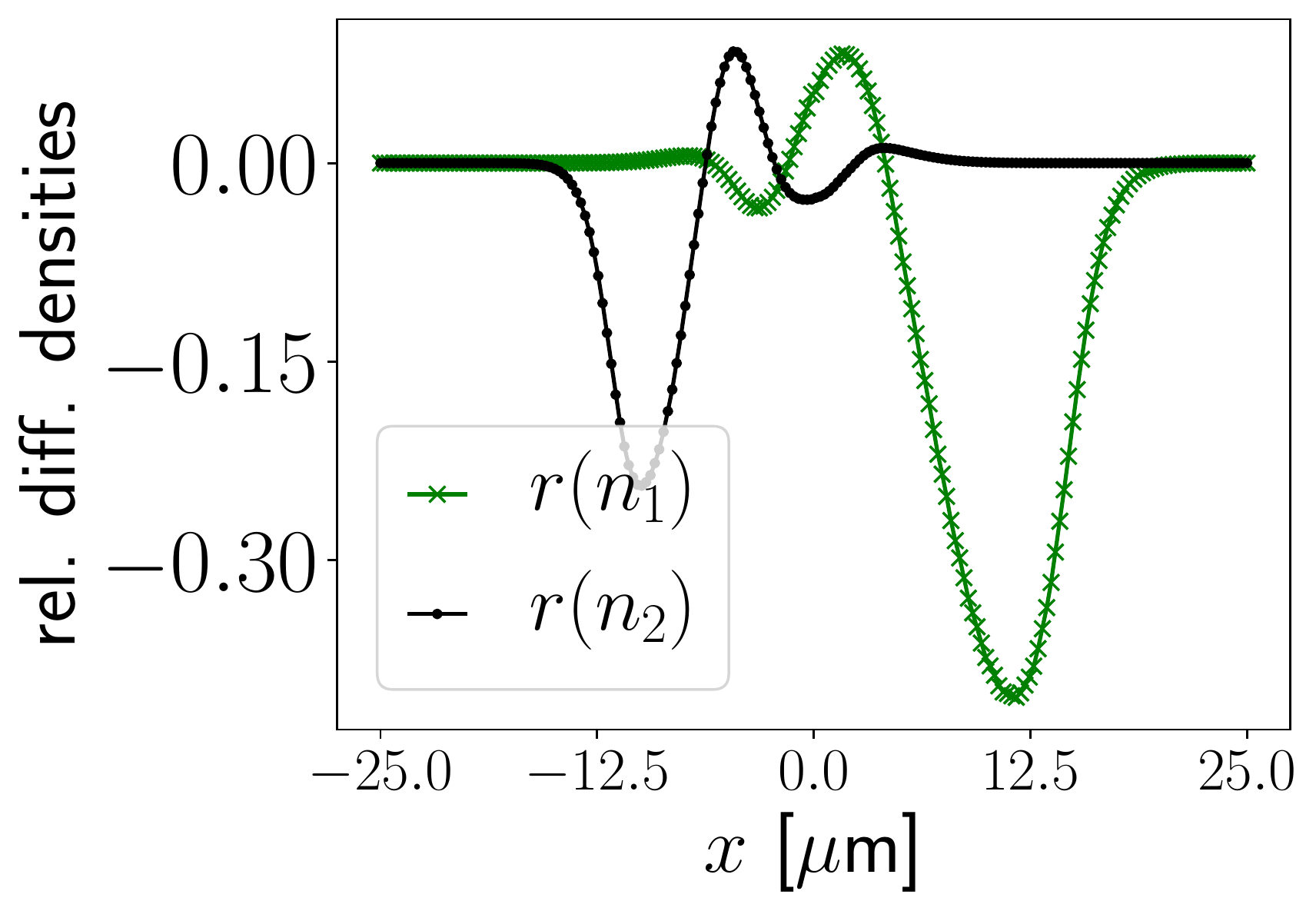}
\end{subfigure}
\begin{subfigure}[c]{0.33\textwidth}
\includegraphics[width=\textwidth]{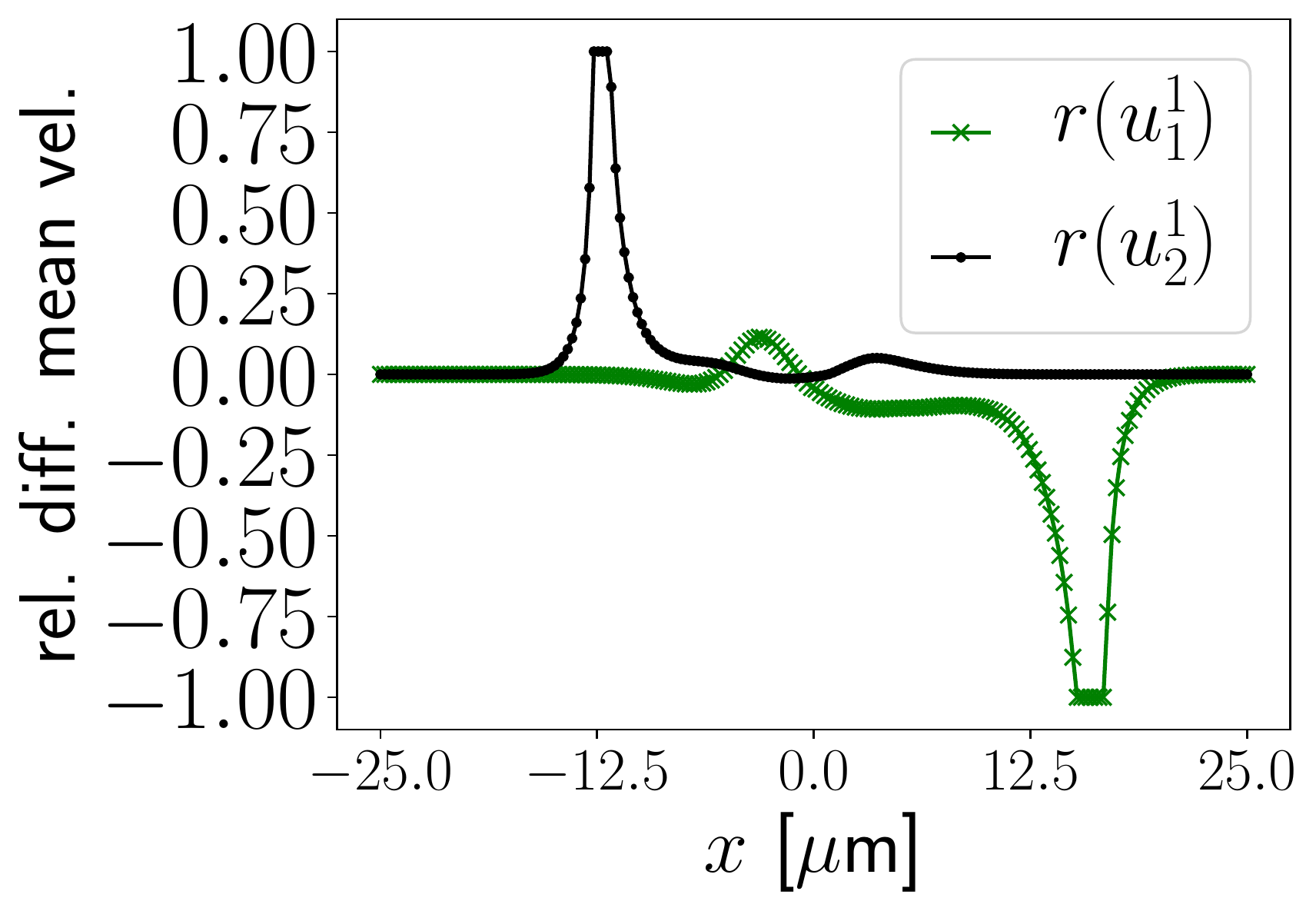}
\end{subfigure}
\begin{subfigure}[c]{0.33\textwidth}
\includegraphics[width=\textwidth]{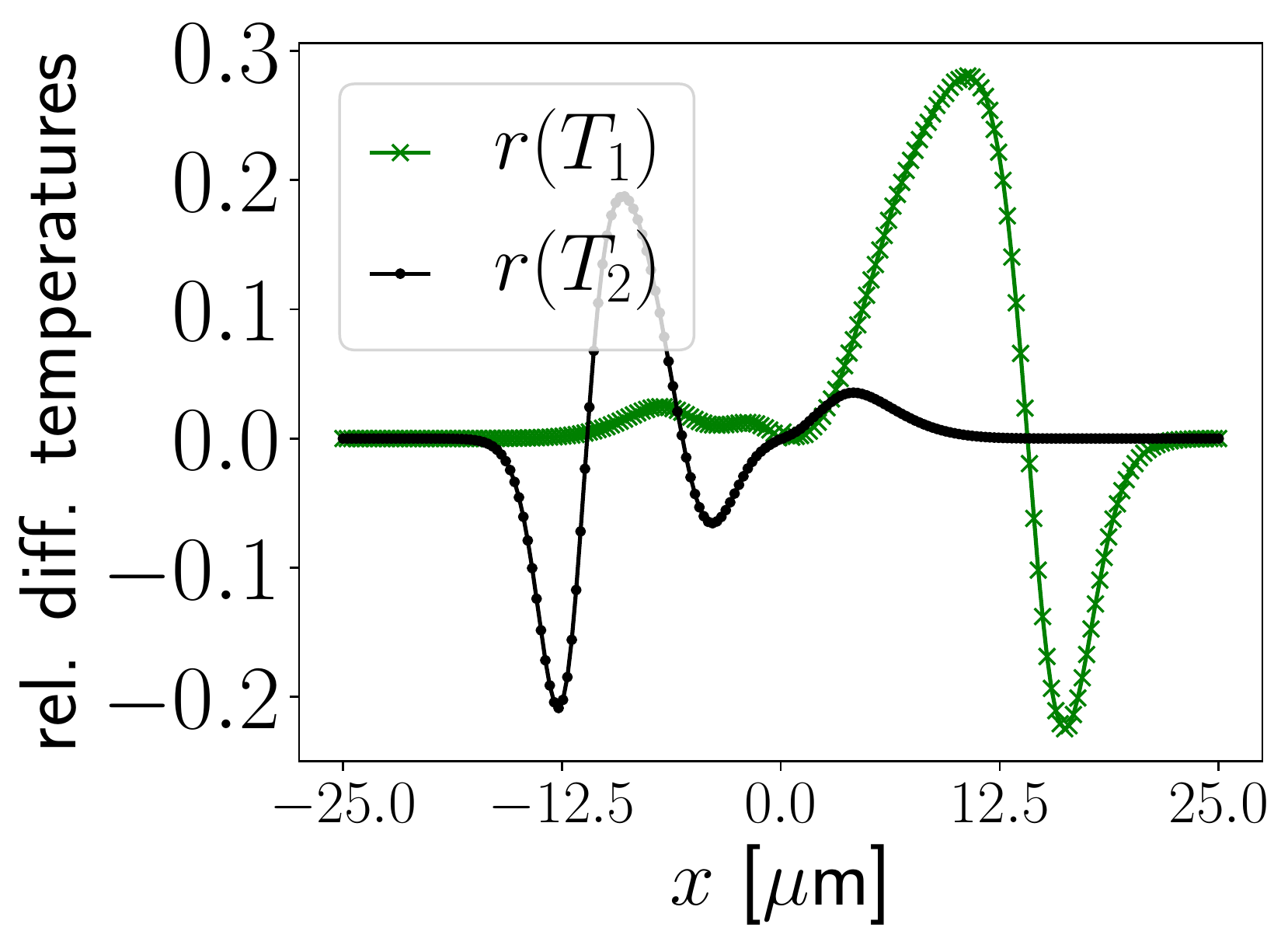}
\end{subfigure}

    \caption{The fluid quantities for the lower density interpenetration problem from Section \ref{test:interpenetration2} are presented at time $t=120.870$ ps. First row: velocity-dependent collision frequencies $\nu_{ij}$, given in \eqref{eq:colfreq_dep}. Second row: the constant collision frequencies $\hat{\nu}_{ij}$, given in \eqref{eq:colfreq_indep_2};. Third row: relative difference between rows 1 and 2 according to \eqref{eq:rel_diff_new}. Red line: hydrogen. Blue line: helium.  As expected, we see more interpenetration than in the high density test case. However, the relative sensitivity of hydrogen to the velocity-dependent collision frequency is less dramatic.
    } 
    \label{fig:interpenetration2}
\end{figure}

We  compare the numerical results  at time $t=120.870$ ps using the velocity-dependent collision frequency $\nu(\mbv)$, given in \eqref{eq:colfreq_dep}, with those using the constant collision frequencies $\hat{\nu}$, given in \eqref{eq:colfreq_indep_2}, in Figure \ref{fig:interpenetration2}. As expected, we see more interpenetration than in the high density test case. As in the higher density test case, we see more significant differences in the lighter species of the mixture; the hydrogen species penetrates more into the right side of the domain when the collision frequency is velocity-dependent. Due to relatively higher mass and charge state, the helium species is more collisional.
Furthermore, the density spike at the interface seen in the high density case has mostly disappeared. 

\section{Conclusions}
\label{sec:conclusions}

We have developed a numerical scheme for the multi-species BGK model with velocity-dependent collision frequency, first proposed in \cite{HHKPW}. The main new contribution is the implicit update in an IMEX formulation.  The dependence of the target functions on the distribution function is only known implicitly for general collision frequencies so that standard approaches for BGK models cannot be used. We find the target function via a convex minimization problem that mimics the dual of the minimization problem that defines the theoretical model in \cite{HHKPW}.  This procedure automatically satisfies a discrete version of the conservation laws for mass, total momentum, and total energy satisfied by the BGK operator. The transport part is discretized by a standard finite volume method. For a first-order scheme, we verify that a discrete entropy dissipation property and positivity of the distribution function hold rigorously.  A second-order version of the method is used for improved accuracy.

We illustrate the properties of the BGK model and our numerical scheme with several test cases, using velocity-dependent collision frequencies that are motivated by Coulomb interactions in plasmas and characterized by slower relaxation in the tails of the kinetic distribution. The simulation results are compared to results with velocity-independent collision frequencies of comparable size. For spatially homogeneous problems, the velocity-dependent collision frequencies induce slower relaxation to equilibrium in the tails of the kinetic distributions.  The convergence of the temperature and mean velocities is also slower and has a distinctly different form than in the velocity-independent case.  

Several Riemann problems are also considered, including the standard Sod shock-tube problem and variations involving mixtures.  In the former, we confirm that the BGK model recovers the general fluid shock structure, but the kinetic effects are more apparent in the case of velocity-dependent collision frequencies. For the mixtures, we observe close agreement between simulations using velocity-dependent and velocity-independent collision frequencies for a Mach 1.7 shock and for a Mach 4 shock, with deviations approaching 2\%. For the interpenetration problems, the profiles differ more significantly. In particular, the effect of the velocity-dependent collision frequencies on the lighter species (in mass and charge state) are substantial.

Allowing for velocity-dependent collision frequencies is not without additional cost.  Indeed, for the velocity-dependent frequencies, the implicit evaluation of the collision operator requires the solution of an optimization problem via a Newton solver.  In particular, the elements of the gradient and Hessian of the objective function  require the evaluation of the integrals in velocity space via a quadrature.  To accelerate the solution procedure, a more efficient implementation of the optimization algorithm is necessary \cite{schaerer2017efficient,AlldredgeHauckOLearyTits2014,AlldredgeHauckTits2012,GarrettHauckHill2015,Abramov2007}.  In spite of the additional cost, the model still has better scaling properties than the original Boltzmann equation. 

The enlarged class of possible collision frequencies is physically motivated, and  makes the extended multi-species BGK model an attractive option for exploring more phenomena in the kinetic regime. However, the model can be improved. For example, when considering charged particles, the model needs to include a force term with an electric (and magnetic) field. The additional transport in velocity space can be easily incorporated in the presented numerical method. 


\section*{Appendix}
We illustrate the convergence behavior of the schemes we presented in Sections \ref{sec:time} and \ref{sec:space}.
For the most part, we use methods that are already established in the literature. More precisely, the basic time discretization techniques can be found in \cite{ARS97} and the numerical fluxes for the spatial discretization are taken from  \cite{MieussensStruchtrup2004}.
The new ingredient is a general implicit solver for determining the target
functions. This can be applied to many different time discretization techniques, provided that the updates can be written in the form of \eqref{eq:update_general}. 

As an illustration, we consider a test case inspired by \cite{hu2021uniform}. In this problem, the spatial domain is $[0,2]$, and both species are initialized in the same way: $f_i(x)=M_i[n_i(x),\mbu_i(x),T_i(x)]$, where 
\begin{align}
    n_1(x)=n_2(x)=1+0.1\sin(\pi x), 
    \quad 
    u_1(x) = u_2(x) =1, 
    \quad
    T_1(x)= T_2(x)\frac{1}{1+0.1\sin(\pi x)},
\end{align}
and the species mass are $m_1 = m_2=1$. 

The simulations are run using a velocity grid with $48^3$ nodes. We incorporate the collision frequencies
\begin{equation} \label{eq:nu-order}
 \nu_{ij}(x,\mbv,t) = \frac{C\, n_j}{\delta_{ij} +\vert \mbv-\mbu_\mix \vert^{3} },
\end{equation}
with the regularization parameter $\delta_{ij} = 0.1 \cdot (\Delta v_{ij})^3$, where $\Delta v_{ij} = \frac{1}{4} \sqrt{T_\mix/(2\mu_{ij})}$ and $\mu_{ij} = m_im_j/(m_i+m_j)$ and $C=1$ or $C=10^4$.

The convergence rates are determined on consecutive meshes with $2^l\cdot 20$ cells, $l=0,\dots,4$. We combine (i) the first-order temporal splitting scheme from Section \ref{subsec:firstordersplit} with the first-order spatial discretization in Section \ref{sec:space} and (ii) the second-order IMEX Runge-Kutta scheme from Section \ref{subsec:secondorderIMEX} with the second-order spatial discretization in Section \ref{sec:space}, with the limiter given in \eqref{eq:second-order-limiter}. In both cases, the time step is set according to the CFL condition in \eqref{eq:CFL}.

\begin{figure}[htb]
    \begin{subfigure}[c]{0.49\textwidth}
    \includegraphics[width=0.9\textwidth]{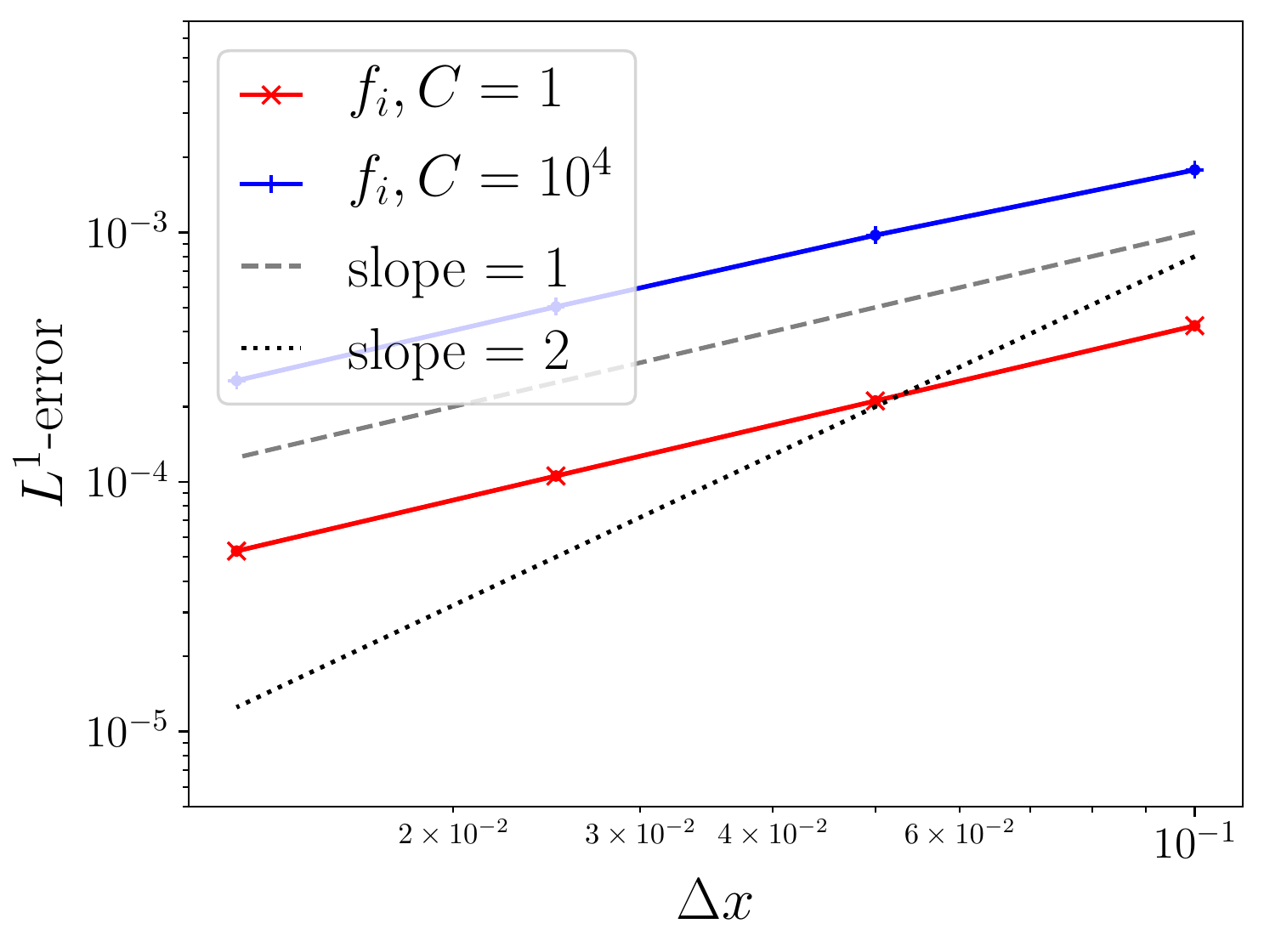}
    \caption{first-order method \ref{subsec:firstordersplit} }
    \end{subfigure}
    \begin{subfigure}[c]{0.49\textwidth}
    \includegraphics[width=0.9\textwidth]{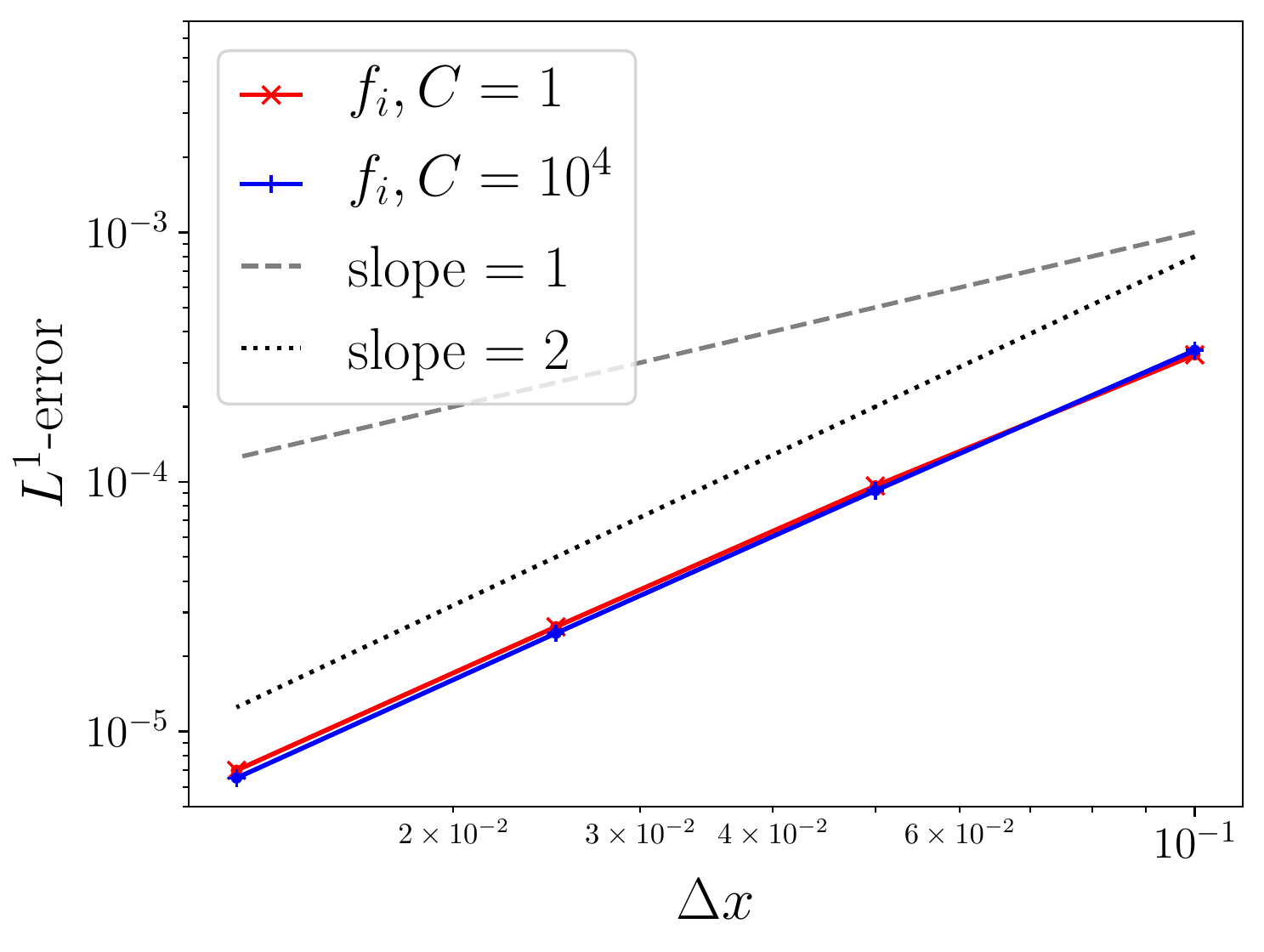}
    \caption{second-order method \ref{subsec:secondorderIMEX} }
    \end{subfigure}
    \caption{Space-time convergence errors measured in the  $L^1$-norm defined in \eqref{eq:L1-error}. The first-order method (left) uses the first-order temporal splitting scheme from Section \ref{subsec:firstordersplit} with the first-order spatial discretization in Section \ref{sec:space}.  The second-order method (right) combines  the second-order temporal IMEX Runge-Kutta scheme from Section \ref{subsec:secondorderIMEX} with the second-order spatial discretization in Section \ref{sec:space}, with the limiter given in \eqref{eq:second-order-limiter}.  The results for $f_1$ and $f_2$ are exactly the same such that the graphs lie on top of each other. The size of the collision frequency \eqref{eq:nu-order}  depends on the factor $C=1$ or $C=10^4$, which does not affect the order of the scheme. }
    \label{fig:order}
\end{figure}
In Figure \ref{fig:order}, the $L^1$-error of the solution $f_i(\Delta t, \Delta x)$ is estimated by 
\begin{align} \label{eq:L1-error}
     ||f_i(\Delta t, \Delta x)-f_i(\Delta t/2, \Delta x/2)||_{L^1_{x,v}}\approx \sum_k \sum_\mbq |f_{i, k\mbq}(\Delta t, \Delta x)-f_{i,k\mbq}(\Delta t/2, \Delta x/2)| \,(\Delta v)^3 \Delta x.
\end{align}
 In both cases, we observe convergence at a rate that agrees with the formal order of the method.

\section*{Acknowledgements} 

We thank Bruno Despr\'{e}s and Michael Murillo for helpful discussions.

We acknowledge support for covering travel costs for Christian Klingenberg and Sandra Warnecke by the Bayerische Forschungsallianz (grant no. BaylntAn UWUE 2019-29). 

The work of Jeff Haack was supported by the U.S. Department of Energy through the Los Alamos National Laboratory. Los Alamos National Laboratory is operated by Triad National Security, LLC, for the National Nuclear Security Administration of U.S. Department of Energy (Contract No. 89233218CNA000001). LA-UR-22-20147.

The work of Cory Hauck is sponsored by the Office of Advanced Scientific Computing Research, U.S. Department of Energy, and performed at the Oak Ridge National Laboratory, which is managed by UT-Battelle, LLC under Contract No. De-AC05-00OR22725 with the U.S. Department of Energy. The United States Government retains, and the publisher, by accepting the article for publication, acknowledges, that the United States Government retains a non-exclusive, paid-up, irrevocable, world-wide license to publish or reproduce the published form of this manuscript, or allow others to do so, for United States Government purposes. The Department of Energy will provide public access to these results of federally sponsored research in accordance with the DOE Public Access Plan (http://energy.gov/downloads/doe-public-access-plan).

Marlies Pirner was supported by the Alexander von Humboldt foundation and the German Science Foundation DFG (grant no. PI 1501/2-1).


\bibliographystyle{abbrv}
\bibliography{bibliography}

\end{document}